\documentclass[final,5p,twocolumn]{elsarticle}

\usepackage{amssymb}
 \usepackage{amsthm}

\usepackage[misc,geometry]{ifsym}
\usepackage{amsmath,amssymb}
\usepackage{subfig}
\usepackage{float}
\usepackage{hyperref}
\usepackage[numbers]{natbib}

\journal{International Journal of Non-Linear Mechanics}

\begin{document}

\begin{frontmatter}


 \title{Model Reduction to Spectral Submanifolds \\in Non-Smooth Dynamical Systems}



\author[]{Leonardo Bettini}
\author[]{Mattia Cenedese}
\author[]{George Haller\corref{cor1}}

\cortext[cor1]{Corresponding author.}
\ead{georgehaller@ethz.ch}
\address{Institute for Mechanical Systems, ETH Z\"urich,\\ Leonhardstrasse 21, 8092 Z\"urich, Switzerland\\}


\begin{abstract}
We develop a model reduction technique for non-smooth dynamical systems using spectral submanifolds. Specifically, we construct low-dimensional, sparse, nonlinear and non-smooth models on unions of slow and attracting spectral submanifolds (SSMs) for each smooth subregion of the phase space and then properly match them. We apply this methodology to both equation-driven and data-driven problems, with and without external forcing.
\end{abstract}




\end{frontmatter}


\section{Introduction}\label{intro}
The ever increasing need to reduce complex, high-dimensional dynamical systems to simple, low-dimensional models has yielded a number of different reduction techniques (see \citet{benner2015,rowley2017,ghadami2022,brunton2020,taira2017,touze2021} for recent reviews). Here, we focus on an extension of one of these methods, spectral submanifold (SSM) reduction, to non-smooth mechanical systems. 

Defined originally for smooth dynamical systems by \citet{NNM}, a primary SSM is the smoothest invariant manifold that is tangent to, and has the same dimension as, a spectral subspace of the linearized system at a steady state.  As such, SSMs mathematically formalise and extend the original idea of nonlinear normal modes (NNMs) introduced in seminal work by \citet{shaw_pierre_93,shaw_pierre_94} and \citet{shaw_pierre_99} (see \citet{mikhlin_2023} for a recent review).

The existence, uniqueness and persistence of SSMs in autonomous and non-autonomous systems have been proven whenever no resonance relationship holds between the linearized spectrum within the spectral subspace and outside that subspace (\citet{NNM,cabre2003,haro2006}). 
The primary SSM tangent to the spectral subspace spanned by the slowest linear modes attracts all nearby trajectories and hence its internal dynamics is an ideal, mathematically justified nonlinear reduced model. 

Recent work has revealed the existence of an additional, infinite family of secondary (or fractional) SSMs near the primary one in $C^\infty$ dynamical systems (\citet{haller2023_fractional}). 
Fractional SSMs are also tangent to the same spectral subspace as their primary counterpart, but they are only finitely many times differentiable. Accordingly, they can only be approximated via fractional-powered polynomial expansions.

When the equations of the system are known, its SSMs can be approximated via Taylor expansion at the stationary state using the SSMTool algorithm developed by \citet{shobit_2022,shobit_2023}. SSMTool can compute a low-dimensional reduced order model even for systems with hundreds of thousands degrees of freedom. The SSM-reduced models can in turn predict the response of the system to small harmonic forcing (\citet{shobit_2022,shobit_2017,Sten_2020}) and its bifurcations (\citet{mingwu_2022_II,mingwu_2022_I}).
In the absence of governing equations, SSMs and their reduced dynamics can also be approximated from data using the SSMLearn algorithm developed in \citet{mattia_2022}. Such data-driven SSM-reduced models have been found to capture the essentially nonlinear features of autonomous systems and also accurately predict nonlinear response under additional external forcing (\citet{mattia_2022,mattia_mechanical_systems_2021,mattia_toappear}). A more recent variant of the same algorithm, fastSSM, provides a simplified, faster SSM-reduction procedure with somewhat reduced accuracy (\citet{joar_2023}). Application of these methodologies has proven successful in a variety of examples, both numerical and experimental, ranging from beam oscillations and sloshing in water tank to structural vibrations and transition in shear flows (\citet{mattia_2022,mattia_mechanical_systems_2021,joar_2023,balint_2022}).\\

As most model reduction methods, SSM reduction also assumes that the full system to be reduced is sufficiently many times differentiable. In the absence of the required smoothness, such reduction methods either fail or apply only under modifications. For example, in \citet{cardin_2013,cardin_2015}, the authors showed the existence of a slow manifold in the context of singular perturbation, seeking an extension of invariant manifold-based model reduction for piecewise smooth systems. In particular, they study how the sliding mode present in such systems is affected by singular perturbation and prove that all hyperbolic equilibria and periodic orbits on the sliding region of the reduced problem persist. Motivated by a possible extension of the center manifold reduction of smooth systems to piecewise smooth systems, \citet{weiss_2012,weiss_2015} and \citet{kupper_2008} identify invariant cones as tools to reduce the dynamics and study bifurcation phenomena, when the equilibrium lies on the switching manifold between regions of smooth behavior. More precisely, invariant cone-like manifolds are found for nonlinear perturbations of linear piecewise smooth systems and they are constructed starting from a fixed point of the Poincaré map. 
A further notable contribution is the study of \citet{szalai_2019} on model reduction for non-densely defined piecewise smooth systems in Banach spaces. This highly technical approach uses singular perturbation techniques to develop meaningful reduced order models on low-dimensional invariant manifolds (including SSMs) across switching surfaces.

In this paper, we pursue a less technical, but more readily applicable objective. We consider finite-dimensional non-smooth dynamical systems with a switching surface that contains a fixed point. By extending the two smooth systems from either side of this surface, we conclude the existence of two smooth SSMs, which are separated from each other by a discontinuous jump. We then track full trajectories by appropriately switching between the two SSM-reduced dynamics. 

The structure of this  paper is as follows. We first describe our model reduction procedure for piecewise smooth system in section \ref{method}, recalling basic concepts of SSM theory. We then apply this procedure to a simple equation-driven example in section \ref{shaw_pierre}, in which we compute the SSMs analytically and compare several switching strategies among different SSMs. Finally, in section \ref{vonkarman_beam} we discuss a data-driven example of a von Kármán beam for which we carry out model reduction under different types of non-smoothness. 

\section{Method}\label{method}
For model reduction in a piecewise smooth system, we consider separately the subregions of the phase space in which the system is smooth and apply the results of primary SSM theory separately in those subregions. In particular, we smoothly extend each subsystem locally across its domain boundary and locate primary, smooth SSMs anchored at fixed points in that boundary for the extended system. Such an SSM will only act as an invariant manifold for the full system over its original subregion of smooth dynamics. These various subsets of SSMs form the skeleton of an attractive set for the full system, with pieces of this skeleton connected by trajectories sliding off from them and converging to other pieces of the skeleton (see Fig. \ref{manifolds_crossing_strategy} and video in the Supplementary Material). The reduced dynamics across different SSM pieces then needs to be connected appropriately by a reduced-order model, as we detail below.
\begin{figure*}
    \centering
    \subfloat[\label{PWS_manifold}]{
        \includegraphics[scale=0.38]{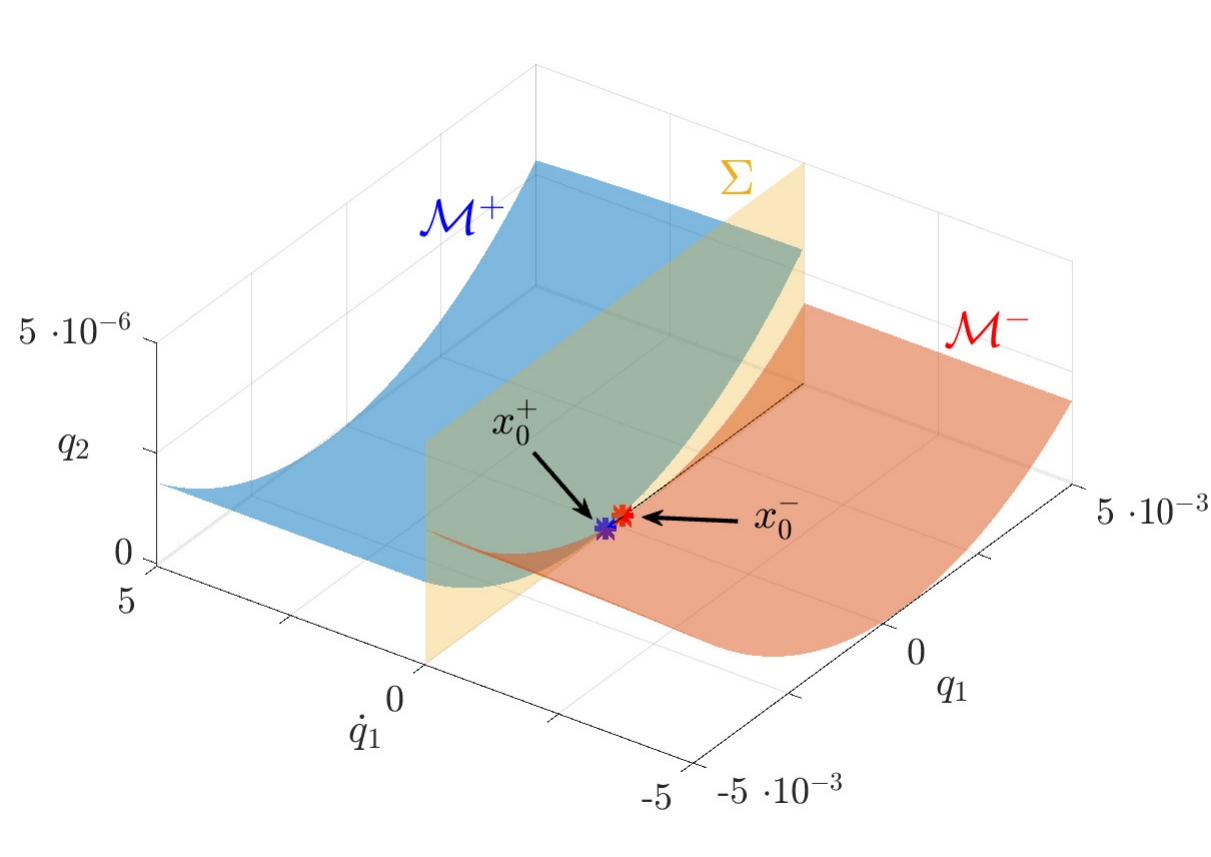}
    }
    \quad
    \subfloat[\label{PWS_crossing}]{
        \includegraphics[scale=0.38]{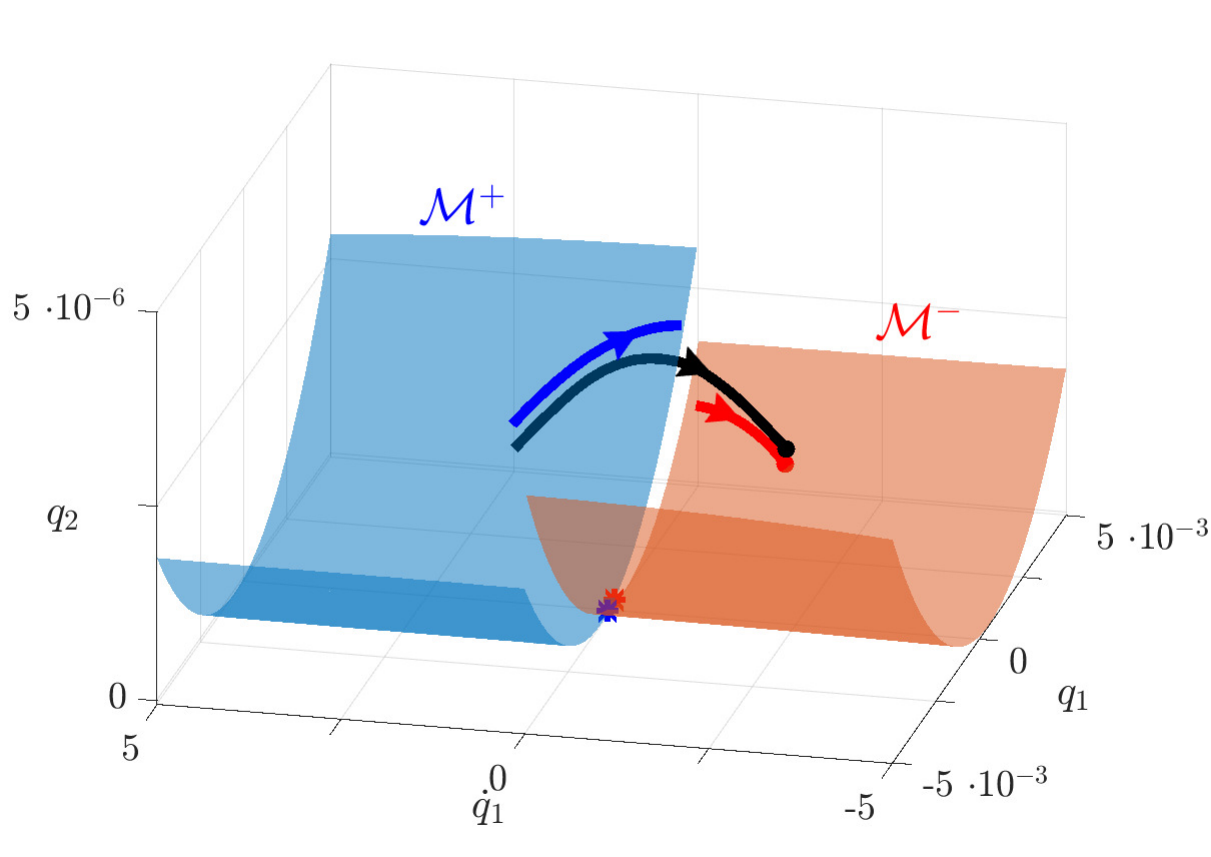}
    }
    \caption{Reduced-order modeling strategy for piecewise smooth dynamical systems using primary SSMs. The plots depict a specific example, described in section \ref{vonkarman_beam} and reported here for motivation. Subfigure \protect\subref{PWS_manifold} represents the two primary SSMs, $\mathcal{M}^+$ and $\mathcal{M}^-$, anchored at their respective equilibrium points, $x_0^+$ and $x_0^-$, and separated by the switching surface $\Sigma$. Subfigure \protect\subref{PWS_crossing} exemplifies the crossing between two subregions. The black curve is the trajectory of the full system, which is approximated by the reduced dynamics (blue and red curves) lying on $\mathcal{M}^+$ and $\mathcal{M}^-$. The full solution quickly synchronizes with the reduced one, even though the latter involves a physical discontinuity.} 
    \label{manifolds_crossing_strategy}
\end{figure*}

In the following, we first introduce some terminology from the theory of piecewise smooth dynamical systems, then give a more thorough description of our construction of a reduced-order model. 

\subsection{Piecewise smooth systems}\label{PWS_system_subsection} Let us consider the $n$-dimensional dynamical system
\begin{equation}\label{initial_system}
    \dot{\mathbf{x}} = \mathbf{f}(\mathbf{x}; \delta),\quad \mathbf{x} \in \mathbb{R}^n,\quad \delta \in \mathbb{R},
\end{equation}
where $\mathbf{x}$ is the state vector and $\mathbf{f}(\mathbf{x}; \delta)$ is a nonlinear and non-smooth right-hand side, depending smoothly on the parameter $\delta$ and the time $t$, but not necessarily on the phase-space variable $\mathbf{x}$. For simplicity, we assume $\mathbf{f}(\mathbf{x}; \delta)$ to be non-smooth across a single hypersurface $\Sigma$ of dimension $n-1$ that contains the origin $\mathbf{x} = \mathbf{0}$. Namely, the phase space is split into two regions separated by the hypersurface $\Sigma$, and $\mathbf{f}(\mathbf{x},\delta)$ is smooth within each of these regions.

The surface $\Sigma$ is usually called a switching surface, defined by a scalar-valued switching function $\sigma(\mathbf{x})$ as
\begin{equation}
    \Sigma = \left\{\mathbf{x} \in \mathbb{R}^n: \sigma(\mathbf{x})= 0 \right\}
\end{equation}
The phase space is then partitioned as $\mathbb{R}^n = \Sigma^+ \cup \Sigma \cup \Sigma^-$, where
\begin{equation}
    \begin{array}{l}
         \Sigma^- = \left\{\mathbf{x} \in \mathbb{R}^n : \sigma(\mathbf{x}) < 0 \right\},\\[\medskipamount]
         \Sigma^+ = \left\{\mathbf{x} \in \mathbb{R}^n : \sigma(\mathbf{x}) > 0 \right\}.
    \end{array}
\end{equation}
The original piecewise smooth system \eqref{initial_system} can now be written as
\begin{equation}\label{pws_system}
    \dot{\mathbf{x}} = \mathbf{f}(\mathbf{x};\delta) =  
    \begin{cases}
    \mathbf{f}^+(\mathbf{x};\delta), & \mathbf{x} \in \Sigma^+,\\
    \mathbf{f}^-(\mathbf{x};\delta), & \mathbf{x} \in \Sigma^-,
    \end{cases}
\end{equation}
where we assume that $\mathbf{f}^\pm(\mathbf{x};\delta)$ both extend to smooth functions of $\mathbf{x}$ in an open neighborhood of $\Sigma$. We also assume that
\begin{equation}
    \begin{array}{c}
          \mathbf{f}(\mathbf{x};0)\equiv \mathbf{f}^\pm(\mathbf{x};0), \quad \mathbf{f}(\mathbf{0};0) = \mathbf{0},\\[\medskipamount]
           0 \notin \text{Re}[\text{Spect}(\text{D}_\mathbf{x} \mathbf{f}(\mathbf{0};0)].
    \end{array}
\end{equation}

In other words, for $\delta = 0$, the discontinuity of system \eqref{initial_system} disappears and $\mathbf{x} = \mathbf{0}$ is a hyperbolic fixed point of \eqref{initial_system}, contained in $\Sigma$ for this value of $\delta$.  Note that by the theory developed by \citet{Filippov}, the dynamics within $\Sigma$ can be approximated by constructing a proper inclusion (see \ref{supp_PWS}). 

We assume that trajectories intersecting $\Sigma$ exhibit two possible behaviors:
 
\begin{figure}[H]
    \centering
    \subfloat[\label{crossing_1}]{
        \includegraphics[scale=0.22]{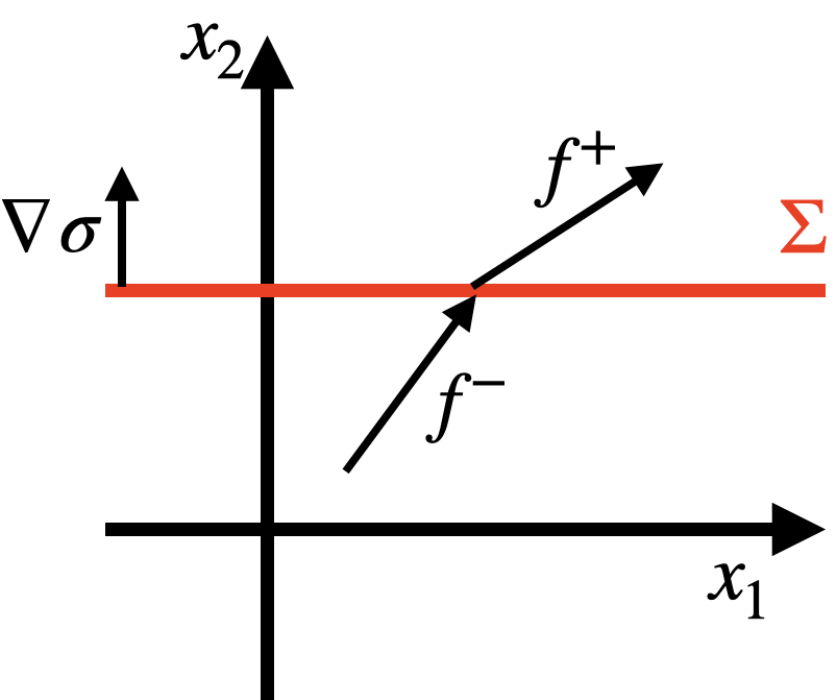}
    }
    \quad
    \subfloat[\label{crossing_2}]{
        \includegraphics[scale=0.22]{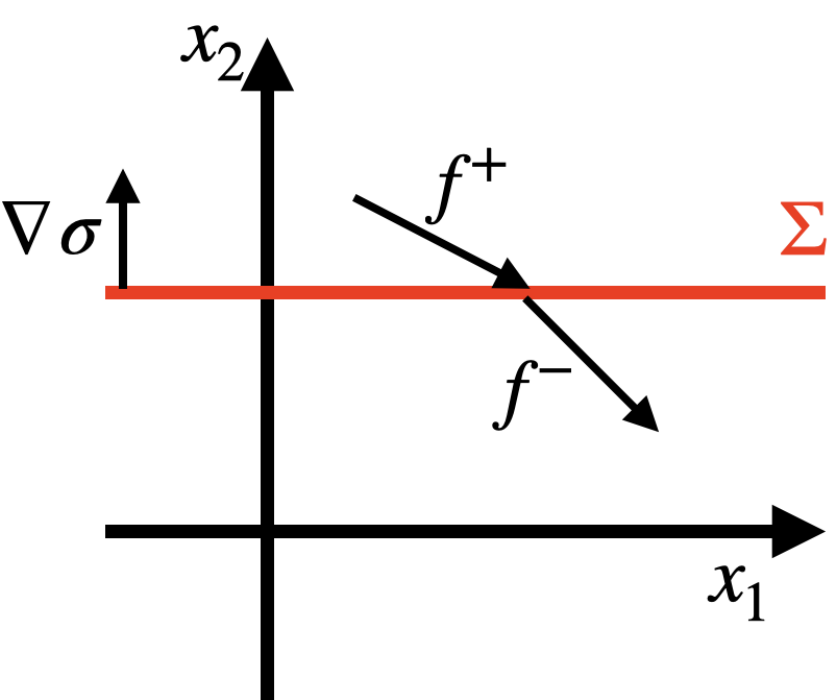}
    }\quad
    \subfloat[\label{attractive_sliding}]{
        \includegraphics[scale=0.22]{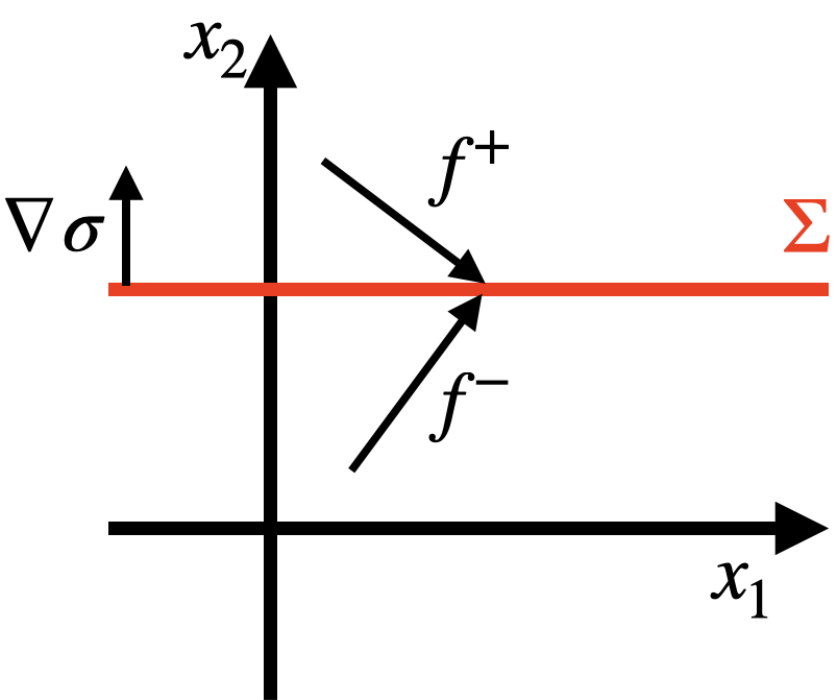}
    }
    \quad
    \subfloat[\label{repulsive_sliding}]{
        \includegraphics[scale=0.22]{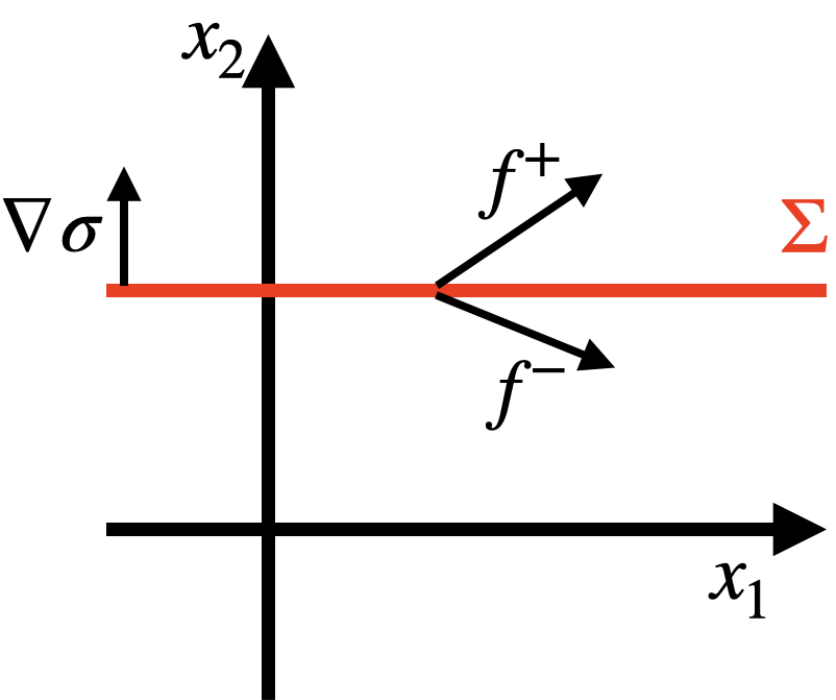}
    }
    \caption{The orbits can cross either from $\Sigma^-$ to $\Sigma^+$ \protect\subref{crossing_1} or from $\Sigma^+$ to $\Sigma^-$ \protect\subref{crossing_2}. Alternatively, they can exhibit an attracting 
  \protect\subref{attractive_sliding} or a repelling sliding mode \protect\subref{repulsive_sliding}. }
    \label{crossing_sliding_figure}
    \end{figure}
\begin{itemize}
    \item \textbf{Crossing}: the orbit crosses trough $\Sigma$ if
    \begin{equation}\label{crossing_condition}
        (\nabla \sigma \cdot \mathbf{f}^+) (\nabla \sigma \cdot \mathbf{f}^-)>0 \text{ on } \Sigma.
    \end{equation}
    \item \textbf{Sliding}: the orbit slides along the switching surface if
    \begin{equation}\label{sliding_condition}
        (\nabla \sigma \cdot \mathbf{f}^+) (\nabla \sigma \cdot \mathbf{f}^-)<0 \text{ on } \Sigma.
    \end{equation}

    In the case of sliding, the vector field is tangent to the switching surface and can be expressed as a convex combination 
    of $\mathbf{f}^+$ and $\mathbf{f}^-$ along $\Sigma$. The switching surface $\Sigma$ then either attracts or repels neighbouring trajectories, depending on the direction of the vector fields in $\Sigma^+$ and $\Sigma^-$, as seen in Fig. \ref{crossing_sliding_figure}. In the former case, the solution generally evolves along the switching surface. In the latter case, a solution still exists, but it is not unique in forward time.
    
\end{itemize}

\subsection{Model reduction strategy}\label{model_reduction_strategy}  We now outline a procedure which allows us to extend SSM-based model reduction to piecewise smooth systems. Our approach utilizes the primary SSMs constructed over slow modes existing on both sides of $\Sigma$ by the theory of SSMs for smooth systems reviewed in \ref{supp_SSM}. 

By our assumptions in section \ref{PWS_system_subsection}, for $\delta = 0$, system \eqref{initial_system} has a hyperbolic fixed point at $\mathbf{x} = \mathbf{0}$. Any nonresonant spectral subspace $E_0$ of the linearized system at $\mathbf{x} = \mathbf{0}$ will then admit a unique, primary SSM, $\mathcal{M}_0$, under the addition to the nonlinear terms of $\mathbf{f}(\mathbf{x};0)$. For $\delta \neq 0$ and small, hyperbolic fixed points will continue to exist $\mathcal{O}(\delta)$-close to $\mathbf{x} = \mathbf{0}$ for the two ODEs $\dot{\mathbf{x}} = \mathbf{f}^\pm(\mathbf{x};\delta)$ by the implicit function theorem. These fixed points, $\mathbf{x}_0^\pm(\delta)$, will have their own SSMs, $\mathcal{M}^\pm_\delta$, that are smooth, $\mathcal{O}(\delta)$ perturbations of $\mathcal{M}_0$ by the smooth dependence of SSMs on parameters (\citet{NNM}).
As a consequence, the set $\mathcal{M}_\delta = \mathcal{M}_\delta^+ \cup \mathcal{M}_\delta^-$ is a piecewise smooth attracting set for system \eqref{initial_system}.\\

We now recall the construction of a general SSM and its reduced dynamics from observable data (see \citet{mattia_2022,joar_2023} for more detail). We seek the parametrization of the manifold $\mathcal{M}$ as a graph over its \textit{a priori} unknown tangent space $T\mathcal{M}$, which is the image of the spectral subspace $E$ in the space of observables $\mathbf{y}\in\mathbb{R}^n$. In particular, we approximate the manifold via a multi-variate Taylor expansion
\begin{equation}\label{parametrization_formula}
    \mathbf{y}(\boldsymbol{\xi}) = \mathbf{M}\boldsymbol{\xi}^{1:m} = \mathbf{V}\boldsymbol{\xi} + \mathbf{M}_{2:m} \boldsymbol{\xi}^{2:m},
\end{equation}
where $\mathbf{M} = [\mathbf{V}, \mathbf{M}_2, ..., \mathbf{M}_m]$, with $\mathbf{M}_i \in \mathbb{R}^{n \times d_i}$ and $d_i$ refers to the number of $d$-variate monomials at order $i$. We refer to  $\boldsymbol{\xi} = \mathbf{V}^T \mathbf{y}$ as reduced coordinates, where the columns of the matrix $\mathbf{V} \in \mathbb{R}^{n \times d}$ are orthonormal vectors that span $T\mathcal{M}$. The notation $(\cdot)^{l:r}$ refers to the vector containing all monomials composed of the entries of the vector $(\cdot)$, with monomials ranging from $l$ to $r$. For instance, if $\boldsymbol{\xi} = [\xi_1 \xi_2]^T$, then
\begin{equation*}
    \boldsymbol{\xi}^{2:3} = [\xi_1^2, \,\xi_1\xi_2,\, \xi_2^2,\, \xi_1^3, \,\xi_1^2\xi_2,\, \xi_1\xi_2^2,\, \xi_2^3]^\mathrm{T}.
\end{equation*}

Learning $\mathcal{M}$ from a set of training data $\{\mathbf{y}_j\}$ in the observable space means finding the optimal matrix $\mathbf{M}^\ast$, such that
\begin{equation}
   		 \left(\mathbf{V^\ast},\mathbf{M}^\ast \right)= \underset{\mathbf{V},\mathbf{M}_{2:m}}{\arg\min} \sum_j \Big\| \mathbf{y}_j - \mathbf{V}\mathbf{V}^T\mathbf{y}_j - \mathbf{M}_{2:m}(\mathbf{V}^T\mathbf{y}_j)^{2:m}\Big\|,
\end{equation}
subject to the constraints
\begin{equation}
    \mathbf{V}\mathbf{V}^T = \mathbb{I}, \quad \mathbf{V}\mathbf{M}_{2:m} = \mathbf{0}.
\end{equation}

The dynamics on the SSM in the reduced coordinates $\boldsymbol{\xi}$ can then be approximated as
\begin{equation}\label{reduced_dynamics}
    \dot{\boldsymbol{\xi}} = \mathbf{R}\boldsymbol{\xi}^{1:r},
\end{equation}
where the elements of $\mathbf{R} \in \mathbb{R}^{d:d_{1:r}}$ are found by solving the minimization problem
\begin{equation}
    \mathbf{R}^\ast = \underset{\mathbf{R}}{\arg\min} \Big\| \dot{\boldsymbol{\xi}}_j - \mathbf{R}\boldsymbol{\xi}_j^{1:r} \Big\|.
\end{equation}

In the present work, we assume to know \textit{a priori} the model of the system that generates the data $\mathbf{y}_j$, as well as the linear parts of the parametrization \eqref{parametrization_formula} and reduced dynamics \eqref{reduced_dynamics}, as in  \citet{mattia_toappear}. More specifically, in the piecewise smooth context, we suppose that the governing equations, the domains $\Sigma^\pm$, and the switching function $\sigma(\mathbf{x})$ defining $\Sigma$ are known. 

With these ingredients, we will use the smoothly extended versions of the right-hand sides $\mathbf{f}^\pm(\mathbf{x};\delta)$ to generate training data $\{\mathbf{y}_j^\pm\}$ from both ODEs $\dot{\mathbf{x}} = \mathbf{f}^\pm(\mathbf{x};\delta)$ in a neighborhood of the switching surface. We then construct the SSMs, $\mathcal{M}^\pm_\delta$, separately, but only keep their subsets falling in the domains $\Sigma^\pm$, respectively. Our reduced-order will then switch between the reduced dynamics of $\mathcal{M}_\delta^\pm$ based on appropriately reduced switching conditions that we will discuss in our upcoming examples. 

\section{Example 1: Shaw-Pierre model with friction}\label{shaw_pierre}
\begin{figure}[H]
    \centering
    \subfloat[\label{mechanical_system}]{
        \includegraphics[scale=0.6]{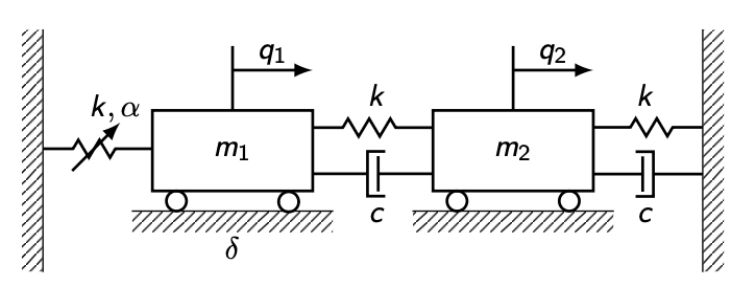}
    }
    \quad
    \subfloat[\label{friction_law}]{
        \includegraphics[scale=0.5]{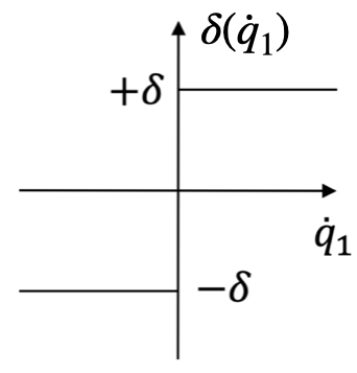}
    }
    \caption{Two-degree-of-freedom modified Shaw-Pierre mechanical system with friction. \protect\subref{mechanical_system} System geometry. \protect\subref{friction_law} Coulomb friction law added to the first mass.}
    \label{SP_figure}
\end{figure}

We add friction to a modified version of the mechanical system studied by \citet{shaw_pierre_93}, as introduced in \citet{NNM}. The resulting system is sketched in Fig. \ref{SP_figure}.

Dry friction is modelled via the classical Coulomb's law (\citet{leine}), where the static friction coefficient $\delta$ (active when $\dot{q}_1 = 0$) is equal to the dynamic one (valid when $\dot{q}_1 \neq 0$). 
We now follow the procedure outlined in section \ref{method} for the computation of a reduced order model for this piecewise smooth system.
\paragraph{Non-smooth formulation} 
The equations of motion of the system in Fig. \ref{SP_figure} are
\begin{equation}\label{eq:physical}
\begin{aligned}
&
    \begin{pmatrix}
    m_1 & 0\\
    0 & m_2
    \end{pmatrix}
    \begin{pmatrix}
    \ddot q_1 \\
    \ddot q_2
    \end{pmatrix}
    +
    \begin{pmatrix}
    c & -c\\
    -c & 2c
    \end{pmatrix}
    \begin{pmatrix}
    \dot q_1 \\
    \dot q_2
    \end{pmatrix}\\&\quad
    +
    \begin{pmatrix}
    2k & -k\\
    -k & 2k
    \end{pmatrix}
    \begin{pmatrix}
    q_1 \\
    q_2
    \end{pmatrix}
    +
    \begin{pmatrix}
    -\alpha q_1^3 \\
    0
    \end{pmatrix}
    = \begin{pmatrix}-\delta({\dot{q}_1})m_1 g\\0\end{pmatrix},
\end{aligned}
\end{equation}
where $g$ denotes the constant of gravity. In the following, $\delta = \delta(\dot{q}_1) g$.
With the notation $x_1 = q_1$, $x_2 = \dot{q}_1$, $x_3 = q_2$ and $x_4 = \dot{q}_2$, we can rewrite eq.\eqref{eq:physical} as a first-order system of ODEs

\begin{equation}\label{eq:x_coordinates}
    \dot {\mathbf{x}} =
    \mathbf{A}_0\mathbf{x} + \mathbf{f_{nl}} \mp \,\mathbf{f_{\delta}},
\end{equation}
with
\begin{align*}
    &
    \mathbf{A}_0 = \begin{pmatrix}
    0 & 1 & 0 & 0 \\
    -\frac{2k}{m_1} & -\frac{c}{m_1} & \frac{k}{m_1} & \frac{c}{m_1} \\
    0 & 0 & 0 & 1 \\
    \frac{k}{m_2} & \frac{c}{m_2} & -\frac{2k}{m_2} & -\frac{2c}{m_2}
    \end{pmatrix},&
    \\[\medskipamount]&
    \mathbf{f_{nl}} = \begin{pmatrix}
    0,
    -\frac{\alpha}{m_1} x_1^3,
    0,0
    \end{pmatrix}^\mathrm{T},
    \quad
    \mathbf{f_{\delta}} = \begin{pmatrix}
    0, \delta ,0 ,0
    \end{pmatrix}^\mathrm{T}.
\end{align*}

The sign of the velocity $x_2$ of the first mass decides whether the friction force is positive or negative, hence this velocity defines the switching function
\begin{equation}
    \sigma(\mathbf{x}) = x_2 = \dot{q}_1.
\end{equation}
Therefore, we have the splitting $\mathbb{R}^4 = \Sigma^+ \cup \Sigma \cup \Sigma^-$, where
\begin{equation}\label{subregions}
    \begin{array}{l}
        \Sigma ^\pm = \left\{\mathbf{x} \in \mathbb{R}^4 : \text{sign}(\sigma(\mathbf{x})) = \pm 1 \right\}.
    \end{array}
\end{equation}

The piecewise smooth system \eqref{initial_system} is here specifically defined as 
\begin{equation}\label{SP_equation_motion_pm}
\begin{aligned}
    \dot{\mathbf{x}} &= \mathbf{f}^\pm(\mathbf{x};\delta) \\&= 
    \begin{pmatrix}
        \displaystyle
        x_2\\\displaystyle
        -\frac{2k}{m_1} x_1 - \frac{c}{m_1} x_2 + \frac{k}{m_1} x_3 + \frac {c}{m_1} x_4 -\frac{\alpha}{m_1} x_1^3 \mp \delta \\
        x_4\\\displaystyle
        \frac{k}{m_2} x_1 + \frac{c}{m_2} x_2 - \frac{2k}{m_2} x_3 - \frac{2c}{m_2} x_4
    \end{pmatrix}.
\end{aligned}
\end{equation}
Given that $\nabla \sigma = (0,1,0,0)^\mathrm{T}$, a crossing of $\Sigma$ at a point $\mathbf{x}\in\Sigma$ takes place if $\nabla \sigma \cdot \mathbf{f}^+$ and $ \nabla \sigma \cdot \mathbf{f}^->0$ have the same nonzero sign at $\mathbf{x}$, i.e.,
\begin{equation}\label{crossing_condition_SP}
\begin{aligned}
&
    \begin{cases}
    \displaystyle
        -\frac{2k}{m_1}x_1  + \frac {k}{ m_1} x_3 + \frac {c}{m_1} x_4 - \frac{\alpha}{m_1} x_1^3 > \delta, \vspace{0.3cm} \\
        \displaystyle
        -\frac{2k}{m_1}x_1  + \frac {k}{ m_1} x_3 + \frac {c}{m_1} x_4 - \frac{\alpha}{m_1} x_1^3 > -\delta,
    \end{cases} \text{or} \\[0.3cm]&
     \begin{cases}
     \displaystyle
        -\frac{2k}{m_1}x_1  + \frac {k}{ m_1} x_3 + \frac {c}{m_1} x_4 - \frac{\alpha}{m_1} x_1^3 < \delta, \vspace{0.3cm} \\
        \displaystyle
        -\frac{2k}{m_1}x_1  + \frac {k}{ m_1} x_3 + \frac {c}{m_1} x_4 - \frac{\alpha}{m_1} x_1^3 < -\delta.
    \end{cases}
\end{aligned}
\end{equation}
Attracting sliding motion along $\Sigma$ ensues when 
\begin{equation}\label{sticking_condition}
    \Big|-\frac{2k}{m_1}x_1  + \frac {k}{ m_1} x_3 + \frac {c}{m_1} x_4 -  \frac{\alpha}{m_1} x_1^3 \Big| <\delta .
\end{equation}
Inside the switching surface $\Sigma$, the first mass is in the state of sticking and the second mass acts as a linear harmonic oscillator, i.e., we have 
\begin{equation}
    \dot{\mathbf{x}} = 
    \begin{pmatrix}
0\\0\\x_4\\\frac{k}{m_2} x_1 - \frac{2k}{m_2}x_3 - \frac{2c}{m_2}x_4
    \end{pmatrix}, \quad \mathbf{x} \in \Sigma.
\end{equation}
The conditions for repulsive sliding mode are
\begin{equation}
    \begin{cases}
    \displaystyle
        -\frac{2k}{m_1}x_1  + \frac {k}{ m_1} x_3 + \frac {c}{m_1} x_4 - \frac{\alpha}{m_1} x_1^3 > \delta \vspace{0.3cm}, \\
        \displaystyle
            -\frac{2k}{m_1}x_1  + \frac {k}{ m_1} x_3 + \frac {c}{m_1} x_4 - \frac{\alpha}{m_1} x_1^3 < -\delta ,
    \end{cases}
\end{equation}
which cannot be satisfied for any $\mathbf{x} \in \Sigma$. As a consequence, the Shaw-Pierre system with friction exhibits either crossing or attracting sliding (or sticking) behavior. 

We will simulate the behavior of the system with the parameter values used by \citet{shaw_pierre_93} and \citet{NNM},
\begin{equation*}
    m = 1, \quad c = 0.3, \quad k = 1, \quad \alpha = 0.5,
\end{equation*}
and will consider a range of $\delta$ values in our analysis.

\paragraph{Analysis of the linearized system} The two different dynamical systems defined in \eqref{SP_equation_motion_pm} have their own fixed points and smooth SSMs anchored at them.
The fixed points are defined by
\begin{equation}
\begin{aligned}
   & \mathbf{x}_0^\pm = \begin{pmatrix}
    q_0^\pm \\
    0 \\
    \frac{1}{2}q_0^\pm \\
    0
    \end{pmatrix},\\ &q_0^\pm =\pm  \sqrt[3]{-\delta+ \sqrt{\delta ^2+1}} \mp \sqrt[3]{\delta + \sqrt{\delta^2 + 1}}.
\end{aligned}
\end{equation}
We shift coordinates so that the origin in the two cases coincides with $\mathbf{x}_0^\pm$, respectively:
\begin{equation}\label{shifted_coordinates}
    \mathbf \xi^\pm = \mathbf x-\mathbf x_0^\pm = \begin{pmatrix}
    q_1 - q_0^\pm\\
    \dot q_1 \\
    q_2 - \frac{1}{2}q_0^\pm \\ 
    \dot q_2
    \end{pmatrix}.
\end{equation}
In these coordinates, system \eqref{SP_equation_motion_pm} becomes 
\begin{equation}\label{eq:xi_coordinates}
    \dot {\mathbf{\xi}}^\pm =
    \mathbf{\tilde{A}}_0 \mathbf{\xi}^{\pm} + \mathbf{f_{nl}^{II}}(\mathbf{\xi}^\pm) + \mathbf{f_{nl}^{III}}(\mathbf{\xi}^\pm)  + \mathbf{f}_0,
\end{equation}
with
\begin{equation}\label{shifted_terms}
\begin{aligned}
    &\mathbf{\tilde{A}}_0 = \begin{pmatrix}
    0 & 1 & 0 & 0 \\
    -2 - \frac {3}{2} {(q_0^\pm)}^2 & -0.3 & 1 & 0.3\\
    0 & 0 & 0 & 1 \\
    1&0.3&-2&-0.6
    \end{pmatrix},\\[\medskipamount]&
    \mathbf{f_{nl}^{II}}(\mathbf{\xi}^\pm) = \begin{pmatrix}
    0,
   - \frac {3}{2} q_0^\pm \left(\xi_1^\pm\right)^2 ,0,0
    \end{pmatrix}^\mathrm{T},\\[\medskipamount]&
    \mathbf{f_{nl}^{III}}(\mathbf{\xi}^\pm) = \begin{pmatrix}
    0,
    -\frac {1}{2} \left(\xi_1^\pm\right)^3,
    0,0
    \end{pmatrix}^\mathrm{T},\\[\medskipamount]&
    \mathbf{f}_0 = 
    \begin{pmatrix}
    0,-\frac{3}{2}q_0^\pm - \frac{1}{2}(q_0^\pm)^3 \mp\delta ,0,0
    \end{pmatrix}^\mathrm{T}.
\end{aligned}
\end{equation}
The matrix $\mathbf{\Tilde{A}}_0$ is unaffected by the sign choice in $q_0^\pm$, therefore both cases have the same spectral properties and hence the same spectral subspaces.\\
The eigenvalues of $\mathbf{\Tilde{A}}_0$ for $\delta = 10^{-1}$ are
\begin{equation*}
    \begin{array}{l}
          \lambda_{1,2} = -0.0741 \pm i\, 1.0027,\\
          \lambda_{3,4} = -0.3759 \pm i\,1.6812 ,
    \end{array}
\end{equation*}
whose eigenvectors give rise to two two-dimensional real invariant subspaces, $E_1$ and $E_2$.

\subsection{Computation of SSMs.} According to the definition given in \ref{supp_SSM}, the relative spectral quotients of $E_1$ and $E_2$ are
\begin{equation}\label{eq:spectral_quotients}
    \begin{array}{l}
          \sigma(E_1) = \text{Int}\left[\frac{\text{Re}\,\lambda_{3,4}}{\text{Re}\,\lambda_{1,2}}\right] = 5 .
    \end{array}
\end{equation}
Changing the friction coefficient $\delta$ has a mild effect on the eigenvalues and hence the spectral quotients will not change for the range of friction coefficients studied here. Considerations about the existence and uniqueness of the SSMs in the positive and negative cases are exactly the same, as they share the same linearized dynamics, even though they are anchored at different points. The nonlinear contribution of the spring $\mathbf{f_{nl}}$ is an autonomous term that is analytic on the whole phase space. 

Based on these facts, the results in \citet{NNM} guarantee the existence and uniqueness of the slow two-dimensional SSM,  $\mathcal{M}_1(0)$, because the required nonresonance conditions among the eigenvalues are satisfied up to order $\sigma(E_1)$ (see \citet{NNM} for details). Therefore, we can state that the analytic SSM $\mathcal{M}_1(0)$ exists and it is unique among all $\mathcal{C}^6$ invariant manifolds tangent to $E_1$ at the origin. We introduce coordinates aligned with $E_1$ and $E_2$ by letting 
\begin{equation}\label{parametrization_matrix}
    \boldsymbol{\xi}^\pm = \mathbf{V}\boldsymbol{\eta}^\pm,
\end{equation}
where
\begin{equation*}
    \boldsymbol{\eta}^\pm = (\mathbf{y}^\pm, \mathbf{z}^\pm) \in E_1 \times E_2,
\end{equation*}
and $\mathbf{V}$ is the matrix whose columns are the eigenvectors of $\mathbf{\tilde{A}_0}$. The reduced coordinates $\mathbf{y}^\pm$ act as the master variables over which we seek the slow SSM as a graph. \\
We report the analytical computation of the SSM in  \ref{supp_SP}. 

\subsection{Combination of SSMs and the switching surface.} 
\begin{figure}[H]
\centering
\includegraphics[scale=0.4]{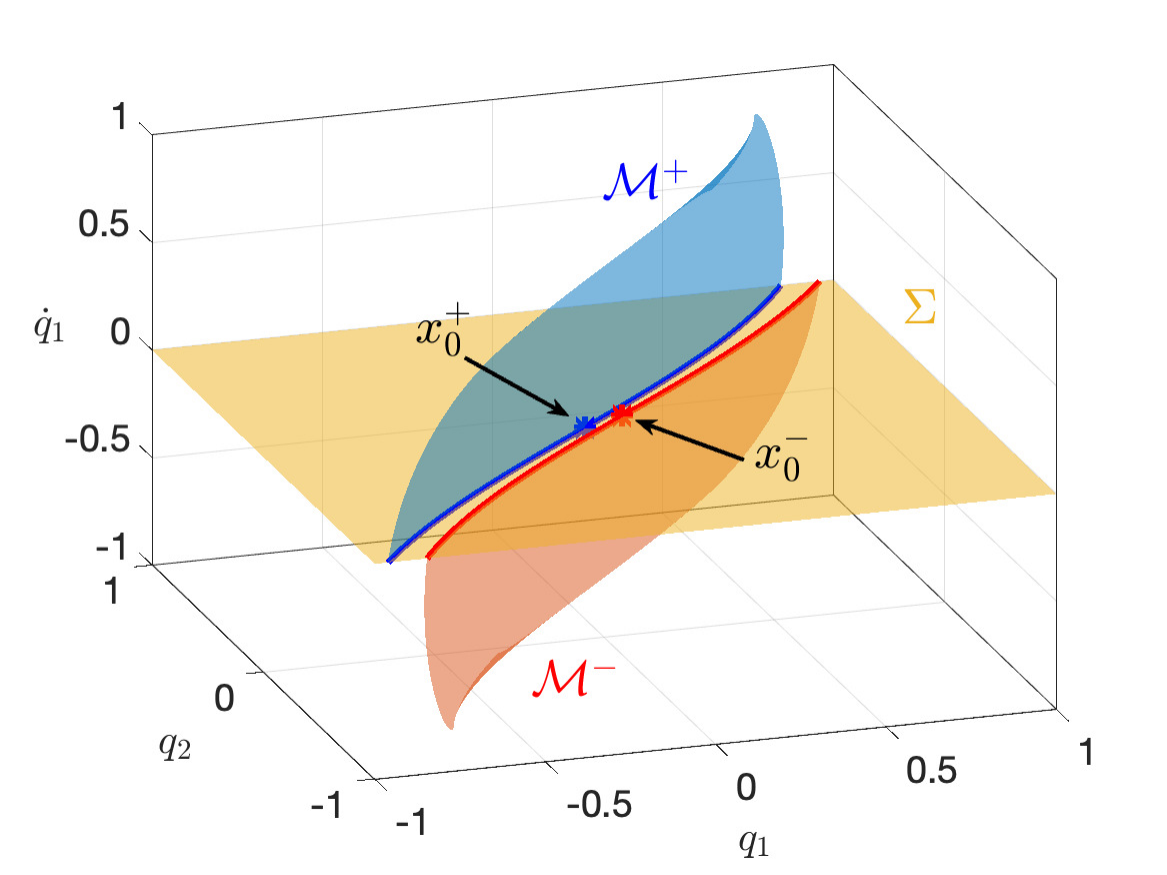}
\caption{The two primary SSMs $\mathcal{M}^+$ and $\mathcal{M}^-$, separated by the switching surface $\Sigma$, for $\delta  = 10^{-1}$.}
\label{SSMs}
\end{figure}

We now restrict the SSMs constructed at $\mathbf{x}_0^\pm$ to their region of validity imposed by condition \eqref{subregions}. The resulting restricted SSMs are separated by the switching surface $\Sigma$, as seen in Fig. \ref{SSMs}.

\begin{figure}[H]
    \centering
    \subfloat[\label{ic_strategies}]{
        \includegraphics[scale=0.6]{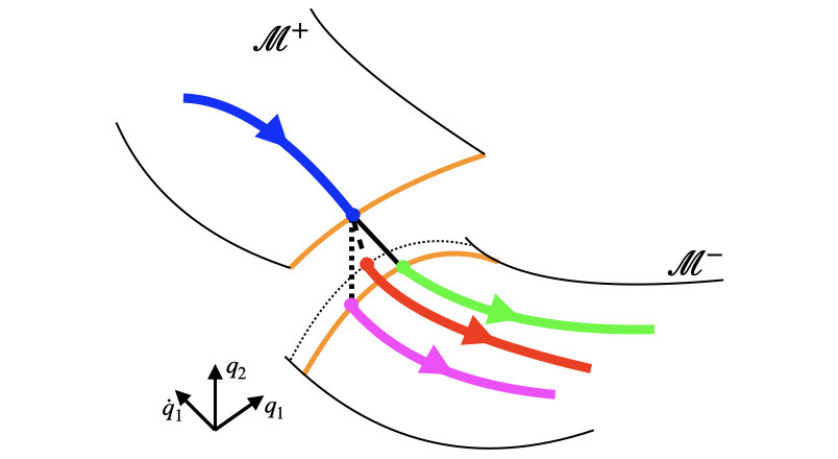}
    }
    \quad
    \subfloat[\label{ic_strategies2}]{
        \includegraphics[scale=0.46]{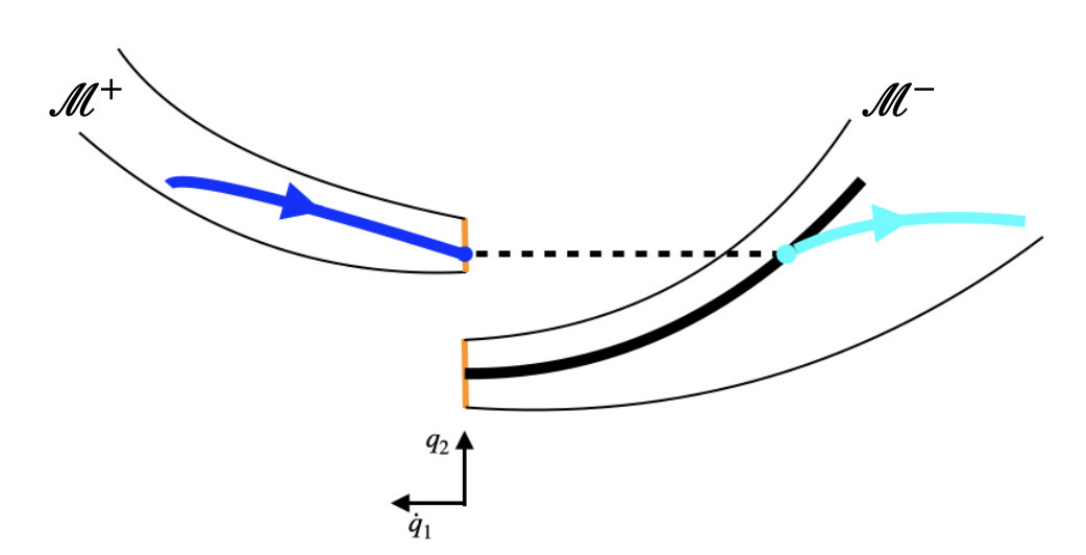}
    }
    \caption{Different strategies for the computation of the new initial condition of the reduced dynamics across the switching surface $\Sigma$. 
    The orange lines represent the intersection of the SSMs with the $\Sigma$. When a trajectory on one SSM (blue line) hits the switching surface, the projection of the intersection point onto the new SSM according to eq. \eqref{new_ic_equation} (red dot) represents a new initial condition. Alternatively, we can enforce that the new initial condition lies on the switching surface, minimizing the different of all physical variables across the jump (green dot) or ensuring the continuity of $q_1$ (purple dot). If we require that both $q_1$ and $q_2$ are continuous (light blue dot), then the new initial condition lands far from the switching surface.
    The black line in \protect\subref{ic_strategies2} depicts the set of points on the manifold $\mathcal{M}^-$ with the same value of $q_1$ as the final point of the blue trajectory on $\mathcal{M}^+$.}
    \label{ic_strategies_all}
\end{figure}

As the two pieces of primary SSMs form an attracting set for the full system, the dynamics restricted to them provides a reduced-order model for nearby initial conditions. We now describe how we connect the dynamics across the different SSM pieces when a trajectory hits the switching surface and satisfies the crossing condition \eqref{crossing_condition_SP}. In this scenario, coming from one SSM, we need a new initial condition lying on the other one. 
By Fig. \ref{ic_strategies_all}, the intersection between the incoming trajectory and the switching surface (blue dot in Fig. \ref{ic_strategies}) is associated to a set of both reduced $\mathbf{y}^\pm = \left(y_1^\pm, y_2^\pm\right)^\mathrm{T}$ and observable $\mathbf{x}= \left(q_1, \dot{q}_1, q_2, \dot{q}_2\right)^\mathrm{T}$ coordinates. We then exploit the relationship between the two sets of reduced coordinates for the two SSMs to find the new initial condition (red dot in Fig. \ref{ic_strategies}) as
\begin{equation}
    \begin{pmatrix}
        y_1^-\\y_2^-
    \end{pmatrix} = 
    \begin{pmatrix}
        y_1^+\\y_2^+
    \end{pmatrix} + \mathbf{V}^{-1}_{12}(\mathbf{x}_0^+ - \mathbf{x}_0^-),
    \label{new_ic_equation}
\end{equation}
where the rectangular matrix $\mathbf{V}^{-1}_{12}$ contains the first two rows of $\mathbf{V}^{-1}$.
Coming from $\mathcal{M}^+$, the new initial condition $(y_1^-, y_2^-)$ on $\mathcal{M}^-$ according to eq. \eqref{new_ic_equation} does not lie in the intersection between the $\mathcal{M}^-$ and the switching surface. As a consequence, all physical variables experience a discontinuity while crossing from one SSM to the other. Discontinuities in the solution are unavoidable, but one can investigate further strategies in order to enforce physical consistency for specific variables, as seen in Fig. \ref{ic_strategies_all}. 

\begin{figure*}
    \centering
    \subfloat[\label{edges_full_scaled}]{
        \includegraphics[scale=0.2]{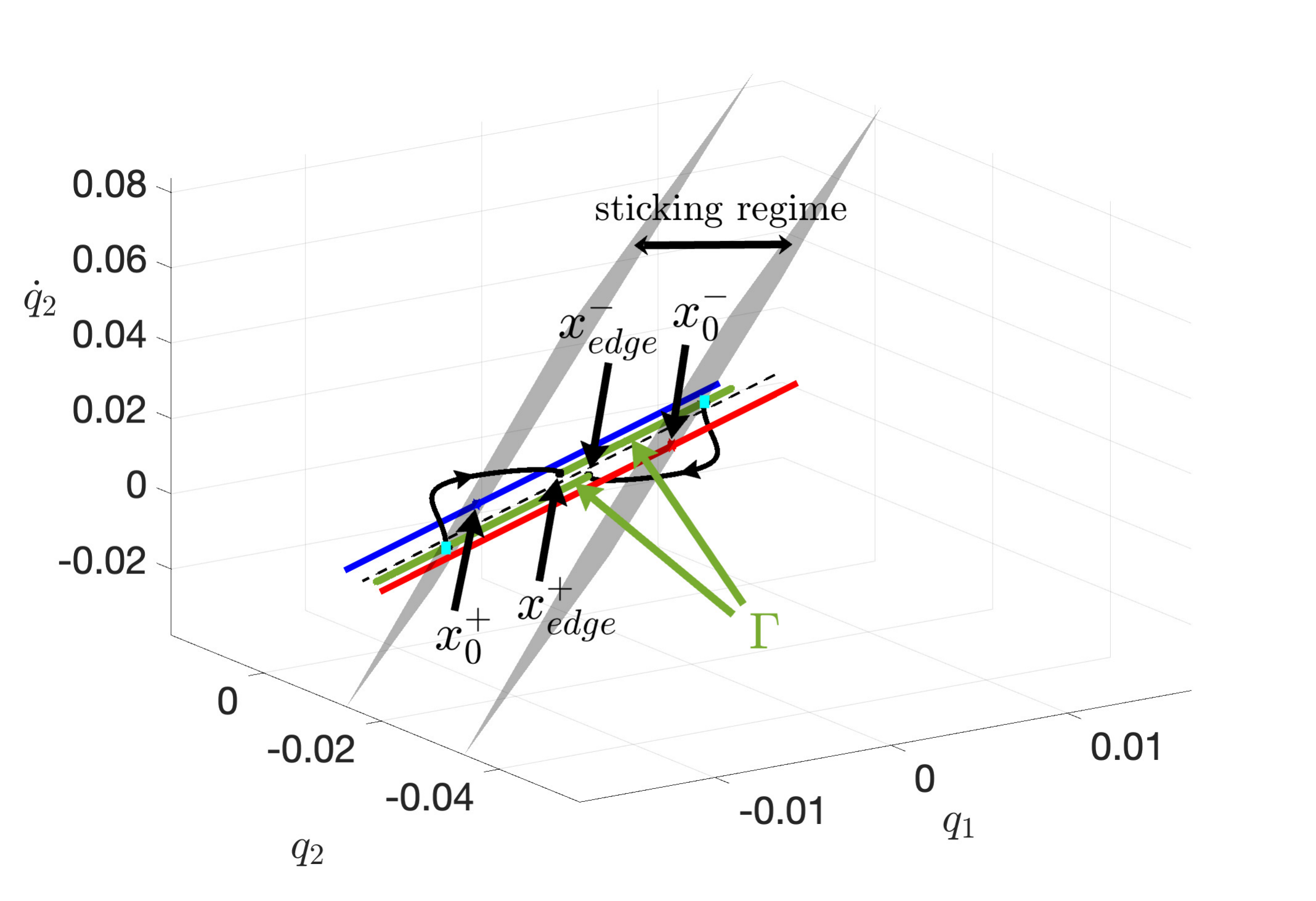}
    }
    \quad
    \subfloat[\label{edges_reduced_scaled}]{
        \includegraphics[scale=0.2]{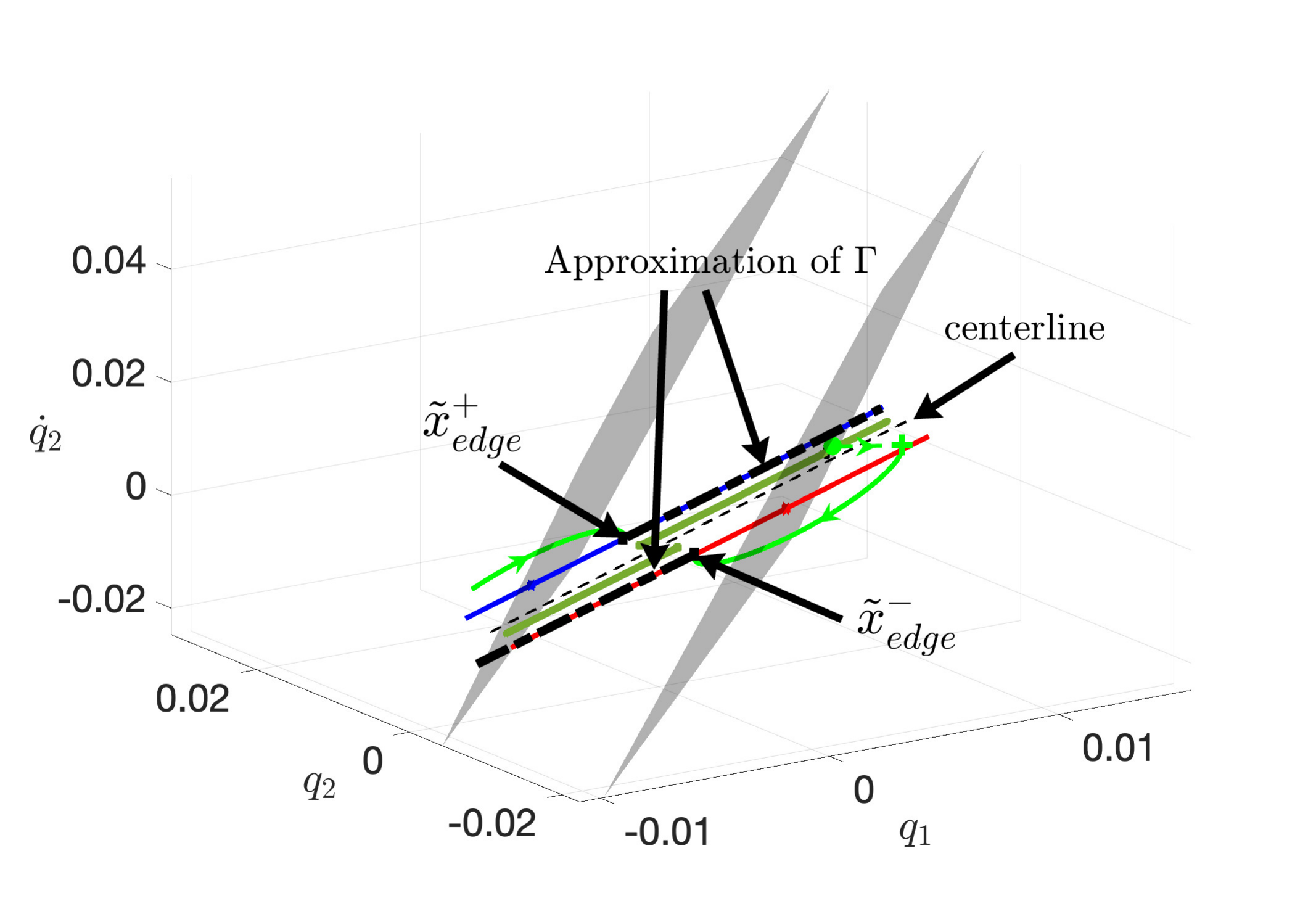}
    }
    \caption{The invariant curve $\Gamma$ of the Poincaré map for the Shaw-Pierre model with friction presents a discontinuity close to the origin. The grey surfaces delimit the region where sticking occurs and contain the two fixed points $x_0^\pm$. \protect\subref{edges_full_scaled}: $x^\pm_{edge}$ are the iterates of the configurations closest to the sticking (green squares), according to the full order model. \protect\subref{edges_reduced_scaled}: an approximation of the invariant curve $\Gamma$ is derived as a portion of the blue and red lines, namely the intersection of the positive and negative SSMs with the switching surface, where $\tilde{x}_{edge}^\pm$ are defined by the reduced-order model.}
    \label{delta=0.01_edge_full}
 \end{figure*}  
 
Studying the reduced dynamics according to these different choices, we find that the initial condition defined by \eqref{new_ic_equation} tracks the full solution most effectively (see \ref{choice_ic}). Therefore, we will use this matching scheme for initial conditions in the following. 

\subsection{Poincaré map and invariant set}\label{SP_poincare_map} We now seek an actual attractor (i.e., a closed invariant set with an open domain of attraction) near the locally invariant, attracting set $\mathcal{M}^+\cup\mathcal{M}^-$. To this end, we consider a set of initial conditions around the intersection of $\mathcal{M}^\pm$ with the switching surface in the $(q_1, q_2, \dot{q}_2)$ space. Recording the subsequent crossing of trajectories through $\Sigma$ defines a Poincaré map, which turns out to have an invariant curve $\Gamma$ with a discontinuity close to the origin (green straight piecewise continuous line in Fig. \ref{delta=0.01_edge_full}). The two segments of this invariant curve lie close to the intersection of $\mathcal{M}^\pm$ with $\Sigma$ (blue and red lines) on either side of the origin. The points at the discontinuity of $\Gamma$ are denoted $x^\pm_{edge}$. They arise from the states closest to the sticking regime (green squares in \ref{edges_full_scaled}). In these states, the sum of forces from the springs and dampers are almost equal to the static friction force, but enough to induce crossing. The invariant curve $\Gamma$ and $x^\pm_{edge}$ have been computed integrating the full model, but they can be directly approximated from the reduced-order model (see \ref{SP_supp_poincare_map}).

\subsection{Non-autonomous problem} 
Under small periodic external force applied to the mass $m_1$, system \eqref{initial_system} becomes 
\begin{equation}\label{eq:physical_nonautonomous}
\begin{aligned}
 &   \begin{pmatrix}
    m_1 & 0\\
    0 & m_2
    \end{pmatrix}
    \begin{pmatrix}
    \ddot q_1 \\
    \ddot q_2
    \end{pmatrix}
    +
    \begin{pmatrix}
    c & -c\\
    -c & 2c
    \end{pmatrix}
    \begin{pmatrix}
    \dot q_1 \\
    \dot q_2
    \end{pmatrix}
   +
    \begin{pmatrix}
    2k & -k\\
    -k & 2k
    \end{pmatrix}
    \begin{pmatrix}
    q_1 \\
    q_2
    \end{pmatrix}
\\&     +
    \begin{pmatrix}
    -\alpha q_1^3 \\
    0
    \end{pmatrix}
    \pm
    \begin{pmatrix}
    \delta m_1\\
    0
    \end{pmatrix}
    = \epsilon 
    \frac{1}{\sqrt{2}}\begin{pmatrix}
        1\\1
    \end{pmatrix}\cos{\Omega t},
\end{aligned}
\end{equation}
where $0 \leq\epsilon< < 1$. In terms of the coordinates introduced in eq. \eqref{shifted_coordinates} and the quantities defined in eq. \eqref{shifted_terms}, we have 
\begin{equation}\label{eq:x_coordinates_nonautonomous}
    \dot {\mathbf{\xi}}^\pm =
    \mathbf{\tilde{A}}_0 \mathbf{\xi}^{\pm} \pm \mathbf{f_{nl}^{II}}(\mathbf{\xi}^\pm) + \mathbf{f_{nl}^{III}}(\mathbf{\xi}^\pm) + \mathbf{f}_0 + \epsilon \,\mathbf{f}_\epsilon(\Omega t) ,
\end{equation}
with
\begin{equation*}
    \mathbf{f}_\epsilon = \frac{1}{\sqrt{2}}\begin{pmatrix}
        0 \\ \frac{1}{m_1} \\ 0 \\ \frac{1}{m_2}
    \end{pmatrix}\cos{\Omega t} = 
    \mathbf{f}_0 \cos{\Omega t}.
\end{equation*}
To obtain an approximation for the time-dependent SSMs on the two sides of the switching surface, we again rely on a cubic Taylor expansion, but with the addition of a $2\pi/\Omega$-periodic time-dependent term (the details are reported in \ref{supp_SP_na}). 

We now compare the SSM-reduced model with the full-order one in terms of forced response curves computed for different forcing amplitudes $|\mathbf{f}_0|$ and friction coefficients $\delta$. For each case, we compute the response of the full system using the numerical continuation software COCO of \citet{COCO} for a range of forcing frequencies. In contrast, results of the reduced model come from the direct integration of the reduced dynamics. 

 \begin{figure*}[]
     \centering
     \subfloat[$\delta = 10^{-3}$ \label{FRC_1e-3.}]{
         \includegraphics[scale=0.145]{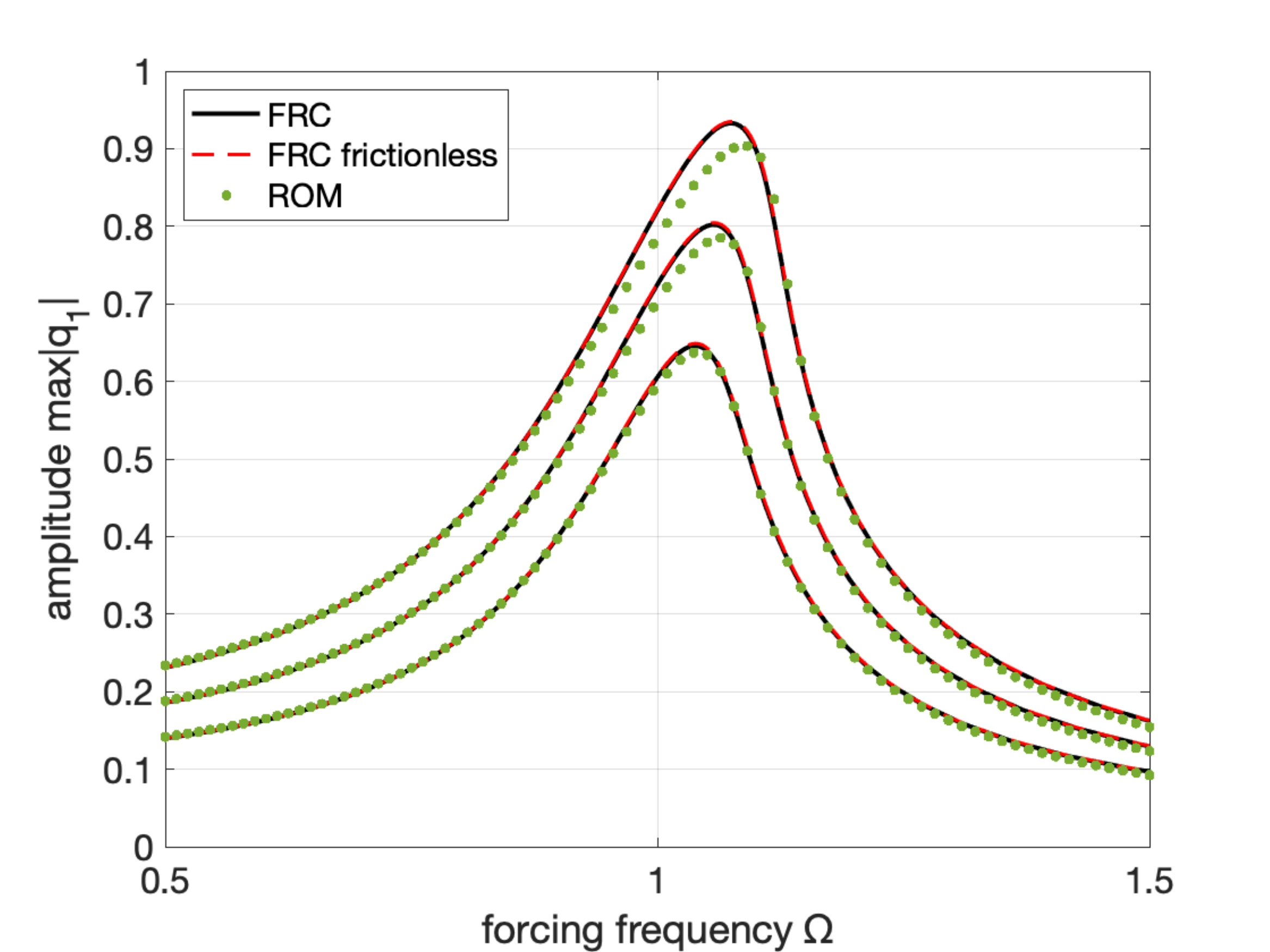}
     }
     \quad
     \subfloat[$\delta = 5 \cdot 10^{-3}$ \label{FRC_5_1e-3}]{
         \includegraphics[scale=0.145]{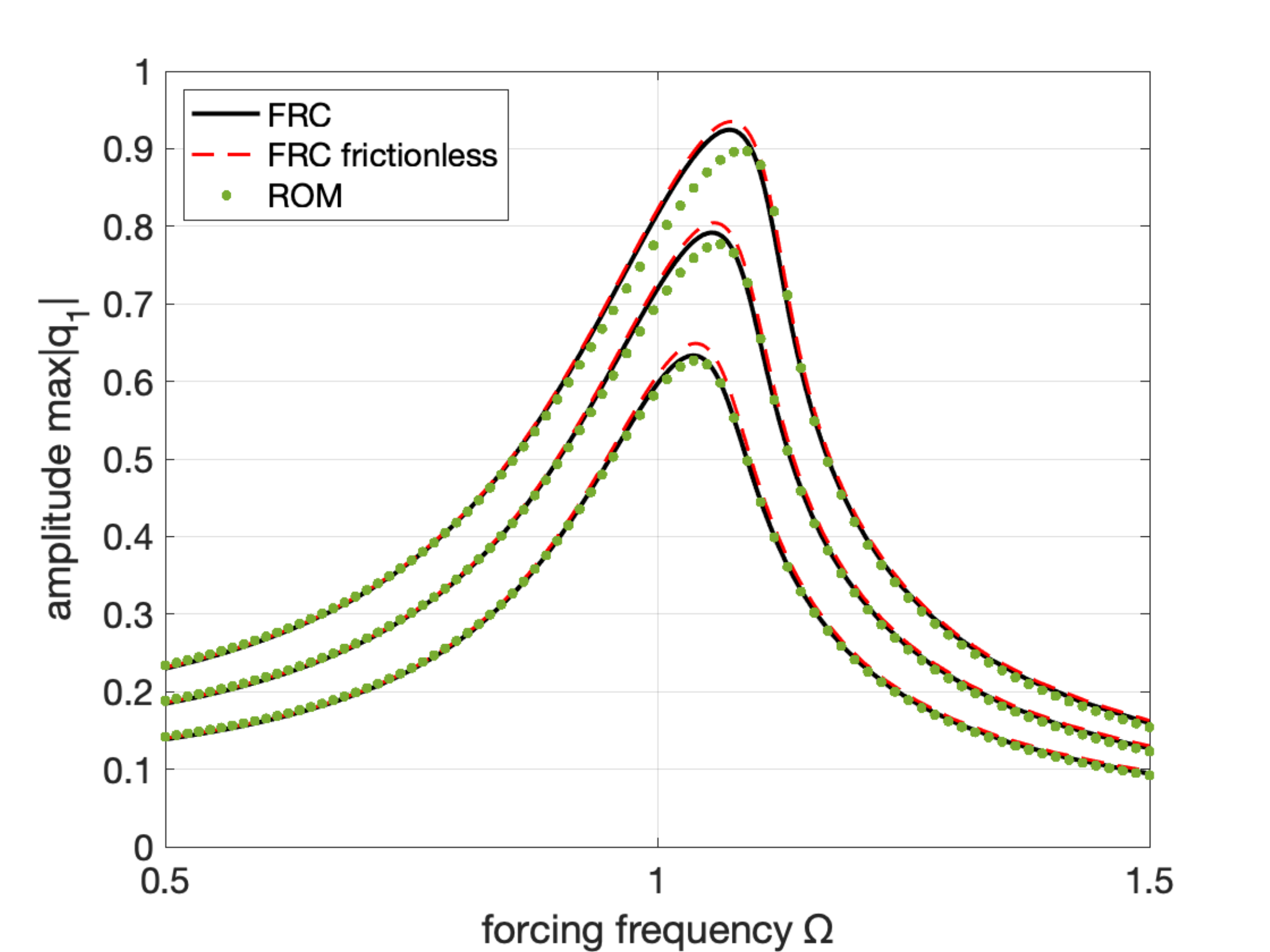}
     }
     \\
     \subfloat[$\delta = 10^{-2}$ \label{FRC_1e-2}]{
         \includegraphics[scale=0.145]{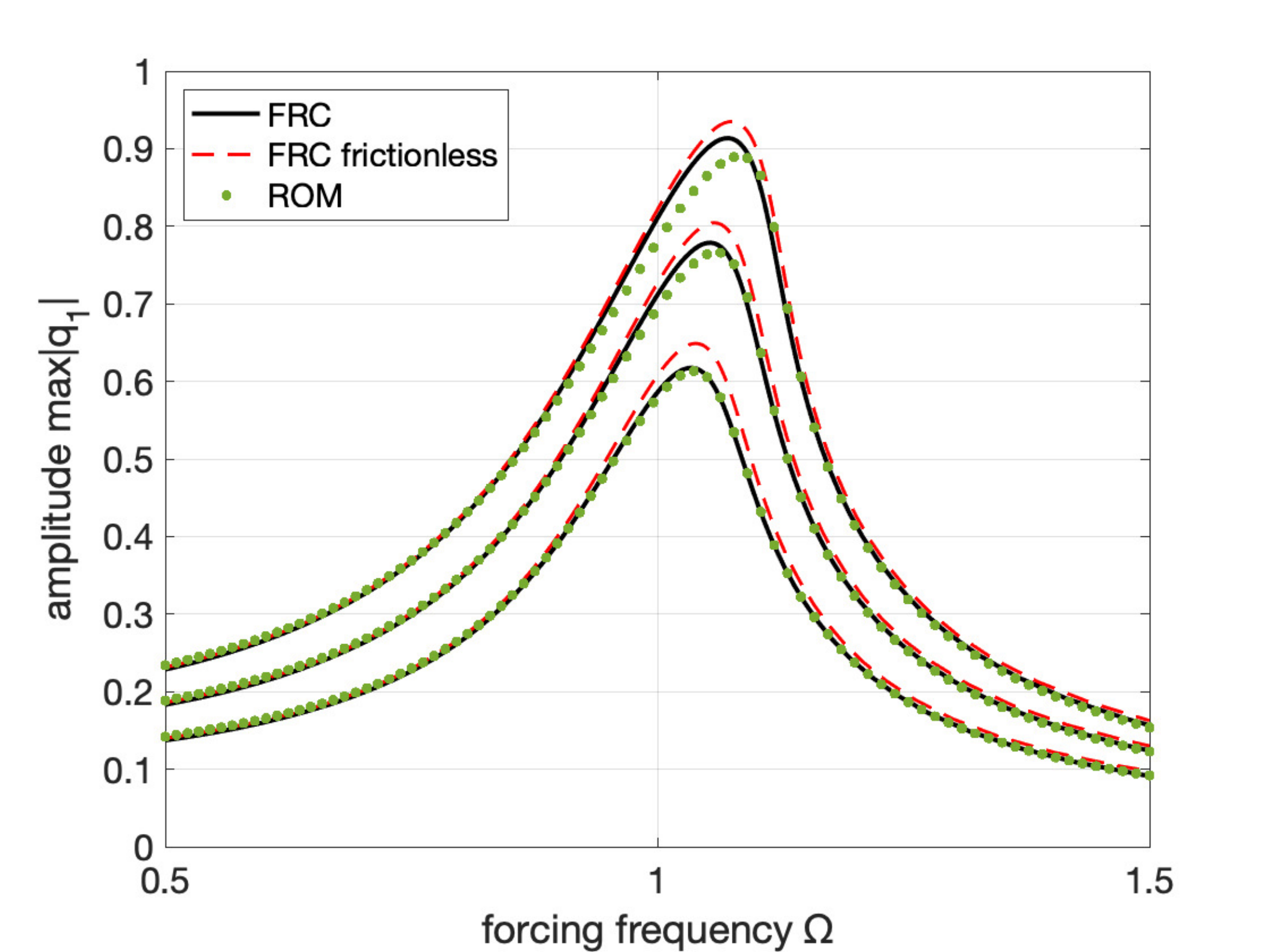}
     }  
     \quad
     \subfloat[$\delta = 5 \cdot 10^{-2}$ \label{FRC_5_1e-2}]{
         \includegraphics[scale=0.145]{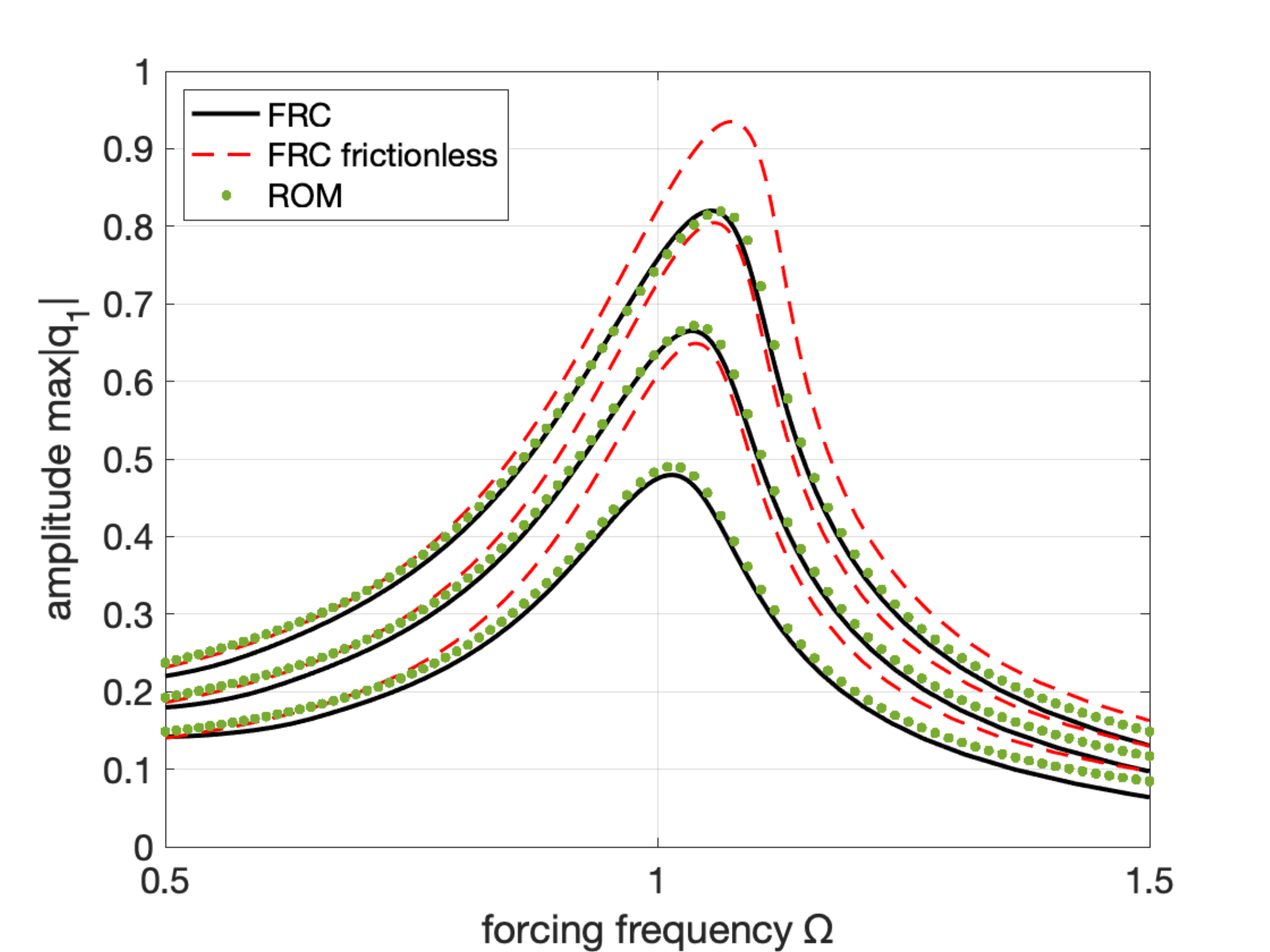}
     }
     \caption{Forced response curves for the forced Shaw-Pierre example with friction, obtained for different values of the friction coefficient $\delta$. For each case, the amplitudes of the forcing ($\epsilon$ in equation \eqref{eq:physical_nonautonomous}) are $0.15$,  $0.2$ and $0.25$. The black curves represent the solution of the full-order model; the red dashed curves indicate the forced response without friction; the set of green dots represents the reduced-order model approximation.}
     \label{FRCs}
 \end{figure*}

Figure \ref{FRCs} shows highly accurate predictions for forced responses for the smallest value of the forcing amplitude $\epsilon$, even for large values of $\delta$. The prediction error, however, grows for larger forcing amplitudes as we keep the order of approximation for the SSMs cubic. Nevertheless, the reduced order model correctly captures the effect of friction in the system, at least for a range of parameters, as visible especially in Fig. \ref{FRC_5_1e-2}. Indeed, in this case the difference between the system with and without friction is the most evident.

\section{Data-driven model reduction of a non-smooth beam model}\label{vonkarman_beam}
We now apply our nonsmooth SSM-reduction procedure in a data-driven setting. We consider a finite-element model of a von Kármán beam with clamped-clamped boundary conditions (\citet{shobit_2017}). 
Each finite element has three degrees of freedom: axial deformation, transverse displacement and transverse rotation. The beam is approximated by 4 elements, resulting in 9 total degrees of freedom, i.e., an 18-dimensional phase space. The material is aluminium, with Young modulus $E = 70 \,\text{GPa}$, density $\rho = 2700 \,\text{Kg/}\text{m}^3$, Poisson ratio $\nu = 0.3$ and material damping modulus $k = 1\times 10^6 \, \text{Pa}\cdot\text{s}$.  We set the length 1 [m], the width 5 [cm] and the thickness 2 [cm]. 

 \begin{figure}[]
     \centering
     \subfloat[\label{vkbeam_figure}]{
         \includegraphics[scale=0.2]{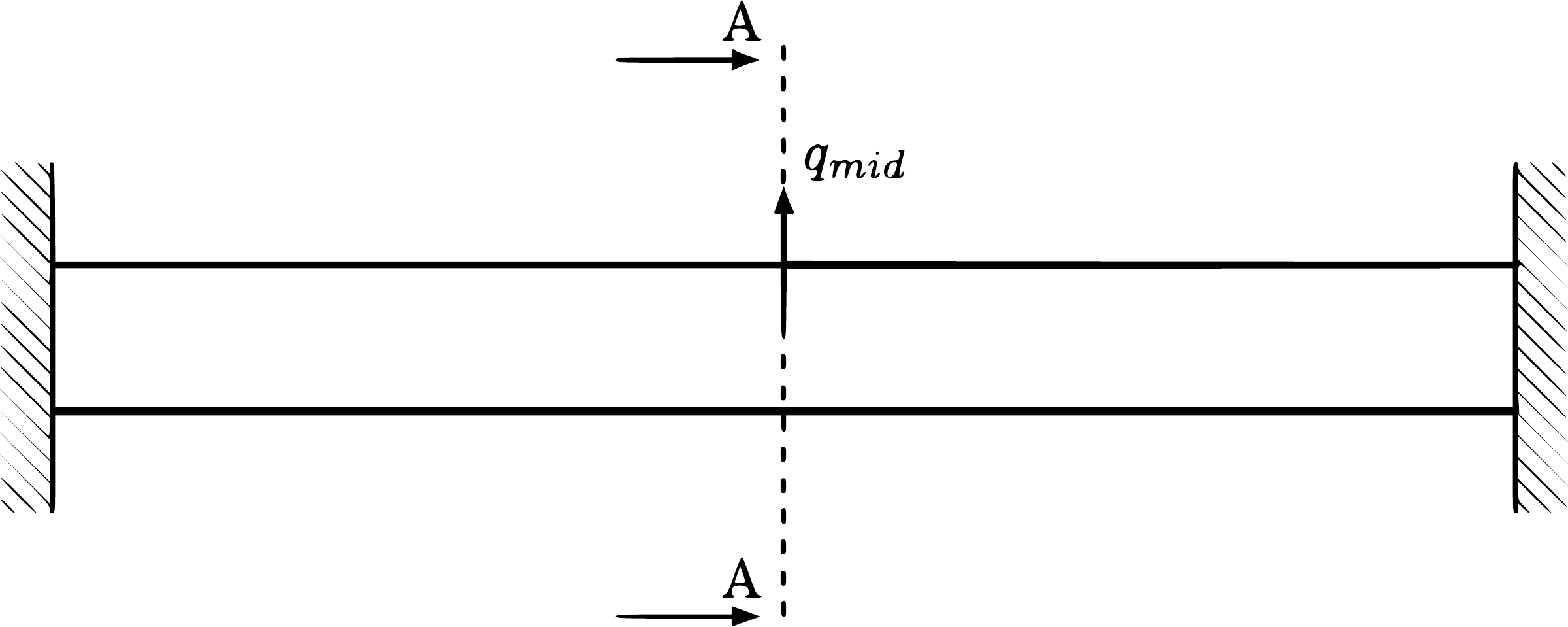}
     }\\
     \subfloat[\label{coulomb_friction}]{
         \includegraphics[scale=0.3]{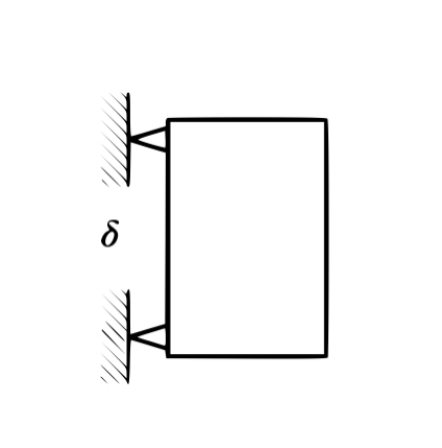}
     }\quad
     \subfloat[\label{soft_impact}]{
         \includegraphics[scale=0.3]{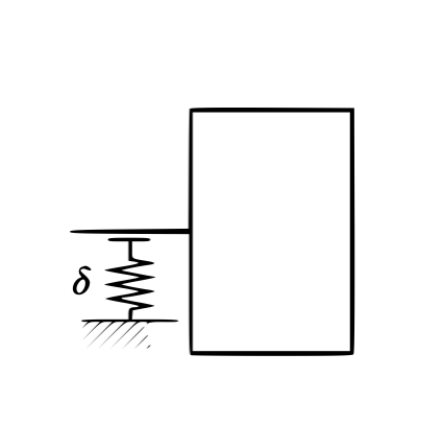}
     }\quad
     \subfloat[\label{moving_ground}]{
         \includegraphics[scale=0.3]{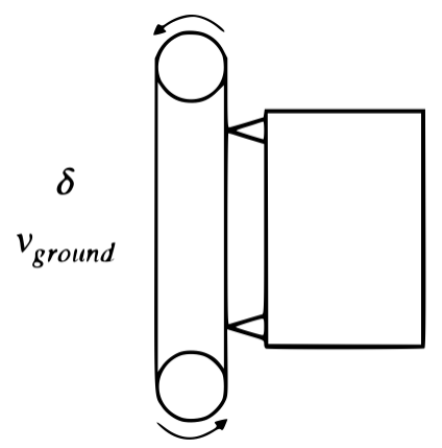}
     }
     \caption{\protect\subref{vkbeam_figure} Geometry of the clamped-clamped von Kármán beam. In each of the cases shown in the subplots \protect\subref{coulomb_friction}, \protect\subref{soft_impact}, \protect\subref{moving_ground}, the non-smoothness is localized at the midpoint section A-A, where $q_{mid}$ represents the vertical displacement of the beam in that section. \protect\subref{coulomb_friction}: Coulomb friction, \protect\subref{soft_impact} soft impact, \protect\subref{moving_ground}  friction on moving ground.}
     \label{vk_beam_scenarios}
 \end{figure}

We consider three different cases of non-smoothness introduced at the middle of the beam (section A-A of Fig. \ref{vkbeam_figure}), with the parameter $\delta$ playing a different role in each case. 

\begin{enumerate}
    \item \textbf{Coulomb friction:} Dry friction between the node at the middle gives the switching function
    \begin{equation}
        \sigma(\mathbf{x}) = \dot{q}_{mid},
    \end{equation}
    and the switching surface
    \begin{equation}
        \begin{array}{l}
             \Sigma^\pm = \left\{\mathbf{x} \in \mathbb{R}^N : \text{sign}(\sigma(\mathbf{x})) = \pm 1\right\},
        \end{array}
    \end{equation}
    with $N$ denoting the phase space dimension of the discretized beam model.
    Once discretization has been performed, the equations read
    \begin{equation}\label{governing_equations_beam}
        M\ddot{\mathbf{q}} + C\dot{\mathbf{q}} + K\mathbf{q} + \mathbf{f}_{nl}(\mathbf{q}, \dot{\mathbf{q}}) = \mathbf{f}_\delta^\pm,
    \end{equation}
    where
    \begin{equation}\label{forcing_beam_coulomb}
        \mathbf{f}_\delta^\pm =
        (0\,\cdot \,\cdot \,\cdot \,0 ,\,\mp\delta ,\, 0 \,\cdot \, \cdot \, \cdot \, 0)^T.
    \end{equation}
    In eq. \eqref{forcing_beam_coulomb}, the forcing acts only on the element located at the midpoint. The behavior within the switching surface is qualitatively the same as in the Shaw-Pierre example of Section \ref{shaw_pierre}, in that both crossing and attractive sliding (sticking) regimes are possible. We define a normalized non-smooth coefficient $\tilde{\delta}$ by dividing $\delta$ by the maximal elastic force at the midpoint.
    
    \item \textbf{Soft impact:} Impact at the section A-A gives the switching function
    \begin{equation}
        \sigma(\mathbf{x}) = q_{mid}.
    \end{equation}
    Formally, eq. \eqref{governing_equations_beam} is still valid, but
    \begin{equation}\label{equation:forcing_beam_coulomb}
        \mathbf{f}_\delta^+ = \mathbf{0}, \quad \text{and} \quad \mathbf{f}_\delta^- = 
        (0\,\cdot \,\cdot \,\cdot \,0 ,\,-\delta x_{mid} ,\, 0 \,\cdot \, \cdot \, \cdot \, 0)^T.
    \end{equation}

    The behavior of the system is asymmetric with respect to the $q_{mid} = 0$ position, due to the presence of the further stiffness on one side, as seen in Fig. \ref{decaying_soft_impact}. 
    
     \begin{figure}[H]
     \centering
     \includegraphics[scale=0.6]{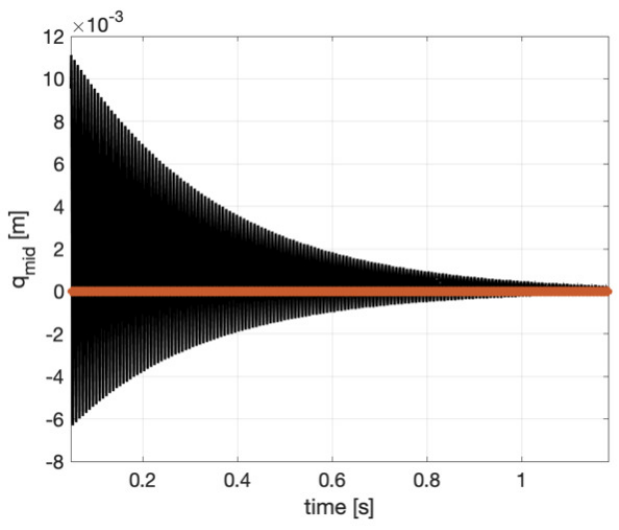}
     \caption{Decaying trajectory of the von-Kármán beam with soft-impact located at the midpoint. The orange line refers to the switching between the two sub-regions.}
     \label{decaying_soft_impact}
     \end{figure}

 According to conditions \eqref{crossing_condition} and \eqref{sliding_condition}, only crossing between the regions $\Sigma^+$ and $\Sigma^-$ is admissible, while sticking is not allowed. In this case, the coefficient $\delta$ represents the increased stiffness acting on the midpoint of the beam for $q_{mid} <0$ and it is naturally normalized by the linear stiffness of the beam.

    \item \textbf{Friction on moving ground: } Assume that the midpoint of the beam rides on a belt moving with constant velocity $v_{ground}$ (Fig. \ref{moving_ground}) and dry friction is present between the beam and the belt. The non-smooth forcing term is then given by
     \begin{equation}
         \mathbf{f}_\delta^\pm = 
    (0\,\cdot \,\cdot \,\cdot \,0 ,\,\delta f_{ns} ,\, 0 \,\cdot \, \cdot \, \cdot \, 0)^\mathrm{T},
     \end{equation}

     where
     \begin{equation}\label{friciton_law_moving_belt}
     \begin{aligned}
         &f_{ns}(\dot{q}_{mid}) = -\text{sign}\,(\dot{q}_1 - v_{ground})\\&\quad\times\left( 1 + \frac{\alpha}{e}\, \text{exp}\, \left( \frac{\beta - |\dot{q}_{mid} - v_{ground}|}{\beta} \right) \right),
     \end{aligned}
     \end{equation}
as shown in Fig. \ref{image_friction_law_moving_belt}. This friction model, similar to that used in \citet{leine}, takes into account that the value of the static friction to be overpowered in order to trigger a relative motion is higher than the kinetic friction force when the relative motion is different from zero.   
\begin{figure}
	\centering
	 \includegraphics[scale=0.13]{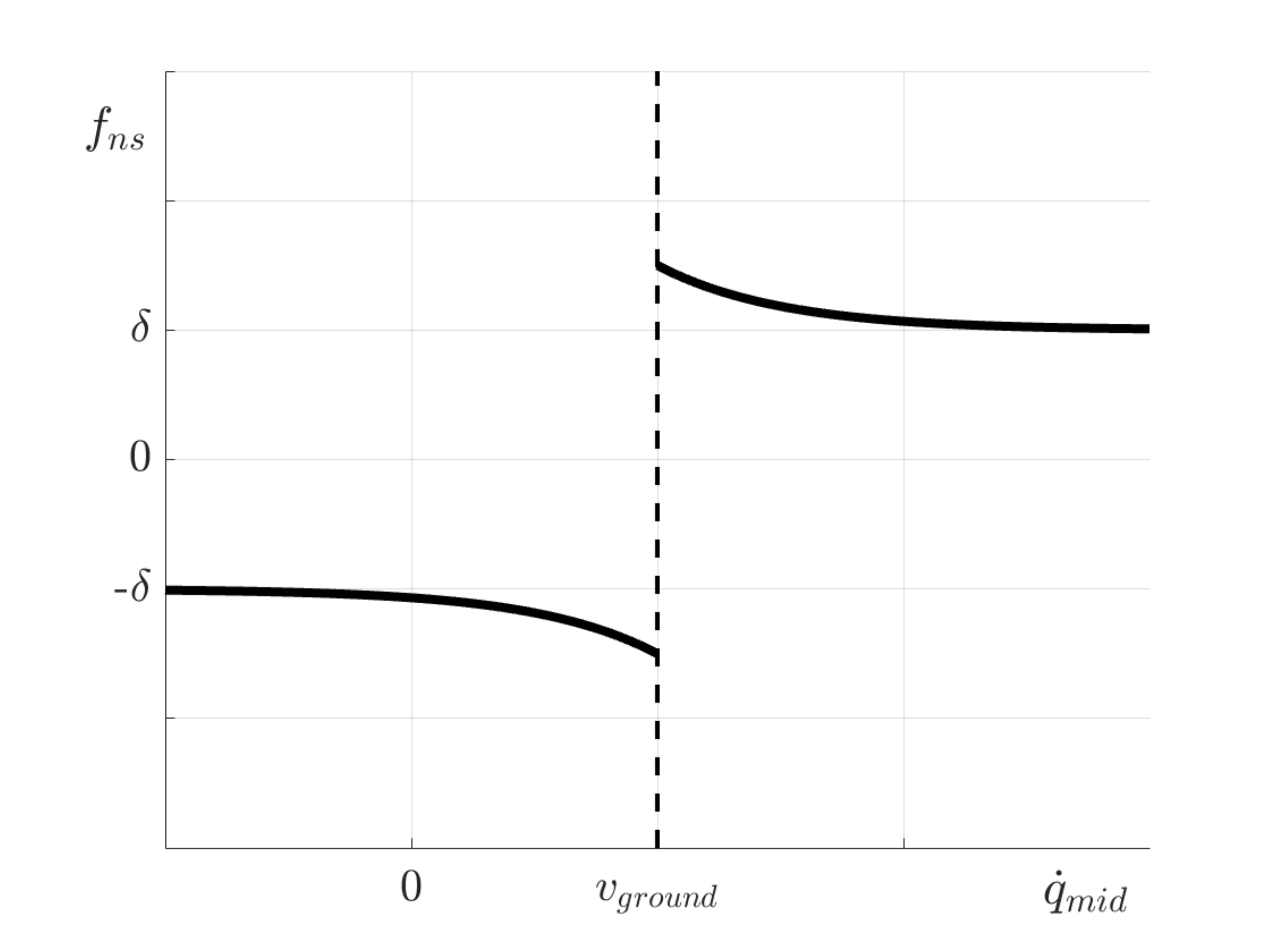}
	 \caption{Friction law dependent on the relative velocity between the midpoint of the von Kármán beam and the moving belt. In contrast to the Coulomb friction model, here the static friction force is different from the dynamic one.}
	 \label{image_friction_law_moving_belt}
\end{figure}
The switching function in this case is
     \begin{equation}
         \sigma(\mathbf{x}) = \dot{q}_{mid} - v_{ground},
     \end{equation}
    which splits the phase space according to the relative motion between the beam element at the midpoint and the moving belt. Interestingly, the difference between static and dynamic friction forces causes the fixed points of the positive and negative systems to be unstable under certain parameter values, triggering a limit cycle. Also in this case both crossing and sticking are allowed: the presence of sticking to the switching surface represents  indeed a crucial factor in sustaining the limit cycle. 
\end{enumerate}

The data-driven procedure for the computation of the reduced order model is similar for all three piecewise-smooth systems. We constructed two different slow SSMs in each case, one for each region $\Sigma^+$ and $\Sigma^-$. Within each of these regions, the beam equations are analytic and hence the SSMs are approximated by a convergent Taylor expansion near the respective equilibrium points.  We generate decaying trajectories which serve as training data, under $\mathbf{f}^+$ and $\mathbf{f}^-$, separately. The initial condition is a deformed configuration caused by transverse static loading of 12 kN at the midpoint. 

\begin{figure*}
     \centering
     \subfloat[Coulomb friction, $\tilde{\delta} = 1\cdot 10^{-3}$\label{VK_coulomb_response}]{
         \includegraphics[scale=0.3]{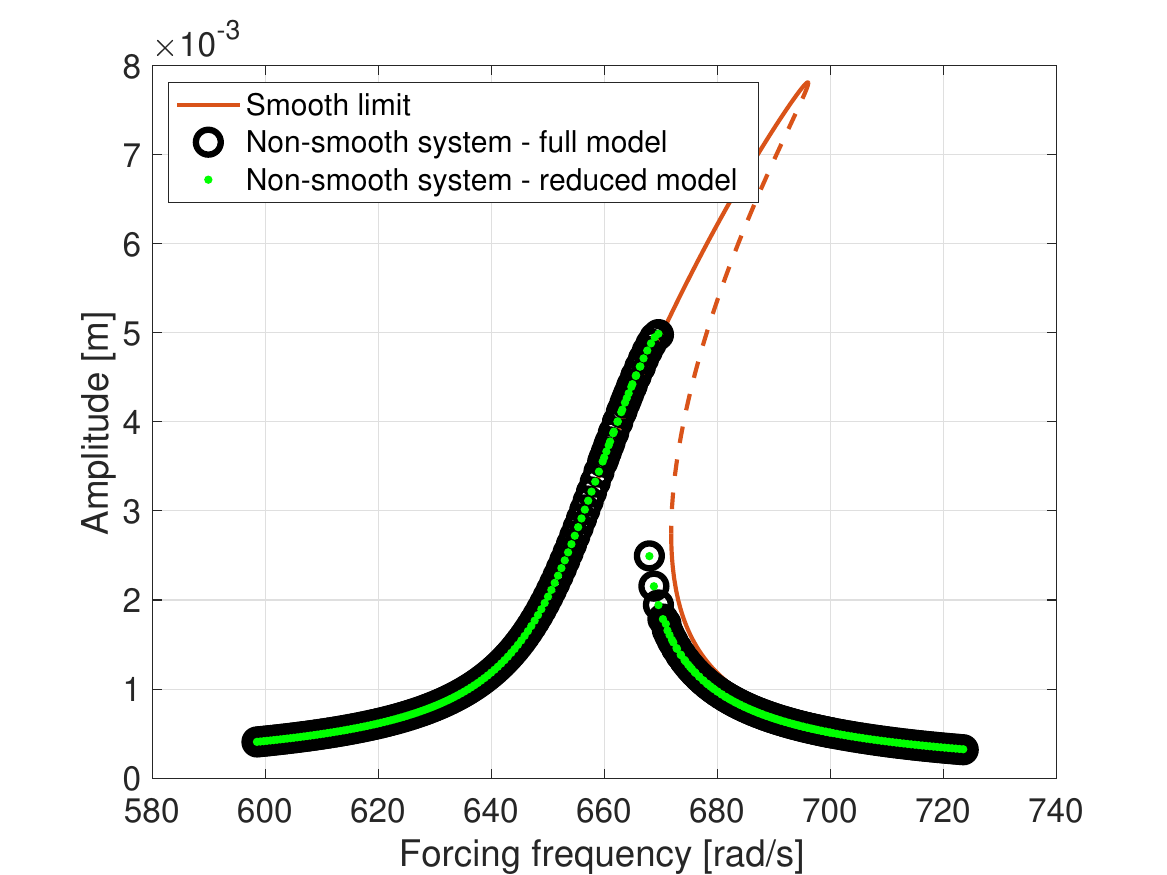}
     }\quad
     \subfloat[soft impact, $\Tilde{\delta} = 5\cdot10^{-4}$\label{VK_soft_imp_response}]{
         \includegraphics[scale=0.3]{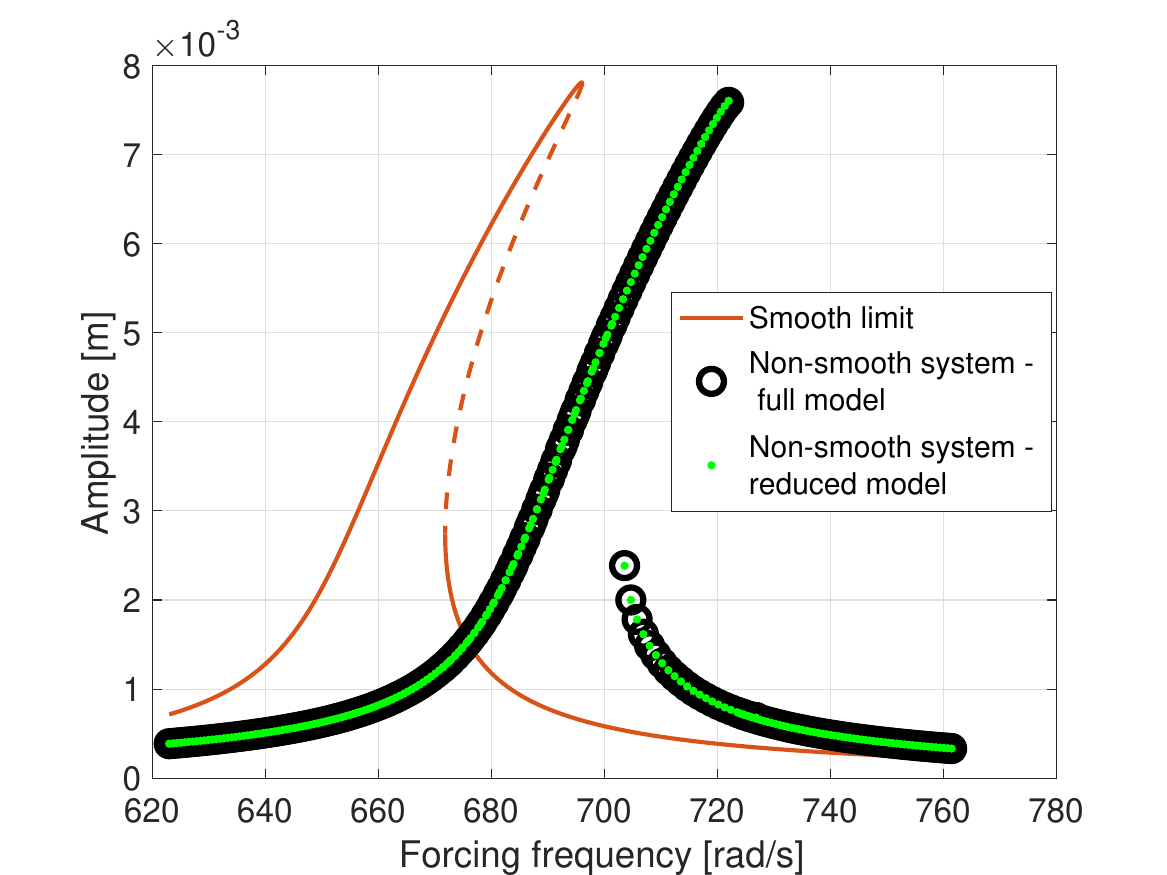}
     }
     \caption{Response curves for the von Kármán beam forced periodically at the midpoint with Coulomb friction \protect\subref{VK_coulomb_response} and soft impact \protect\subref{VK_soft_imp_response} with the same forcing level (35 kN). In both cases, the effect of non-smoothness is evident. The smooth case ($\tilde{\delta} = 0$) and full model cases are computed via the continuation software COCO \cite{COCO}, while the response of the reduced model is retrieved by numerical integration.}
     \label{VK_response}
 \end{figure*}

We seek a two-dimensional SSM, approximated by a 5th-order expansion and a 5th-order approximation of the reduced dynamics. According to Whitney's embedding theorem \cite{whitney_embed}, the minimal embedding dimension for a two-dimensional SSM is $p=5$, which is satisfied for our specific case, as we observe the full system. The error of the reconstructed SSM-based reduced model is quantified by the normalized mean-trajectory error (NMTE, see \citet{mattia_2022}) for each region $\Sigma^+$ and $\Sigma^-$. For a data set of $P$ instances of observable points $\mathbf{y}_j \in \mathbb{R}^p$ and their reconstruction $\mathbf{\hat{y}}_j$:
\begin{equation}
    \text{NMTE} = \frac{1}{\|\mathbf{\underline{y}}\|}\frac{1}{P}\sum_{j=1}^P \|\mathbf{y}_j - \mathbf{\hat{y}}_j \|,
\end{equation}
where $\mathbf{\underline{y}}$ is a relevant normalization vector.
Finally, we combine the two SSM-based reduced models according to their individual ranges of validity dictated by the switching function. We then employ the model trained on decaying trajectories to predict force response. Based on the arguments of \citet{mattia_2022}, moderate forcing terms can simply be added to the reduced model constructed from unforced data.
\begin{figure}\label{reconstruction}
	\centering
	 \includegraphics[scale=0.4]{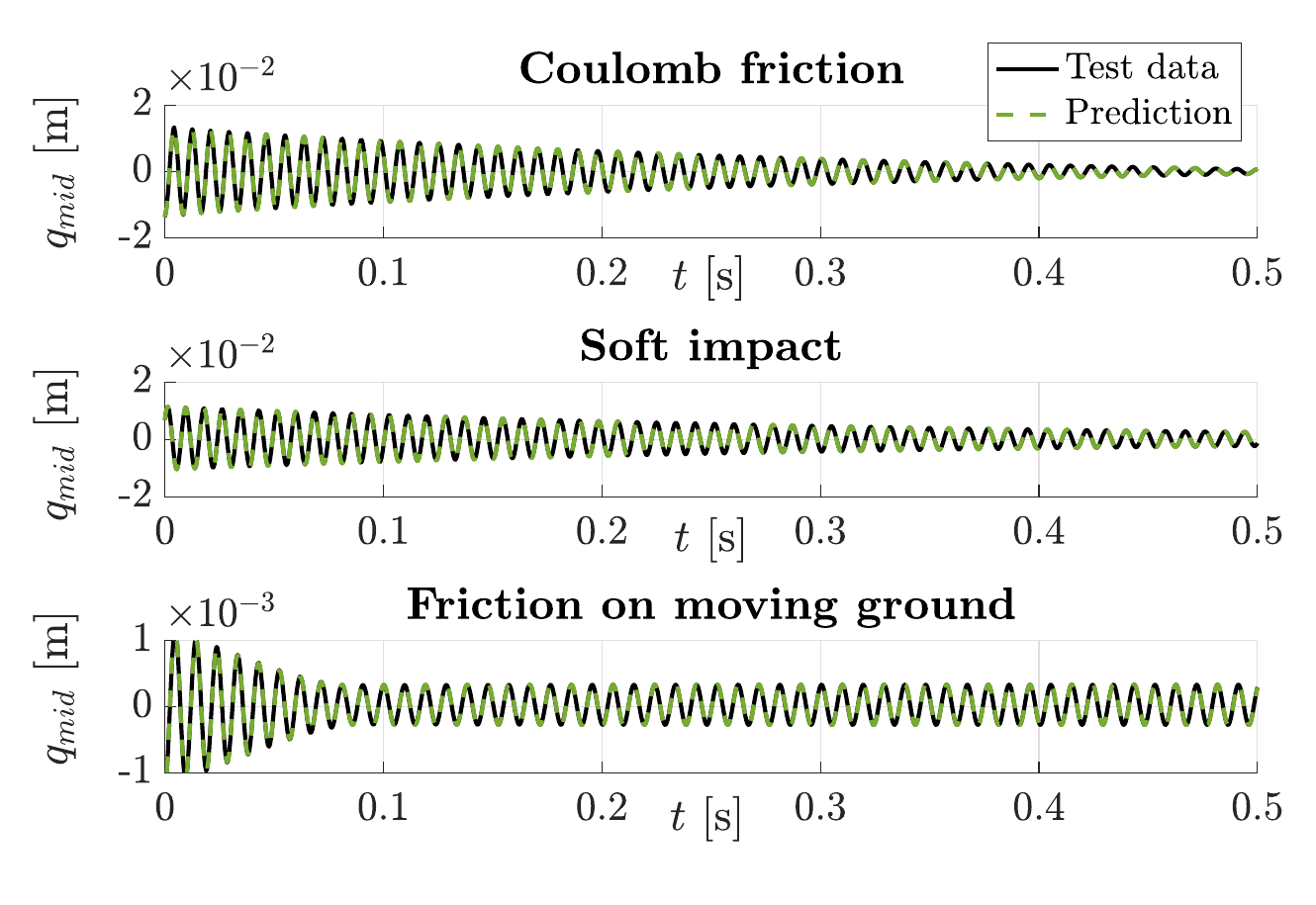}
	 \caption{Model-testing trajectories of the unforced von Kármán beam with different kinds of non-smoothness, and their reconstructions from an unforced, nonsmooth, SSM-based model.}
\end{figure}
Here, we consider a transversal periodic forcing applied to the midpoint. As seen in Fig. \ref{VK_response} for the Coulomb friction and soft impact cases, the piecewise-smooth reduced order model closely tracks the response for a range of frequencies around the first mode of the system, even when the effect of the non-smoothness factor $\tilde{\delta}$ is not negligible. Increasing $\Tilde{\delta}$ intensifies the dissipation in the system, entailing a lower amplitude of the response. In the soft impact case, $\Tilde{\delta}$ strengthens the stiffness in some operating region, which leads to an overall shift of the resonance towards higher frequencies. 
 
\begin{figure}[H]
     \centering
     \subfloat[\label{p_map_phase}]{
         \includegraphics[scale=0.2]{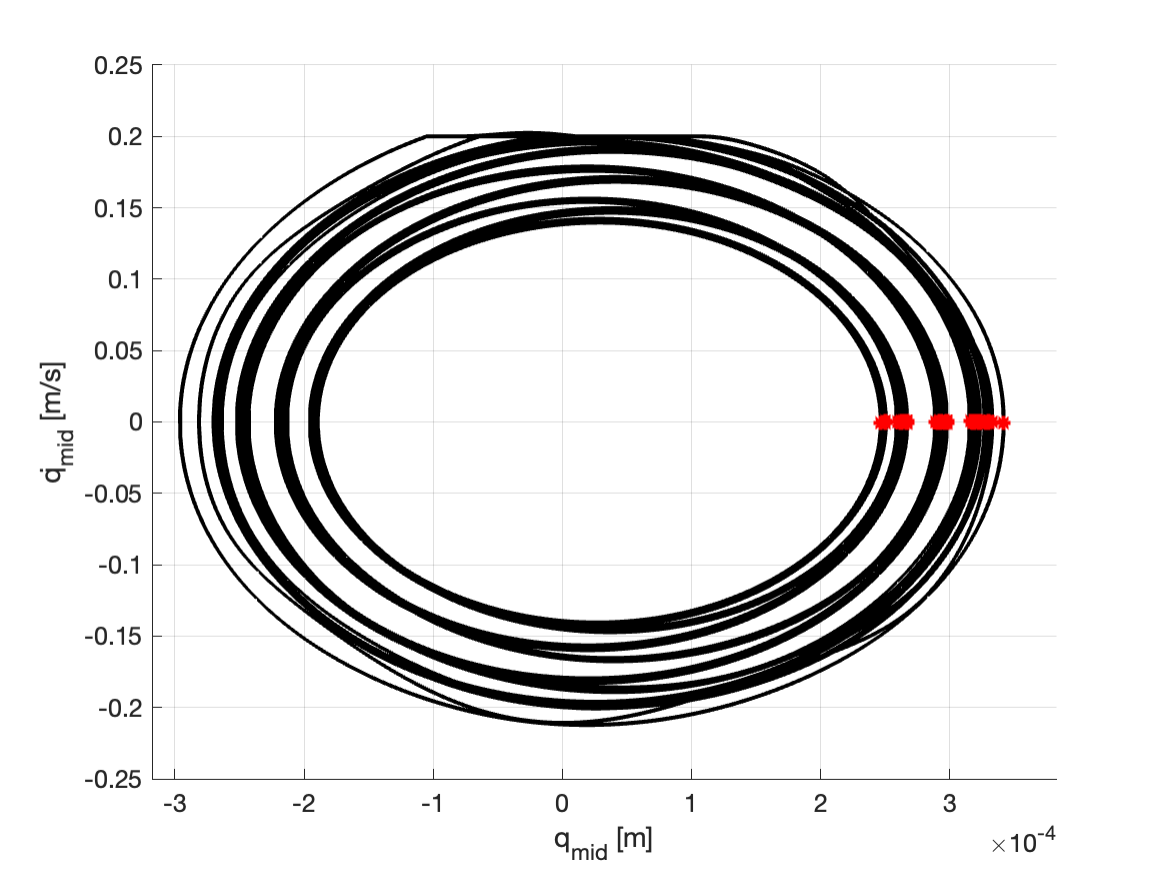}
     }\quad
     \subfloat[\label{p_map}]{
         \includegraphics[scale=0.2]{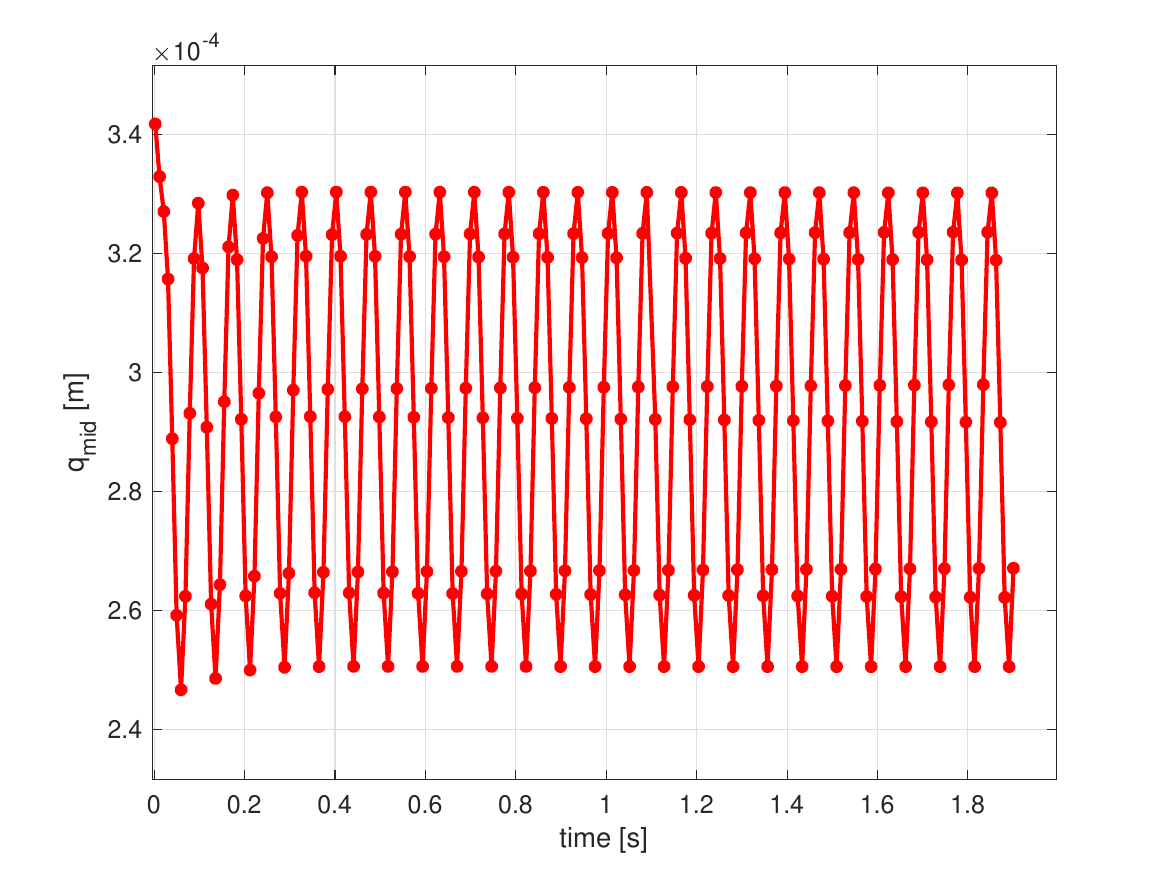}
     } \\
     \subfloat[\label{forced_reduced_phase}]{
         \includegraphics[scale=0.2]{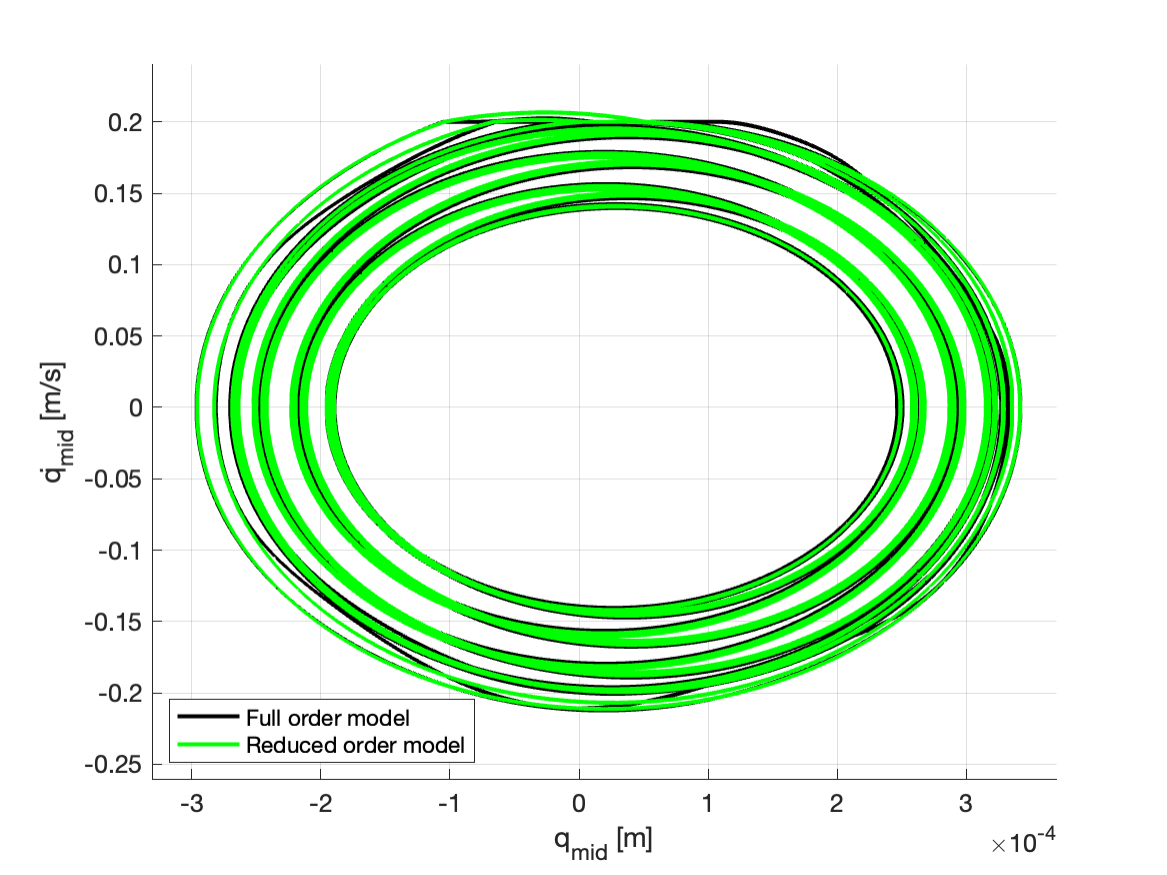}
     }\quad
     \subfloat[\label{forced_reduced_time_colored}]{
         \includegraphics[scale=0.1]{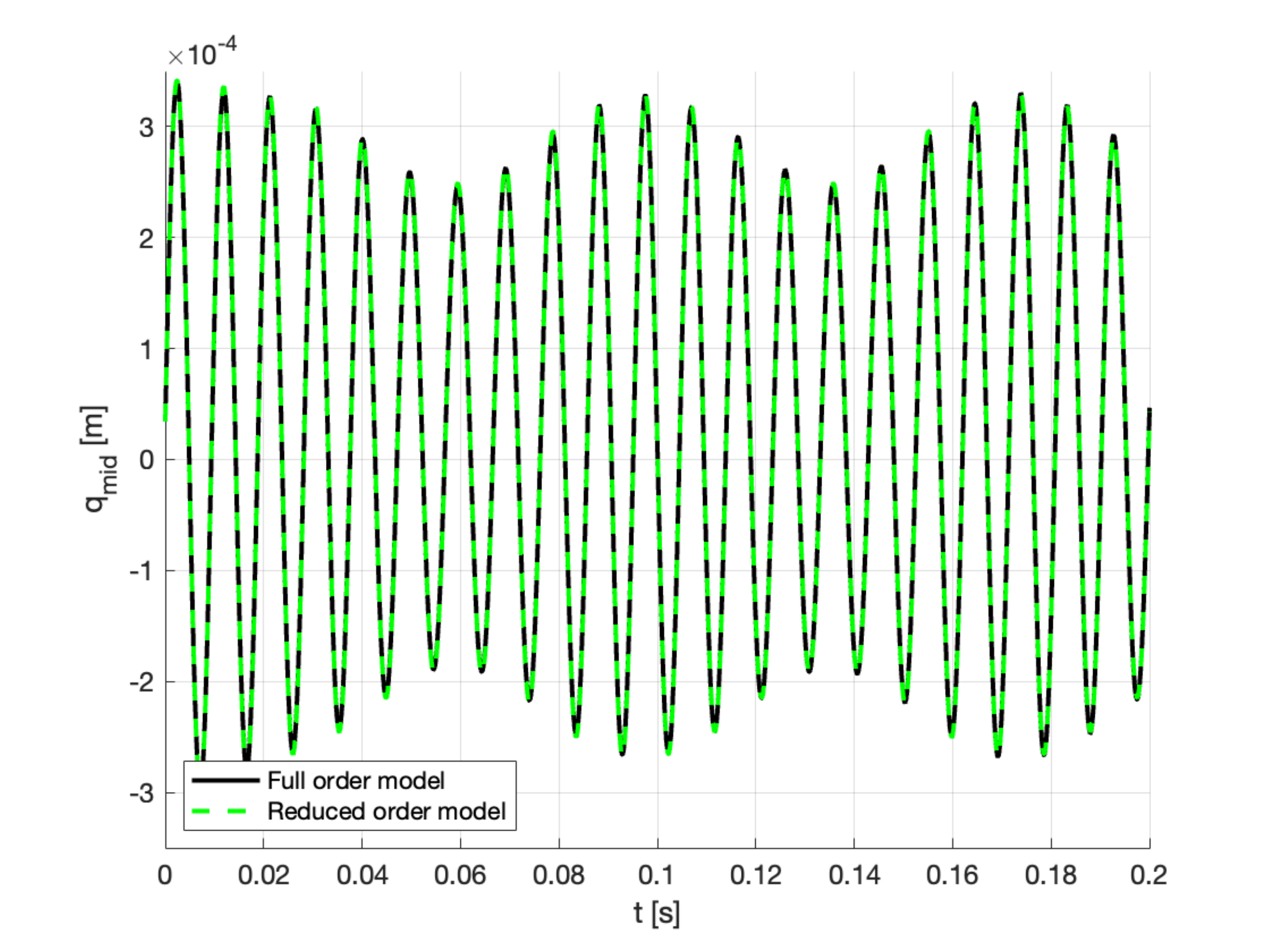}
     }
     \caption{Periodic orbit lying on a torus generated by external forcing ($50 N$, $f_{F}/f_{LC} = 1.125$ ) applied to a von Kármán beam with friction on a moving ground, with unstable slowest eigenvalues. Subfigure \protect\subref{p_map} shows the temporal history of the red dots in subfigure \protect\subref{p_map_phase}. Subfigures \protect\subref{forced_reduced_phase} and \protect\subref{forced_reduced_time_colored} compare the reduced-order model (in green) with the full order one (in black).}
     \label{VK_moving_ground_forcing}
 \end{figure}

 Modeling the sticking regime is delicate, as the constraints confining the trajectory within the switching surface must be enforced. This becomes especially relevant when the dynamics enter and exit the sticking regime multiple times, as in the case of the limit cycle that characterizes the beam with friction on moving ground, for a range of parameters. Since the condition for the sticking regime to occur is expressed in terms of physical variables, choosing the physical coordinates (vertical displacement and velocity at the midpoint) as the reduced variables for the SSM parametrization is advantageous. This is not the standard choice, as we usually employ a pair of coordinates spanning the slowest spectral subspace in order to describe a two-dimensional SSM. More details are given in the \ref{alternative_param}. 
Let us consider a configuration of parameters wherein the limit cycle is triggered for the autonomous system already: adding forcing in this case generates an invariant periodic or quasiperiodic torus. The behavior of the response depends on the relationship between the frequency of the limit cycle of the autonomous system ($f_{LC}$) and that of the forcing term ($f_{F}$). In Fig. \ref{VK_moving_ground_forcing} we report the case for $f_{F}/f_{LC} = 1.125$, where the reduced order model (in green) is able to accurately capture the multi-frequency periodic solution of the full order system (in black).
 
 \section{Conclusions}\label{conclusions}
We have presented a nonlinear, SSM-based model reduction procedure applicable to piecewise smooth dynamical systems. Our approach extends the applicability of the spectral submanifolds to a class of problems where non-smoothness plays a significant role in the dynamics of the system. In particular, the classic SSM theory is applied separately to different regions of the phase space where the system is smooth, near hyperbolic fixed points. Proper matching conditions are then enforced as crossing between those regions occurs.

The method proposed here accurately captured the non-smooth dynamics in both equation- and data-driven examples. The equation-driven examples included forced and unforced versions of a 2-DOF nonlinear oscillator with Coulomb friction. In this case, the SSMs involved were computed analytically. 

The data-driven examples included forced and unforced von Kármán beam models with different types of non-smooth elements at their midpoints. In this case, the SSMs involved were calculated using an extended version of the SSMLearn algorithm.
Not only did the piecewise-smooth reduced order model correctly capture single trajectories for all these different sources of non-smoothness, but it also accurately predicted the forced response. When the sticking regime played a crucial role in the system dynamics (a limit cycle existed for the beam with friction on moving ground), we proposed a specific parametrization of the SSMs based on physical coordinates. We did this to automatically verify the sticking conditions, as they are expressed in those coordinates.

In the present work, we assumed that the phase space is split in two regions, where the system is smooth. We then applied SSM reduction procedure with only one switching surface. In principle, the same procedure can be generalized to multiple switching surfaces, but it would require a careful choice of the reduced physical coordinates for the parametrization of SSMs. Indeed, if one wants to correctly track the sticking regime to different switching surfaces, different relevant physical coordinates might be needed. A further limitation of our method lies in the small value of the normalized non-smoothness parameter, for which the two equilibria, and hence their respective SSMs, lie close to each other. 

Possible future work could investigate a more invasive effect of non-smooth elements to the system dynamics. In such a case, the condition relating points across a switching surface we have presented might not be able to furnish a proper initial condition in order for the reduced trajectory to synchronize with the solution of the full system. A new strategy taking into account the dynamics normal to the SSMs may then be required.\\

\section*{Data availability}
The data discussed in this work are publicly available at \url{https://github.com/leonardobettini/nonsmoothSSM}.

\section*{Code availability}
The code supporting the results of this work is publicly available at \url{https://github.com/leonardobettini/nonsmoothSSM}.

\section*{Author contributions}
MC and GH designed the research. MC, LB and GH carried out the research. MC and LB developed the software and analyzed the examples. LB wrote the paper. GH lead the research team.

\section*{Competing interests}
The authors declare no competing interests. 

 \appendix
\section{Piecewise smooth systems}
\subsection{Filippov's closure for piecewise-smooth systems}\label{supp_PWS}

In order to extend the validity of the governing equations \eqref{pws_system} to the switching surface $\Sigma$, we define an extended system, according to the Filippov's convex inclusion method
\begin{equation}\label{pws_system_extended}
    \dot{\mathbf{x}} = \mathbf{F}( \mathbf{x}) =  
    \begin{cases}
    \mathbf{f}^+(t,\mathbf{x}), & \mathbf{x} \in \Sigma^+\\
    \text{co}\,\left\{ \mathbf{f}^-(\mathbf{x}), \mathbf{f}^+(\mathbf{x})\right\} & \mathbf{x} \in \Sigma \\
    \mathbf{f}^-(\mathbf{x}), & \mathbf{x} \in \Sigma^-
    \end{cases},
\end{equation}
where
\begin{equation}
    \begin{array}{l}
    \displaystyle
        \text{co}\,\{ \mathbf{f}^-, \mathbf{f}^+ \} = \left\{ \frac{1 + \lambda}{2} \mathbf{f}^+ + \frac{1-\lambda}{2}\mathbf{f}^-\, , \lambda \in (-1,1)\right\},\\
    \lambda = \text{sign}\,(\sigma(\mathbf{x})).
    \end{array}
\end{equation}
Filippov's inclusion allows us to derive a model for the dynamics within the switching surface $\Sigma$. Indeed, exploiting the invariance relationship $\dot{\sigma}(\mathbf{x}) = 0$, we obtain 
\begin{equation}
    \lambda^\Sigma = \frac{(\mathbf{f}^- + \mathbf{f}^+)\cdot \nabla\sigma}{(\mathbf{f}^- - \mathbf{f}^+)\cdot \nabla\sigma},
\end{equation}
\begin{equation}\label{evolution_attractive_sliding}
    \dot{\mathbf{x}} = \mathbf{f}^\Sigma = \frac{(\mathbf{f}^- \cdot \nabla \sigma)\mathbf{f}^+ - (\mathbf{f}^+ \cdot \nabla \sigma)\mathbf{f}^-}{(\mathbf{f}^- - \mathbf{f}^+)\cdot \nabla \sigma}.
\end{equation}

\subsection{Model reduction strategy}\label{supp_SSM}

We recall some fundamental concepts about primary smooth SSMs and their relevance for model reduction (\citet{NNM}). Restricting ourselves to a single region in which governing equations are smooth, we write
\begin{equation}\label{standard_eq}
\begin{aligned}
        \dot{\mathbf{x}} = \mathbf{f}(\mathbf{x}, \mathbf{\Omega}t; \epsilon) &= \mathbf{A}\mathbf{x} + \mathbf{f}_0 + \epsilon\mathbf{f}_1(\mathbf{x},\mathbf{\Omega} t), \\&\quad \mathbf{f}_0(\mathbf{x}) = \mathcal{O}(|\mathbf{x}|^2), \quad 0\leq\epsilon\ll 1,
\end{aligned}
\end{equation}
where $\mathbf{A} \in \mathbb{R}^n$ is the constant matrix of the linear part, $\mathbf{f}_0: \mathbb{R}^n \to \mathbb{R}^n$ and $\mathbf{f}_1: \mathbb{R}^n\times\mathbb{T}^l \to \mathbb{R}^n$, where $\mathbb{T}^\textit{l} = S^1\times ... \times S^\textit{l}$ is the \textit{l}-dimensional torus and with $\mathbf{f}_0, \mathbf{f}_1$ being class $\mathcal{C}^r$ functions. The degree of smoothness of the right-hand side is $r \in \mathbb{N}^+ \cup \{\infty, a\}$, where $a$ refers to analiticity.
The eigenvalues $\lambda_j = \alpha_j + i\omega_j \in \mathbb{C}$ of $\mathbf{A}$ are counted with their multiplicites and addressed in descending order according to their real part
\begin{equation}   
    \text{Re}\lambda_n \le  \text{Re}\lambda_{n-1} \le ... \le \text{Re}\lambda_1 < 0.
\end{equation}
The real and imaginary parts of the eigenvectors or generalized eigenvectors of $\mathbf{A}$ relative to the j-th eigenvalue gives rise to the j-th real modal subspace $E_j \subset \mathbb{R}^n$. 
The direct sum of modal subspaces defines a spectral subspace
\begin{equation}
    E_{j_1,...,j_q} = E_{j_1} \oplus E_{j_2} \oplus ... \oplus E_{j_q}.
\end{equation}
If we further assume that $\text{Re}\lambda_j < 0 \,\forall j$, projecting the linearized system onto the nested hierarchy of slow spectral subspaces \begin{equation}
    E^1 \subset E^2 \subset ..., \quad E^k := E_{1,...,k} \, \text{for} \, k = 1,...,n
\end{equation} 
defines a strategy for reducing the order of the linearized dynamics with increasing accuracy as $k$ is increased (Galerkin projection method).  
We want now to study the existence of a nonlinear continuation of such a spectral subspace, in order to be able to reduce the order of the model, but in the presence of nonlinear and also time-dependent terms. 
In particular, let us fix a specific spectral subspace $E = E_{j_1,...,j_q}$. If the non-resonance condition
\begin{equation}\label{nonresonance_condition_nonauto}
    \langle m, \text{Re}\lambda \rangle_E \neq \text{Re}\lambda_l,\quad \lambda_l \notin \text{Spect}(A|_E),\quad 2\leq |m| \leq \Sigma(E)
\end{equation}
holds, then infinitely many nonlinear continuations of $E$ exist, for $\epsilon$ small enough (\citet{haller2023_fractional}). 
Here $|m| := \sum_{j = 1}^q m_j$.
In \eqref{nonresonance_condition_nonauto}, $\Sigma(E)$ is the absolute spectral quotient defined as
\begin{equation}
    \Sigma(E) = \text{Int}\left[\frac{\displaystyle\min_{\lambda\in \text{Spect}(A) - \text{Spect}(A|_E)} \text{Re}\lambda }{\displaystyle\max_{\lambda\in\text{Spect}(A|_E)} \text{Re}\lambda }\right].
\end{equation}
Among all $\mathcal{C}^{\Sigma(E)+1}$ invariant manifolds there is a unique smoothest one, the primary SSM $\mathcal{W}(E)$, which can therefore be approximated more accurately than the other infinitely many nonlinear continuations of $E$.
In the case of autonomous systems ($\epsilon = 0$), the non-resonance condition can be relaxed to
\begin{equation}
    \langle m, \lambda \rangle_E \neq \lambda_l,\quad \lambda_l \notin \text{Spect}(A|_E),\quad 2\leq |m| \leq \sigma(E),
\end{equation}
where $\sigma(E)$ is the relative spectral quotient, defined as
\begin{equation}\label{relative_spectral_quotient}
    \sigma(E) = \text{Int}\left[\frac{\displaystyle\min_{\lambda\in \text{Spect}(A)} \text{Re}\lambda }{\displaystyle\max_{\lambda\in\text{Spect}(A|_E)} \text{Re}\lambda }\right].
\end{equation}

\subsection{Equation-driven model reduction for the Shaw-Pierre model}\label{supp_SP}

Once we obtain the governing equations of the system centered at the two fixed points $\mathbf{x}_0^\pm$ (eq. \ref{eq:xi_coordinates}), we perform a linear change of coordinates,
\begin{equation}\label{parametrization_matrix}
    \boldsymbol{\xi}^\pm = \mathbf{V}\boldsymbol{\eta}^\pm,
\end{equation}
where
\begin{equation*}
    \boldsymbol{\eta}^\pm = (\mathbf{y}^\pm, \mathbf{z}^\pm) \in E_1 \times E_2
\end{equation*}
and $\mathbf{V}$ is the matrix composed of the eigenvectors of $\mathbf{\tilde{A}_0}$ as columns. For ease of notation, the superscript $\pm$ is dropped in the following equations. The equations of motion now read
\begin{equation}\label{master_dynamics}
\begin{aligned}
&
    \dot{\mathbf{y}} = \mathbf{\tilde{A}_y} \mathbf{y} + \left(\mathbf{r} \,\mathbf{f_{nl}^{III}}\right)_\mathbf{y} \pm \left(\mathbf{r} \,\mathbf{f_{nl}^{II}}\right)_\mathbf{y} = 
    \begin{pmatrix}
    -0.0789 & 1.0342 \\
    -1.0342 & -0.0789
    \end{pmatrix}
    \mathbf{y}  \\& \quad
   + \begin{pmatrix}
    0.9968\\
    -0.0761
    \end{pmatrix}
    q^{III}(\boldsymbol{\eta})
    \pm \begin{pmatrix}
    0.9968\\
    -0.0761
    \end{pmatrix}
    q^{II}(\boldsymbol{\eta})
\end{aligned}
\end{equation}
and
\begin{equation}\label{slave_dynamics}
\begin{aligned} &
    \dot{\mathbf{z}} = \mathbf{\tilde{A}_z} \mathbf{z} + \left(\mathbf{r} \,\mathbf{f_{nl}^{III}}\right)_\mathbf{z} \pm \left(\mathbf{r} \,\mathbf{f_{nl}^{II}}\right)_\mathbf{z} = 
    \begin{pmatrix}
    -0.3711 & 1.6987 \\
    -1.6987 & -0.3711
    \end{pmatrix}
    \mathbf{z} \\&\quad+ 
    \begin{pmatrix}
    -0.8278\\
    0.1808
    \end{pmatrix}
    q^{III}(\boldsymbol{\eta})
    \pm 
    \begin{pmatrix}
    -0.8278\\
    0.1808
    \end{pmatrix}
    q^{II}(\boldsymbol{\eta}),
\end{aligned}
\end{equation}
where
\begin{equation*}
    q^{III}(\boldsymbol{\eta}) = -\frac{1}{2} (0.1645 z_1  + 0.3070 z_2  + 0.0441 y_1 -0.4828 y_2)^3
\end{equation*}
and
\begin{equation*}
    q^{II}(\boldsymbol{\eta}) = \frac{3}{2} q_0 (0.1645 z_1  + 0.3070 z_2  + 0.0441 y_1 -0.4828 y_2)^2.
\end{equation*}
For each case, we seek a cubic Taylor expansion as an approximation of the slow primary SSM anchored at $\boldsymbol{\xi}^\pm = 0$ in the form
\begin{equation}\label{parametrization}
    \mathbf{z} = \sum_{|\mathbf{p}| = 2}^3 \mathbf{h}_\mathbf{p}\mathbf{y}^\mathbf{p}, \qquad \mathbf{p} = (p_1, p_2) \in \mathbb{N}^2,
\end{equation}
where
\begin{equation*}
    \mathbf{y}^\mathbf{p} = y_1^{p_1}\,y_2^{p_2}, \qquad \mathbf{h}_\mathbf{p} = \begin{pmatrix}
        h_{1,\mathbf{p}} \\ h_{2,\mathbf{p}}
    \end{pmatrix} \in \mathbb{R}^2 .
\end{equation*}
The expanded formula is reported in equation (\ref{parametrization_expansion}) below:
\begin{equation}\label{parametrization_expansion}
	\begin{cases}
	\begin{aligned}
            &z_1^\pm = h_1^{\pm, (2,0)} y_1^2 + h_1^{\pm, (1,1)} y_1 y_2 + h_1^{\pm, (0,2)} y_2^2 \\&\quad + h_1^{\pm, (3,0)} y_1^3 + h_1^{\pm, (2,1)} y_1^2 y_2 + h_1^{\pm, (1,2)} y_1^1 y_2^2  \\&\quad + h_1^{\pm, (0,3)} y_2^3 ,\\
            &z_2^\pm = h_2^{\pm, (2,0)} y_1^2 + h_2^{\pm, (1,1)} y_2 y_2 + h_2^{\pm, (0,2)} y_2^2 \\&\quad + h_2^{\pm, (3,0)} y_1^3 + h_2^{\pm, (2,1)} y_1^2 y_2 + h_2^{\pm, (1,2)} y_1^1 y_2^2  \\&\quad + h_2^{\pm, (0,3)} y_2^3. 
            \end{aligned}
	\end{cases}
\end{equation}
Also the reduced dynamics can be written in the form of a Taylor expansion:
\begin{equation}\label{reduced_dynamics_expansion}
	\begin{cases}
        \begin{aligned}
             & r_1^\pm =r_1^{\pm, (1,0)} y_1 + r_1^{\pm, (0,1)} y_2 + r_1^{\pm, (2,0)} y_1^2 + r_1^{\pm, (1,1)} y_1 y_2 
             \\&\quad+ r_1^{\pm, (0,2)} y_2^2 +r_1^{\pm, (3,0)} y_1^3 + r_1^{\pm, (2,1)} y_1^2 y_2
              \\&\quad+ r_1^{\pm, (1,2)} y_1^1 y_2^2 + r_1^{\pm, (0,3)} y_2^3 \\
             &r_2^\pm =r_2^{\pm, (1,0)} y_1 + r_2^{\pm, (0,1)} y_2 + r_2^{\pm, (2,0)} y_1^2 + r_2^{\pm, (1,1)} y_2 y_2 
             \\&\quad+ r_2^{\pm, (0,2)} y_2^2 +r_2^{\pm, (3,0)} y_1^3 + r_2^{\pm, (2,1)} y_1^2 y_2 
             \\&\quad+ r_2^{\pm, (1,2)} y_1^1 y_2^2 + r_2^{\pm, (0,3)} y_2^3.
        \end{aligned}   
	\end{cases}
\end{equation}
The invariance equation on the primary SSM reads
\begin{equation}\label{invariance_equation}
\begin{aligned}
&    D_y\mathbf{h}(\mathbf{y})\, A_y \,\mathbf{y} + D_y \mathbf{h}(\mathbf{y})\,\mathbf{r}_y\, \left[ q^{III}\left(\mathbf{y},\mathbf h(\mathbf{y})\right) \pm q^{II}\left(\mathbf{y},\mathbf h(\mathbf{y})\right) \right] \\&\quad = \,A_z \, \mathbf{h}(\mathbf{y}) + \mathbf{r}_z \, \left[q^{III}(\mathbf{y}, \mathbf{h}\left(\mathbf{y} )\right)\pm q^{II}(\mathbf{y}, \mathbf{h}\left(\mathbf{y})\right) \right].
\end{aligned}
\end{equation}
Equating powers of $\mathbf{y}$ leads to a set of linear equations with the parametrization coefficients as unknowns (\ref{order2_inv_eq} and \ref{order3_inv_eq}) that can be solved as long as nonresonance conditions are satisfied.  The values of the coefficients in \eqref{parametrization_expansion} and \eqref{reduced_dynamics_expansion} are reported in Tables \ref{positive_parametrization},  \ref{negative_parametrization}, \ref{positive_reduced_dynamics} and \ref{negative_reduced_dynamics}, for $\delta = 0.1$.\\
\begin{equation}\label{order2_inv_eq}
	\begin{aligned}	
	&\mathcal{O}(2,0): \quad D_y \mathbf{h}^{(2)} A_y \mathbf {y} =  A_z \mathbf{h}^{(2)} \pm \mathbf{r}_z q^{II, (2)} \\ 
	&\begin{cases}
	 &a_1  \left(2  \text{Re}(\lambda_1)  - \text{Re}(\lambda_2)\right)     - a_2  \text{Im}(\lambda_1) - b_1  \text{Im}(\lambda_2)  \\&\quad=  \pm3  \alpha  q_0  r_{z,1}  p_3^2 \\
        &a_2 \left(2  \text{Re}(\lambda_1) - \text{Re}(\lambda_2)\right)     + 2  a_1  \text{Im}(\lambda_1) - 2  a_3  \text{Im}(\lambda_1)  \\&\quad - b_2  \text{Im}(\lambda_2)= \pm 6  \alpha  p_3  p_4  q_0  r_{z,1} \\
        &a_3  \left(2  \text{Re}(\lambda_1) - \text{Re}(\lambda_2)  \right)       + a_2  \text{Im}(\lambda_1) - b_3  \text{Im}(\lambda_2) \\&\quad=  \pm 3  \alpha  q_0  r_{z,1}  p_4^2 \\
        &b_1  \left( 2  \text{Re}(\lambda_1) - \text{Re}(\lambda_2)\right )     + a_1  \text{Im}(\lambda_2) - b_2  \text{Im}(\lambda_1) \\&\quad=  \pm 3  \alpha  q_0  r_{z,2}  p_3^2 \\
        &b_2\left(2  \text{Re}(\lambda_1) - \text{Re}(\lambda_2)\right)     + a_2  \text{Im}(\lambda_2) + 2  b_1  \text{Im}(\lambda_1) \\&\quad - 2  b_3  \text{Im}(\lambda_1)  = \pm 6  \alpha  p_3  p_4  q_0  r_{z,2} \\
        &b_3  \left(2  \text{Re}(\lambda_1) - \text{Re}(\lambda_2)  \right)     + a_3  \text{Im}(\lambda_2) + b_2  \text{Im}(\lambda_1) \\&\quad= \pm 3  \alpha  q_0  r_{z,2}  p_4^2 \\
	\end{cases}
	\end{aligned}
\end{equation}

\begin{equation}\label{order3_inv_eq}
	\begin{aligned}
	&\mathcal{O}(3,0): \quad D_y \mathbf {h}^{(3)} A_y \mathbf {y} \pm D_y \mathbf {h}^{(2)} \mathbf {r}_y q^{II,(2)} =  A_z \mathbf {h}^{(3)} \\&\qquad+ \mathbf {r}_z \left( q^{III, (3)} \pm q^{II,(3)} \right)\\
	&\begin{cases}
	& a_4\left(3  \text{Re}(\lambda_1) - \text{Re}(\lambda_2)\right)  \mp r_{z,1}(6  \alpha  q_0  (a_1  p_1 + b_1  p_2)  p_3)  \\
		&\quad - a_5  \text{Im}(\lambda_1)- b_4  \text{Im}(\lambda_2) \pm 3  \alpha  p_3^2  q_0  (2  a_1  r_{y,1} + a_2  r_{y,2}) \\
		&\quad = -r_{z,1}  \alpha  p_3^3\\
		
		& a_5\left(3  \text{Re}(\lambda_1) - \text{Re}(\lambda_2)\right)+3  a_4  \text{Im}(\lambda_1) - 2  a_6  \text{Im}(\lambda_1)
		\\&\quad  - b_5  \text{Im}(\lambda_2) \mp r_{z,1}(3  \alpha  q_0  (2  p_4  (a_1  p_1 + b_1  p_2) 
		\\&\quad + 2  p_3  (a_2  p_1 + b_2  p_2))) \pm 3  \alpha  p_3^2  q_0  (a_2  r_{y,1} + 2  a_3  r_{y,2}) 
		\\&\quad \pm 6  \alpha  p_3  p_4  q_0  (2  a_1  r_{y,1} + a_2  r_{y,2}) = -r_{z,1}  3  \alpha  p_3^2  p_4 \\
		
		&a_6\left(3  \text{Re}(\lambda_1) - \text{Re}(\lambda_2)\right)+ 2  a_5  \text{Im}(\lambda_1) - 3  a_7  \text{Im}(\lambda_1) 
		\\&\quad - b_6  \text{Im}(\lambda_2) \mp r_{z,1} (3  \alpha  q_0  (2  p_4  (a_2  p_1 + b_2  p_2) 
		\\&\quad + 2  p_3  (a_3  p_1 + b_3  p_2))) \pm 3  \alpha  p_4^2  q_0  (2  a_1  r_{y,1} + a_2  r_{y,2}) 
		 \\&\quad \pm 6  \alpha  p_3  p_4  q_0  (a_2  r_{y,1}+ 2  a_3  r_{y,2}) = -r_{z,1}  3  \alpha  p_3  p_4^2 \\
		 
		 &a_7\left(3  \text{Re}(\lambda_1) - \text{Re}(\lambda_2)\right) \mp r_{z,1}(6  \alpha  q_0  (a_3  p_1 + b_3  p_2)  p_4) 
		 \\&\quad + a_6  \text{Im}(\lambda_1) - b_7  \text{Im}(\lambda_2) \pm 3  \alpha  p_4^2  q_0  (a_2  r_{y,1} + 2  a_3  r_{y,2}) 
		\\&\quad = -r_{z,1}  \alpha  p_4^3\\
		
       		& b_4\left(3  \text{Re}(\lambda_1) - \text{Re}(\lambda_2)\right)\mp r_{z,2}(6  \alpha  q_0  (a_1  p_1 + b_1  p_2)  p_3) 
		\\&\quad + a_4  \text{Im}(\lambda_2) - b_5  \text{Im}(\lambda_1) \pm 3  \alpha  p_3^2  q_0  (2  b_1  r_{y,1} + b_2  r_{y,2}) 
		\\&\quad = -r_{z,2}  \alpha  p_3^3\\
		
		&b_5\left(3  \text{Re}(\lambda_1) - \text{Re}(\lambda_2)\right)+ a_5  \text{Im}(\lambda_2) + 3  b_4  \text{Im}(\lambda_1)  
		\\&\quad - 2  b_6  \text{Im}(\lambda_1)  \mp r_{z,2}(3  \alpha  q_0  (2  p_4  (a_1  p_1 + b_1  p_2) 
		\\&\quad + 2  p_3  (a_2  p_1 + b_2  p_2))) \pm 3  \alpha  p_3^2  q_0  (b_2  r_{y,1} + 2  b_3  r_{y,2}) 
		 \\&\quad \pm 6  \alpha  p_3  p_4  q_0  (2  b_1  r_{y,1} + b_2  r_{y,2}) = -r_{z,2}  3  \alpha  p_3^2  p_4 \\
		 
		 & b_6\left(3  \text{Re}(\lambda_1) - \text{Re}(\lambda_2)\right) + a_6  \text{Im}(\lambda_2) + 2  b_5  \text{Im}(\lambda_1) 
		 \\&\quad - 3  b_7  \text{Im}(\lambda_1)  \mp r_{z,2}(3  \alpha  q_0  (2  p_4  (a_2  p_1 + b_2  p_2) 
		 \\&\quad + 2  p_3  (a_3  p_1 + b_3  p_2)))\pm 3  \alpha  p_4^2  q_0  (2  b_1  r_{y,1} + b_2  r_{y,2}) 
		 \\& \pm 6  \alpha  p_3  p_4  q_0  (b_2  r_{y,1} \quad+ 2  b_3  r_{y,2}) = - r_{z,2}  3  \alpha  p_3  p_4^2 \\

		&b_7\left(3  \text{Re}(\lambda_1) - \text{Re}(\lambda_2)\right)\mp r_{z,2}(6  \alpha  q_0  (a_3  p_1 + b_3  p_2)  p_4) 
		\\&\quad + a_7  \text{Im}(\lambda_2) + b_6  \text{Im}(\lambda_1) \pm 3  \alpha  p_4^2  q_0  (b_2  r_{y,1} + 2  b_3  r_{y,2})     
		\\&\quad  = -r_{z,2}  \alpha  p_4^3
		
	\end{cases}	
	\end{aligned}
\end{equation}

\begin{table*}[]
    \caption*{}
    \centering 
    \begin{tabular}{|p{1em}| c | c | c | c | c | c | c |}
    \hline
     & (2,0)& (1,1) & (0,2) & (3,0) & (2,1) & (1,2) & (0,3) \\
    \hline \hline
    $h_1^+$ & 8.2 $10^{-3}$ & -2.4 $10^{-2}$ & -7.3 $10^{-3}$ &  2.7 $10^{-2}$ & -1.5 $10^{-3}$ & 2.3 $10^{-3}$ & -1 $10^{-3}$\\
    \hline
    $h_2^+$ & 1.5 $10^{-2}$  &   1.7 $10^{-2}$  & -2.8 $10^{-3}$ &    3.4 $10^{-3}$  &  4.6 $10^{-2}$  &  7.2 $10^{-3}$  &   3.2$10^{-2}$\\
    \hline
    \end{tabular}
    \\[10pt]
    \caption{Coefficients of the graph-style parametrization for the positive case.}
    \label{positive_parametrization}
\end{table*}

\begin{table*}[]
    \caption*{}
    \centering 
    \begin{tabular}{|p{1em}| c | c | c | c | c | c | c |}
    \hline
     & (2,0)& (1,1) & (0,2) & (3,0) & (2,1) & (1,2) & (0,3) \\
    \hline \hline
    $h_1^+$ & -8.2 $10^{-3}$ & 2.4 $10^{-2}$ & 7.3 $10^{-3}$ &  2.7 $10^{-2}$ & -1.5 $10^{-3}$ & 2.3 $10^{-3}$ & -1 $10^{-3}$\\
    \hline
    $h_2^+$ & -1.5 $10^{-2}$  &   -1.7 $10^{-2}$  & 2.8 $10^{-3}$ &    3.4 $10^{-3}$  &  4.6 $10^{-2}$  &  7.2 $10^{-3}$  &   3.2 $10^{-2}$\\
    \hline
    \end{tabular}
    \\[10pt]
    \caption{Coefficients of the graph-style parametrization for the negative case.}
    \label{negative_parametrization}
\end{table*}

\begin{table*}[]
    \caption*{}
    \centering 
    \begin{tabular}{|p{1em}| c | c | c | c | c | c | c | c | c |}
    \hline
     & (1,0)&(0,1)&(2,0)& (1,2) & (0,2) & (3,0) & (2,1) & (1,2) & (0,3) \\
    \hline \hline
    $r_1^+$ & -0.074   &    1.004   &  1.4 $10^{-4}$  &  3.8 $10^{-3}$  &   2.6 $10^{-2}$ & -1.8 $10^{-5}$  & 4.5 $10^{-4}$  &   1.4 $10^{-2}$ &   6.5 $10^{-2}$\\
    \hline
    $r_2^+$ & -1.004 &   -0.074 & -3.0 $10^{-5}$ & -8.1 $10^{-4}$   & -5.5 $10^{-3}$ &  3.9 $10^{-6}$ & -9.7 $10^{-5}$ &   -3.1 $10^{-3}$  &  -1.4 $10^{-2}$\\
    \hline
    \end{tabular}
    \\[10pt]
    \caption{Coefficients of the reduced dynamics for the positive case.}
    \label{positive_reduced_dynamics}
\end{table*}

\begin{table*}[]
    \caption*{}
    \centering 
    \begin{tabular}{|p{1em}| c | c | c | c | c | c | c | c | c |}
    \hline
     & (1,0)&(0,1)&(2,0)& (1,2) & (0,2) & (3,0) & (2,1) & (1,2) & (0,3) \\
    \hline \hline
    $r_1^+$ & -0.074   &    1.004   &  -1.4 $10^{-4}$  &  -3.8 $10^{-3}$  &   -2.6 $10^{-2}$ & -1.8 $10^{-5}$  & 4.5 $10^{-4}$  &   1.4 $10^{-2}$ &   6.5 $10^{-2}$\\
    \hline
    $r_2^+$ & -1.004 &   -0.074 & 3.0 $10^{-5}$ & 8.1 $10^{-4}$   & 5.5 $10^{-3}$ &  3.9 $10^{-6}$ & -9.7 $10^{-5}$ &   -3.1 $10^{-3}$  &  -1.4 $10^{-2}$\\
    \hline
    \end{tabular}
    \\[10pt]
    \caption{Coefficients of the reduced dynamics for the negative case.}
    \label{negative_reduced_dynamics}
\end{table*}

\begin{figure*}
    \centering
    \includegraphics[scale=0.35]{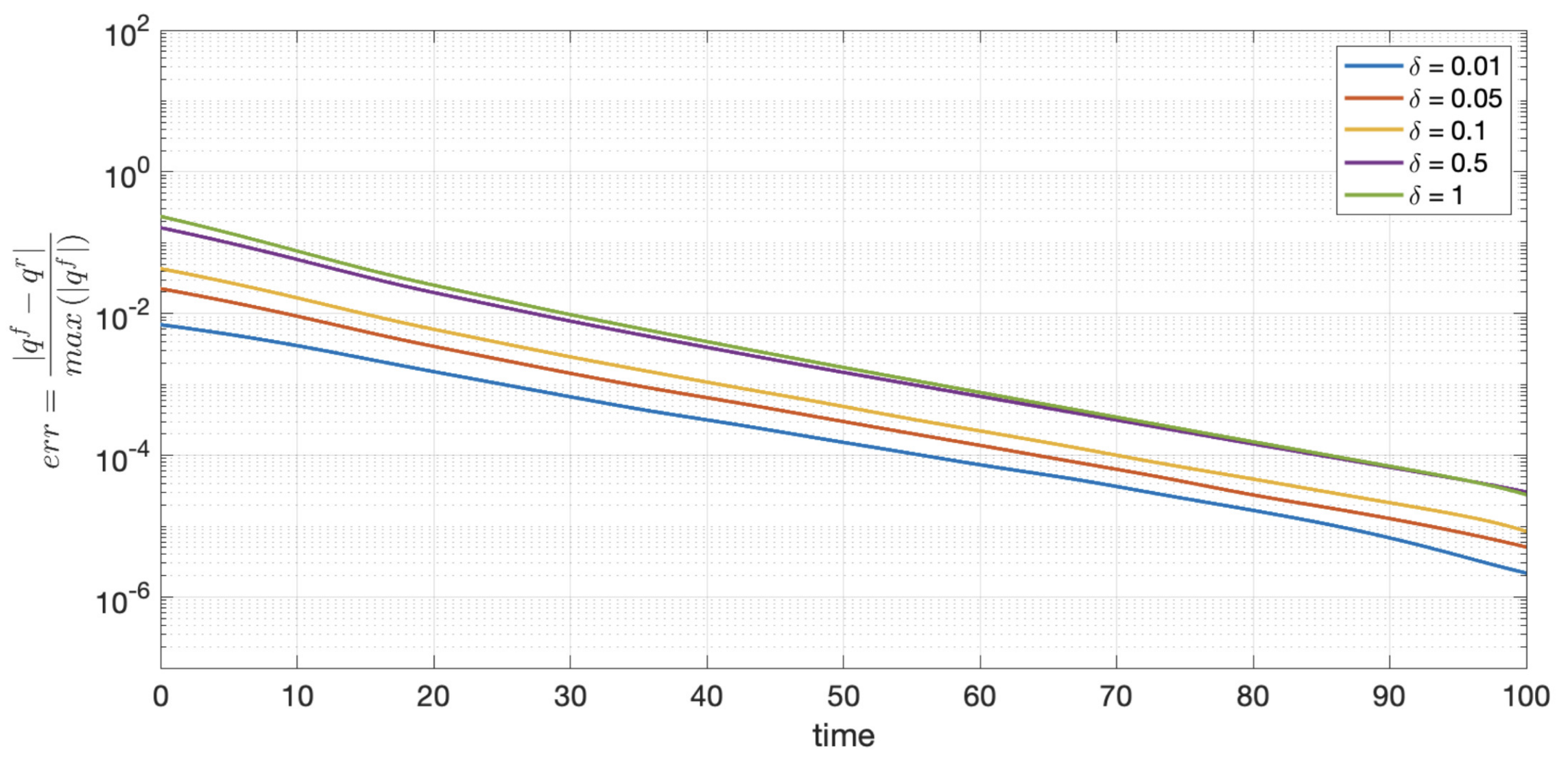}
    \caption{Relative error on the displacement of the first mass for different values of $\delta$.}
    \label{subproblems_error_positive}
\end{figure*}

\subsubsection{Errors}
Considering the positive case, we report in Fig. \ref{subproblems_error_positive} the relative errors between the dynamics of the full system and the reduced one. The initial condition is given slightly outside the SSM: the attracting properties of the SSM are evident as the error reduces significantly with time. Moreover, increasing the value of $\delta$ means intensifying the amplitude of the constant external forcing and therefore the error of the reduced dynamics is shifted upwards. 

The invariance error represents the accuracy of the Taylor expansion in approximating the invariant manifold and the reduced dynamics computed on it. It is quantified as 
\begin{equation}
\begin{aligned}
    &E_{inv} = \frac{1}{N}\sum_{i = 1}^N\frac {\| \mathrm{LHS}(y_i)- \mathrm{RHS}(y_i)\|_2}{\|\mathrm{RHS}(y_i) \|_2}, 
    \\&\quad y_i \in Y_s \subset D(\rho) = \{(y_1, y_2) \in \mathbb{R}^2 : \|y\| = \rho \},
\end{aligned}
\end{equation}
where $\mathrm{LHS}$ and $\mathrm{RHS}$ refer to the left-hand and right-hand sides of equation \eqref{invariance_equation} and $\rho$ denotes the distance from the fixed point in reduced coordinates $(y_1, y_2)$.
\begin{figure}[H]
    \centering
    \subfloat[\label{invariance_error_concept}]{
        \includegraphics[scale=0.28]{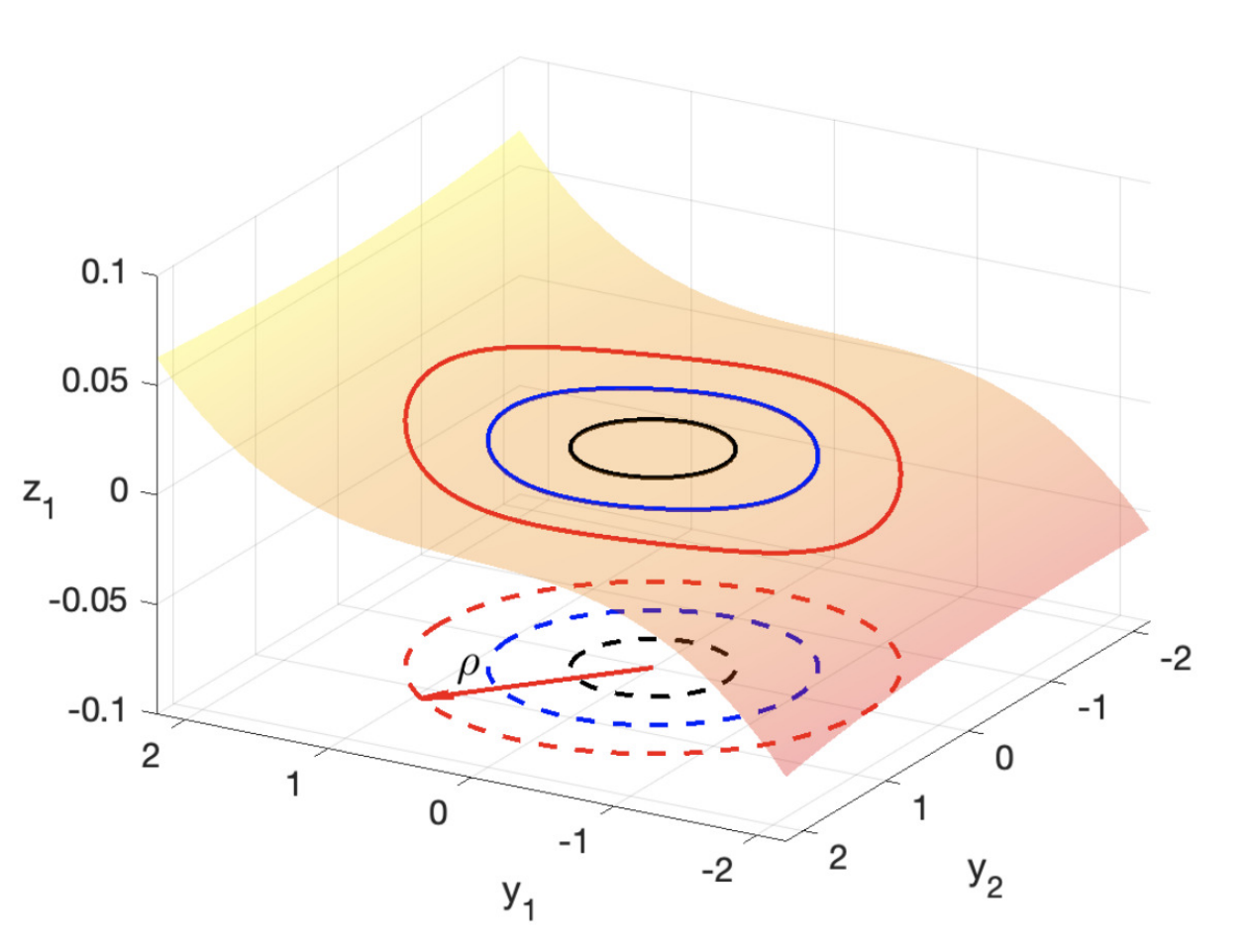}
    }
    \quad
    \subfloat[\label{invariance_error}]{
        \includegraphics[scale=0.15]{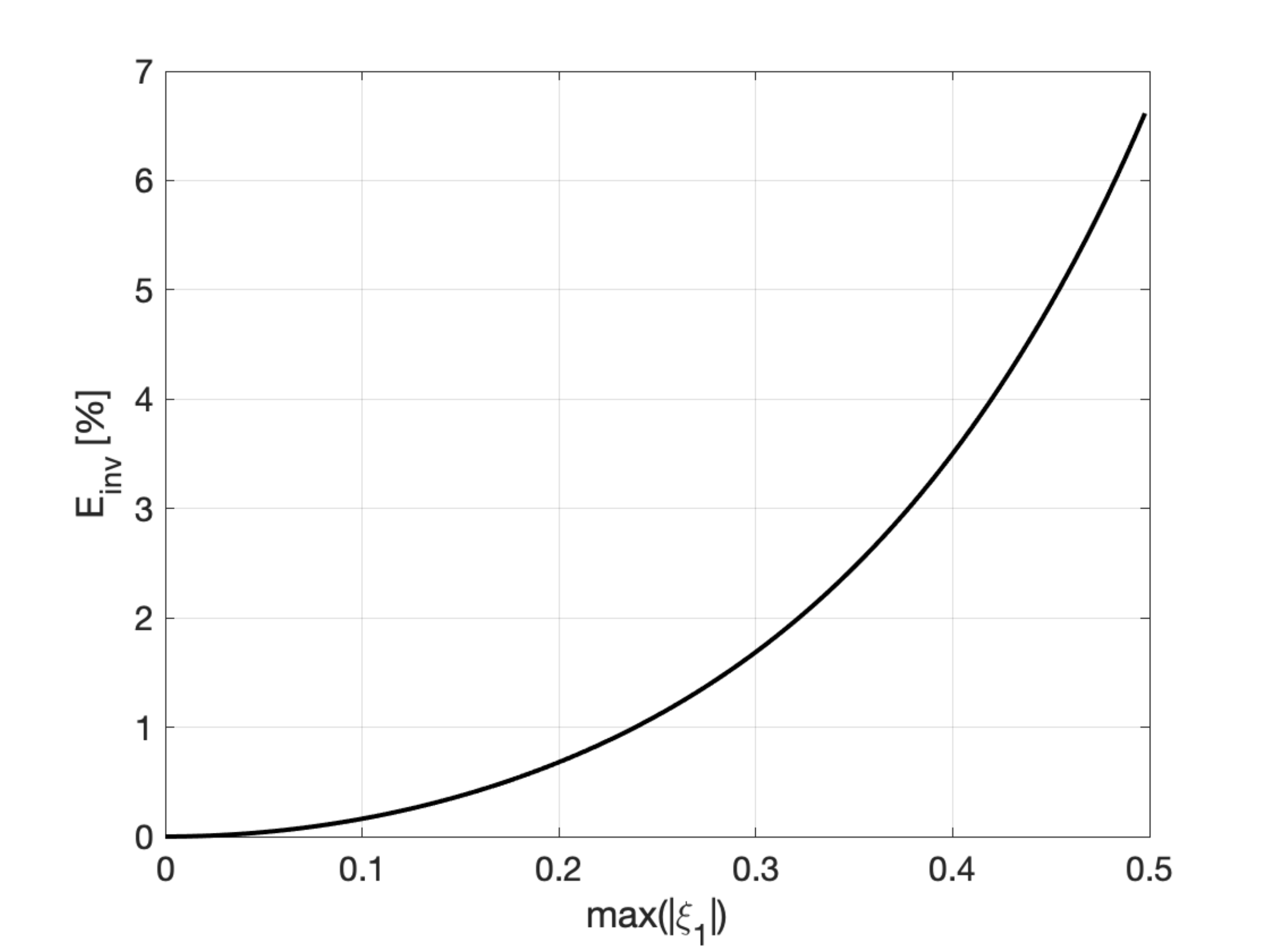}
    }
    \caption{Invariance error for different values of distance from the fixed point $\mathbf{x}_0^+$.}
    \label{friction_separately}
\end{figure}

\subsubsection{Choice of initial conditions across SSMs}\label{choice_ic}
\begin{figure}[H]
    \centering
    \includegraphics[scale=0.25]{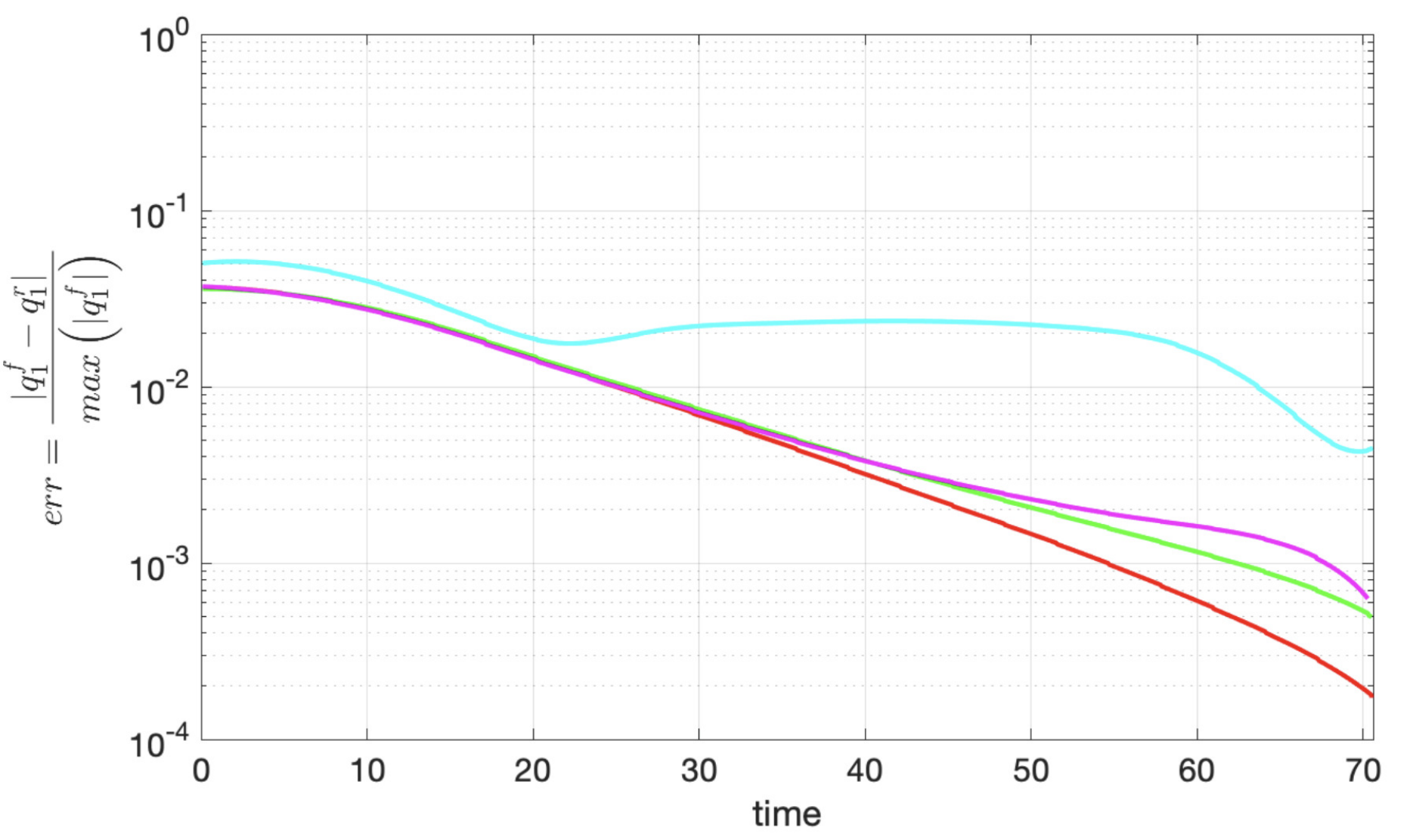}
    \caption{Relative errors between the solution of the reduced model according to the different strategies and the full model, evaluated for the displacement of the first mass and for $\delta = 10^{-3}$. The colors refer to the strategies depicted in Fig. \ref{ic_strategies_all}: projection onto the new SSM (red), minimisation of all physical variables (green), continuity of $q_1$ (purple) and continuity of both $q_1$ and $q_2$ (light blue).}
    \label{ic_errors}
\end{figure}

With reference to Fig. \ref{ic_strategies_all}, the different methods investigate how enforcing continuity of certain physical variables across the switching surface affects the evolution of the reduced trajectory. The most effective strategy in tracking the full dynamics is the one corresponding to the red line in Fig. \ref{ic_errors}, where no continuity constraint is enforced.  

\subsubsection{Approximation of Poincaré map and invariant set. }\label{SP_supp_poincare_map}

\begin{figure}
    \centering
    \includegraphics[scale=0.24]{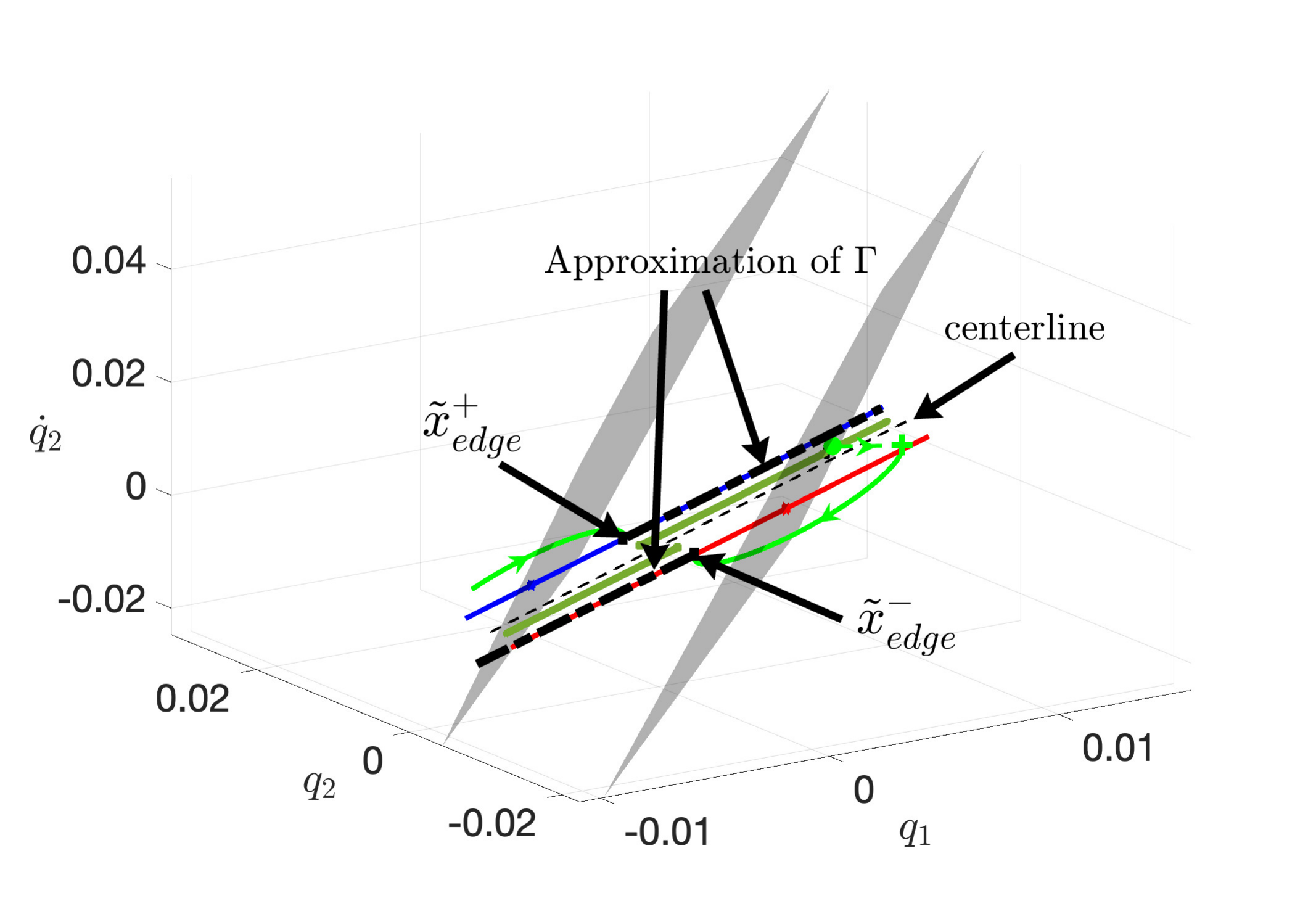}
    \caption{Approximation of the invariant set induced by the Poincaré map (piecewise continuous green line) exploiting information coming from the reduced order model only.}
    \label{edges_reduced_scaled_explained}
\end{figure}

The Poincaré map of the full order model introduced in section \ref{SP_poincare_map} can be approximated by exploiting information from the reduced-order model only. Since the Poincaré map lies close to the intersection between the positive and negative SSMs and the switching surface (blue and red lines in Fig. \ref{edges_reduced_scaled_explained}), the invariant set can be approximated by those curves up to specific points, $\Tilde{x}^\pm_{edge}$, which need to be estimated and play the same role as $x^\pm_{edge}$ for the full order model. We call the average of the blue and red lines the centerline (dashed line in Fig \ref{edges_reduced_scaled_explained}). Recalling that the $x_{edge}^\pm$ are the evolution of the closest points to the sticking regime, we want to approximate such configurations for the reduced order model as well. The strategy is described for the computation of $\Tilde{x}^-_{edge}$, but the same applies to $\Tilde{x}^+_{edge}$. Referring to Fig. \ref{edges_reduced_scaled_explained}, we consider the intersection of the centerline and that of the blue line with the grey surface, representing the configurations on such lines closest to the sticking regime. We take then the average of such intersections to obtain the green circle. Such a point is projected (green cross) onto the negative SSM (not shown in the figure), representing the initial condition sought. The reduced order model is then advected on the negative SSM until the switching surface is hit, identifying $\tilde{x}^-_{edge}$. We note that this procedure requires information coming from the reduced order model only. 

\subsubsection{Non-autonomous problem}\label{supp_SP_na}
\begin{figure}[]
        \centering
        \includegraphics[scale=0.2]{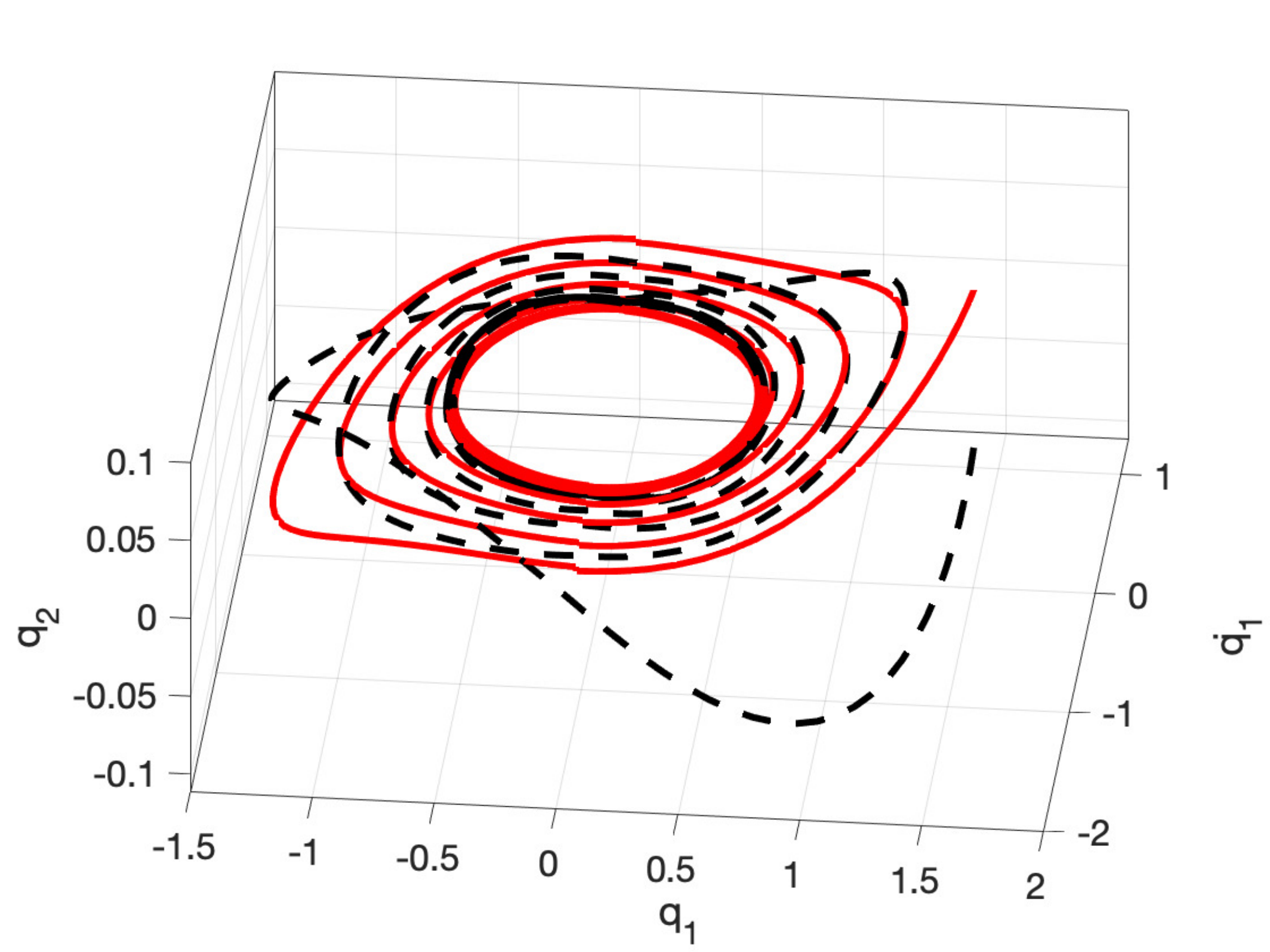}
        \caption{Trajectories of the full (black dashed) and reduced (red) models in the space $(q_1, \dot{q}_1),q_2$, starting from the initial condition $\mathbf{x}_0 = (0, 0, 0.76, 0.76)$, with friction coefficient $\delta = 10^{-2}$ and forcing parameter $\epsilon = 10^{-1}$}
        \label{forced_traj_phys}
\end{figure}

The time-dependent SSMs related to the positive and negative cases are sought in the form of a cubic Taylor expansion
\begin{equation}\label{parametrization_nonauto}
    \mathbf{z} = \sum_{|\mathbf{p}| = 2}^3 \mathbf{h}_\mathbf{p}\mathbf{y}^\mathbf{p} + \epsilon\mathbf{h}_\epsilon (\Omega t), \qquad \mathbf{p} = (p_1, p_2) \in \mathbb{N}^2.
\end{equation}
The invariance equation is now slightly modified to read
\begin{equation}\label{invariance_equation_nonautonomous}
    \begin{aligned}
         &D_y\mathbf{h}(\mathbf{y})\, A_y \,\mathbf{y} + D_y \mathbf{h}(\mathbf{y})\,\mathbf{r}_y\, \left[ q^{III}\left(\mathbf{y},\mathbf h(\mathbf{y})\right) \pm q^{II}\left(\mathbf{y},\mathbf h(\mathbf{y})\right) \right] 
         \\&\quad+ \epsilon\left(V^{-1}\mathbf{f}_\epsilon \right) _y + \epsilon D_t \mathbf{h}_\epsilon = \,A_z \, \mathbf{h}(\mathbf{y}) + \mathbf{r}_z \, \Big[q^{III}(\mathbf{y}, \mathbf{h}\left(\mathbf{y} )\right)  
         \\&\quad \pm q^{II}(\mathbf{y}, \mathbf{h}\left(\mathbf{y})\right) \Big] + \epsilon\left(V^{-1}\mathbf{f}_\epsilon \right)_z .
    \end{aligned}
\end{equation}
Collecting the terms of order $\epsilon$, we obtain a set of ordinary differential equations, which can be solved by a Fourier representation of the periodic forcing and time dependent terms in the parametrization, whereby the forcing term becomes
\begin{equation}\label{fourier_series}
    \mathbf{h}_\epsilon = \sum_{n} \mathbf{\hat{h}}_{\epsilon,n} e^{i n \Omega t} =  \mathbf{\hat{h}}_{-1}e^{-i\Omega t} +  \mathbf{\hat{h}}_1 e^{i \Omega t}.
\end{equation}
The reduced order model well reproduces the dynamics of the system, as long as the forcing amplitude $\epsilon$ remains small and the assumption of small friction coefficient holds. Figure \ref{forced_traj_phys}  shows how the full trajectory rapidly approaches the reduced one, as the initial condition lies outside the primary SSM. They both eventually land on the attracting limit cycle induced by forcing.

\section{Parametrization of invariant manifolds}\label{alternative_param}
In this section, we recall some concepts for the construction of reduced order models with invariant manifolds in a data-driven setting. In particular, referring to \cite{mattia_toappear}, we show that one can choose arbitrary reduced coordinates, as long as they describe the invariant manifold as a graph. This justifies the choice of physical coordinates as reduced coordinates in the von Kárman beam example with friction on moving ground, as it allows us to check the sticking condition in a convenient way .\\
In the setting of eq. \eqref{standard_eq}, we consider a generic spectral subspace
\begin{equation}
    E^{2q} = E_{j_1} \oplus E_{j_2} \oplus ... \oplus E_{j_q},
\end{equation}
and denote as $\text{Spectr}(\mathbf{A}|_{E^{2q}})$ the set of the eigenvalues related to $E^{2q}$. This spectral subspace can be represented by the span of the eigenvectors of $\mathbf{A}|_{E^{2q}}$. The eigenvectors are collected as columns of the matrix $\mathbf{V}_{E^{2q}} \in \mathbb{C}^{2n\times2q}$ satisfying $\mathbf{A}\mathbf{V}_{E^{2q}} = \mathbf{V}_{E^{2q}}\mathbf{R}_{E^{2q}}$, where $\mathbf{R}_{E^{2q}}$ is a diagonal matrix containing $\text{Spectr}(\mathbf{A}|_{E^{2q}})$ as diagonal. Matrix $\mathbf{V}_{E^{2q}}$ defines also its adjoint  $\mathbf{W}_{E^{2q}} \in \mathbb{C}^{2q\times2n}$, such that $\mathbf{W}_{E^{2q}}\mathbf{A} = \mathbf{R}_{E^{2q}}\mathbf{W}_{E^{2q}}$ and normalized according to $\mathbf{W}_{E^{2q}}\mathbf{V}_{E^{2q}} = \mathbf{I}$. Assuming that the spectral subspace $E^{2q}$ satisfies the nonresonance conditions \eqref{nonresonance_condition_nonauto}, then the phase space is characterized by SSMs tangent to $E^{2q}$, among which the primary one is the smoothest and addressed as $\mathcal{W}(E^{2q})$. If we want to study the dynamics on the manifold $\mathcal{W}(E^{2q})$, we need the coordinate chart $\mathbf{y} = \mathbf{w}(\mathbf{x}, \mathbf{\Omega t}; \epsilon)$ and its inverse, i.e. the parametrization of the manifold $\mathbf{x} = \mathbf{v}(\mathbf{y}, \mathbf{\Omega} t; \epsilon)$
\begin{equation}\label{invertibility}
    \mathbf{y} = \mathbf{w(\mathbf{v}(\mathbf{y},\mathbf{\Omega}t;\epsilon), \mathbf{\Omega}t;\epsilon)},
    \quad
    \mathbf{x} = \mathbf{v(\mathbf{w}(\mathbf{x},\mathbf{\Omega}t;\epsilon), \mathbf{\Omega}t;\epsilon)}.
\end{equation}
Hence, the reduced dynamics read
\begin{equation}
    \dot{\mathbf{y}}=\mathbf{r}(\mathbf{y},\mathbf{\Omega}t; \epsilon).
\end{equation}
Since the mappings $\mathbf{v}, \mathbf{w}$ are invariant under the full and the reduced dynamics $\mathbf{f}$ and $\mathbf{r}$, we write
\begin{equation}\label{invariance_equations}
    \begin{cases}
    \begin{aligned}
        &D_\mathbf{y} \mathbf{v}(\mathbf{y},\mathbf{\Omega}t;\epsilon)\,\mathbf{r}(\mathbf{y},\mathbf{\Omega}t; \epsilon) + D_t\mathbf{v}(\mathbf{y},\mathbf{\Omega}t;\epsilon) 
        \\&\quad = \mathbf{f}(\mathbf{v}(\mathbf{y},\mathbf{\Omega}t;\epsilon),\mathbf{\Omega}t;\epsilon),\\
        &D_\mathbf{x} \mathbf{w}(\mathbf{x},\mathbf{\Omega}t;\epsilon)\,\mathbf{f}(\mathbf{x},\mathbf{\Omega}t; \epsilon) + D_t\mathbf{w}(\mathbf{x},\mathbf{\Omega}t;\epsilon) 
        \\&\quad= \mathbf{r}(\mathbf{w}(\mathbf{x},\mathbf{\Omega}t;\epsilon),\mathbf{\Omega}t;\epsilon).
        \end{aligned}
    \end{cases}
\end{equation}
Moreover, as the SSM depends smoothly on the parameter $\epsilon$, we seek the expansions for $\mathbf{w}, \mathbf{v}$ and $\mathbf{r}$ in the form
\begin{equation}\label{expansions}
    \begin{cases}
        \mathbf{w}(\mathbf{x},\mathbf{\Omega}t;\epsilon) = \mathbf{W}_0 + \mathbf{w}_{nl}(\mathbf{x}) + \epsilon\mathbf{w}_1(\mathbf{\Omega}t) + \mathcal{O}(\epsilon\|\mathbf{x}\|), \\
        \mathbf{v}(\mathbf{y},\mathbf{\Omega}t;\epsilon) = \mathbf{V}_0 + \mathbf{v}_{nl}(\mathbf{y}) + \epsilon\mathbf{v}_1(\mathbf{\Omega}t) + \mathcal{O}(\epsilon\|\mathbf{y}\|),\\
        \mathbf{r}(\mathbf{y},\mathbf{\Omega}t;\epsilon) = \mathbf{R}_0 + \mathbf{r}_{nl}(\mathbf{y}) + \epsilon\mathbf{r}_1(\mathbf{\Omega}t) + \mathcal{O}(\epsilon\|\mathbf{y}\|),\\
    \end{cases}
\end{equation}
with $\mathbf{W}_0 \in \mathbb{C}^{2q\times2n}$, $\mathbf{V}_0 \in \mathbb{C}^{2n\times2q}$ and $\mathbf{R}_0 \in \mathbb{C}^{2q\times2q}$ and such that $\text{range}(\mathbf{V}_0) = E^{2q}$ and $\text{Spectr}(\mathbf{R}_0) = \text{Spectr}(\mathbf{A}|_{E^{2q}})$. Substituting the expansions \eqref{expansions} into \eqref{invertibility} and taking the linear contributions for $\epsilon = 0$, we obtain
\begin{equation}\label{identity_relationship}
    \mathbf{W}_0\mathbf{V}_0 = \mathbf{I},
\end{equation}
while considering equations \eqref{invariance_equations} gives:
\begin{equation}\label{eq_from_invariance}
    \mathbf{W}_0\mathbf{A}\mathbf{V}_0 = \mathbf{R}_0.
\end{equation}
Since the spectra of  $\mathbf{R}_0$ and that of $\mathbf{A}$ are the same, 
\begin{equation}
	\mathbf{R}_0 = \mathbf{P}\mathbf{R}_{E^{2q}}\mathbf{P}^{-1},
\end{equation}
holds, which, coupled with equation \eqref{eq_from_invariance}, yields:
\begin{equation}\label{P_relationship}
    \mathbf{W}_0 = \mathbf{P}\mathbf{W}_{E^{2q}}\quad \text{and} \quad \mathbf{V}_0 = \mathbf{V}_{E^{2q}} \mathbf{P}^{-1}.
\end{equation}
In principle, the coordinate chart can be arbitrarily chosen as a projection to the reduced coordinates defined by $\mathbf{W}_0$, as long as it is able to describe the manifold as a graph. In other words, once we define $\mathbf{W}_0$ and we know the eigenvectors spanning the tangent plane to the manifold $V_{E^{2q}}$, then from \eqref{identity_relationship} and the second relation of \eqref{P_relationship} we compute the matrix $\mathbf{P}$ as
\begin{equation}\label{computation_P}
    \mathbf{P} = \mathbf{W}_0\mathbf{V}_{E^{2q}}.
\end{equation}
This means that the vectors spanning the plane where the new reduced coordinates live (columns of $\mathbf{V}_0$) define an alternative coordinate chart of the same manifold, as long as the matrix $\mathbf{P}$ is nonsingular. When we add forcing ($\epsilon \neq 0$), we need to solve 

\begin{equation}\label{}
    \begin{aligned}
         &\mathbf{r}_1(\boldsymbol{\Omega}t)  = \mathbf{W}_0\mathbf{A}\mathbf{v}_1(\boldsymbol{\Omega}t) + \mathbf{W}_0\mathbf{f}_1\left(\mathbf{0}, \boldsymbol{\Omega}t; 0 \right),\\[\medskipamount]
         &\mathbf{v}_1(\boldsymbol{\Omega}t) = \left(\mathbf{I} - \mathbf{V}_0\mathbf{W}_0 \right)\mathbf{A}\mathbf{v}_1(\boldsymbol{\Omega} t) \\[\smallskipamount]& \qquad+ \left(\mathbf{I} - \mathbf{V}_0 \mathbf{W}_0 \right) \mathbf{f}_1\left(\mathbf{0},\boldsymbol{\Omega}t; 0 \right).
    \end{aligned}
\end{equation}
which takes into account additional non-modal contributions to the forcing in reduced coordinates, as explained in \citet{mattia_toappear}.


\newpage
\bibliographystyle{elsarticle-harv} 

\bibliography{bibliography}

\begin{thebibliography}{35}
\expandafter\ifx\csname natexlab\endcsname\relax\def\natexlab#1{#1}\fi
\providecommand{\url}[1]{\texttt{#1}}
\providecommand{\href}[2]{#2}
\providecommand{\path}[1]{#1}
\providecommand{\DOIprefix}{doi:}
\providecommand{\ArXivprefix}{arXiv:}
\providecommand{\URLprefix}{URL: }
\providecommand{\Pubmedprefix}{pmid:}
\providecommand{\doi}[1]{\href{http://dx.doi.org/#1}{\path{#1}}}
\providecommand{\Pubmed}[1]{\href{pmid:#1}{\path{#1}}}
\providecommand{\bibinfo}[2]{#2}
\ifx\xfnm\relax \def\xfnm[#1]{\unskip,\space#1}\fi
\bibitem[{Ax{\aa}s et~al.(2023)Ax{\aa}s, Cenedese and Haller}]{joar_2023}
\bibinfo{author}{Ax{\aa}s, J.}, \bibinfo{author}{Cenedese, M.},
  \bibinfo{author}{Haller, G.}, \bibinfo{year}{2023}.
\newblock \bibinfo{title}{{Fast data-driven model reduction for nonlinear
  dynamical systems}}.
\newblock \bibinfo{journal}{Nonlinear Dyn.} \bibinfo{volume}{111},
  \bibinfo{pages}{7941--7957}.
\bibitem[{Benner et~al.(2015)Benner, Gugercin and Willcox}]{benner2015}
\bibinfo{author}{Benner, P.}, \bibinfo{author}{Gugercin, S.},
  \bibinfo{author}{Willcox, K.}, \bibinfo{year}{2015}.
\newblock \bibinfo{title}{A survey of projection-based model reduction methods
  for parametric dynamical systems}.
\newblock \bibinfo{journal}{SIAM Review} \bibinfo{volume}{57},
  \bibinfo{pages}{483--531}.
\bibitem[{Brunton et~al.(2020)Brunton, Noack and Koumoutsakos}]{brunton2020}
\bibinfo{author}{Brunton, S.L.}, \bibinfo{author}{Noack, B.R.},
  \bibinfo{author}{Koumoutsakos, P.}, \bibinfo{year}{2020}.
\newblock \bibinfo{title}{Machine learning for fluid mechanics}.
\newblock \bibinfo{journal}{Annu. Rev. Fluid Mech.} \bibinfo{volume}{52},
  \bibinfo{pages}{477--508}.
\bibitem[{Cabr{\'e} et~al.(2003)Cabr{\'e}, Fontich and de~la Llave}]{cabre2003}
\bibinfo{author}{Cabr{\'e}, X.}, \bibinfo{author}{Fontich, E.},
  \bibinfo{author}{de~la Llave, R.}, \bibinfo{year}{2003}.
\newblock \bibinfo{title}{The parameterization method for invariant manifolds
  i: manifolds associated to non-resonant subspaces}.
\newblock \bibinfo{journal}{Indiana Univ. Math. J.} ,
  \bibinfo{pages}{283--328}.
\bibitem[{Cardin et~al.(2013)Cardin, Da~Silva and Teixeira}]{cardin_2013}
\bibinfo{author}{Cardin, P.T.}, \bibinfo{author}{Da~Silva, P.R.},
  \bibinfo{author}{Teixeira, M.A.}, \bibinfo{year}{2013}.
\newblock \bibinfo{title}{{On singularly perturbed Filippov systems}}.
\newblock \bibinfo{journal}{Eur. J. Appl. Math.} \bibinfo{volume}{24},
  \bibinfo{pages}{835--856}.
\bibitem[{Cardin et~al.(2015)Cardin, {de Moraes} and {da Silva}}]{cardin_2015}
\bibinfo{author}{Cardin, P.T.}, \bibinfo{author}{{de Moraes}, J.R.},
  \bibinfo{author}{{da Silva}, P.R.}, \bibinfo{year}{2015}.
\newblock \bibinfo{title}{{Persistence of periodic orbits with sliding or
  sewing by singular perturbation}}.
\newblock \bibinfo{journal}{J. Math. Anal. Appl.} \bibinfo{volume}{423},
  \bibinfo{pages}{1166--1182}.
\bibitem[{Cenedese et~al.(2022)Cenedese, Ax{\aa}s, B{\"a}uerlein, Avila and
  Haller}]{mattia_2022}
\bibinfo{author}{Cenedese, M.}, \bibinfo{author}{Ax{\aa}s, J.},
  \bibinfo{author}{B{\"a}uerlein, B.}, \bibinfo{author}{Avila, K.},
  \bibinfo{author}{Haller, G.}, \bibinfo{year}{2022}.
\newblock \bibinfo{title}{{Data-driven modeling and prediction of
  non-linearizable dynamics via spectral submanifolds}}.
\newblock \bibinfo{journal}{Nat. Commun.} \bibinfo{volume}{13},
  \bibinfo{pages}{872}.
\bibitem[{Cenedese et~al.(2021)Cenedese, Ax{\aa}s, Yang, Eriten and
  Haller}]{mattia_mechanical_systems_2021}
\bibinfo{author}{Cenedese, M.}, \bibinfo{author}{Ax{\aa}s, J.},
  \bibinfo{author}{Yang, H.}, \bibinfo{author}{Eriten, M.},
  \bibinfo{author}{Haller, G.}, \bibinfo{year}{2021}.
\newblock \bibinfo{title}{{Data-driven Nonlinear Model Reduction to Spectral
  Submanifolds in Mechanical Systems}}.
\newblock \bibinfo{journal}{Phil. Trans. R. Soc. A.} \bibinfo{volume}{380}.
\bibitem[{Cenedese et~al.(2023)Cenedese, Jain, Marconi and
  Haller}]{mattia_toappear}
\bibinfo{author}{Cenedese, M.}, \bibinfo{author}{Jain, S.},
  \bibinfo{author}{Marconi, J.}, \bibinfo{author}{Haller, G.},
  \bibinfo{year}{2023}.
\newblock \bibinfo{title}{{Data-Assisted Non-Intrusive Model Reduction for
  Forced Nonlinear Finite Elements Models}}.
\newblock \bibinfo{journal}{arXiv:2311.17865} .
\bibitem[{Dankowicz and Schilder(2013)}]{COCO}
\bibinfo{author}{Dankowicz, H.}, \bibinfo{author}{Schilder, F.},
  \bibinfo{year}{2013}.
\newblock \bibinfo{title}{{Recipes for Continuation}}.
\newblock \bibinfo{publisher}{Society for Industrial and Applied Mathematics},
  \bibinfo{address}{Philadelphia, PA}.
\bibitem[{Filippov(1988)}]{Filippov}
\bibinfo{author}{Filippov, A.F.}, \bibinfo{year}{1988}.
\newblock \bibinfo{title}{{Differential equations with discontinuous righthand
  sides}}. volume~\bibinfo{volume}{18} of \textit{\bibinfo{series}{Mathematics
  and its Applications (Soviet Series)}}.
\newblock \bibinfo{publisher}{Kluwer Academic Publishers Group},
  \bibinfo{address}{Dordrecht}.
\bibitem[{Ghadami and Epureanu(2022)}]{ghadami2022}
\bibinfo{author}{Ghadami, A.}, \bibinfo{author}{Epureanu, B.I.},
  \bibinfo{year}{2022}.
\newblock \bibinfo{title}{Data-driven prediction in dynamical systems: recent
  developments}.
\newblock \bibinfo{journal}{Philos. Trans. Royal Soc. A.}
  \bibinfo{volume}{380}, \bibinfo{pages}{20210213}.
\bibitem[{Haller et~al.(2023)Haller, Kasz{\'a}s, Liu and
  Ax{\aa}s}]{haller2023_fractional}
\bibinfo{author}{Haller, G.}, \bibinfo{author}{Kasz{\'a}s, B.},
  \bibinfo{author}{Liu, A.}, \bibinfo{author}{Ax{\aa}s, J.},
  \bibinfo{year}{2023}.
\newblock \bibinfo{title}{Nonlinear model reduction to fractional and
  mixed-mode spectral submanifolds}.
\newblock \bibinfo{journal}{Chaos} \bibinfo{volume}{33}.
\bibitem[{Haller and Ponsioen(2016)}]{NNM}
\bibinfo{author}{Haller, G.}, \bibinfo{author}{Ponsioen, S.},
  \bibinfo{year}{2016}.
\newblock \bibinfo{title}{{Nonlinear normal modes and spectral submanifolds:
  Existence, uniqueness and use in model reduction}}.
\newblock \bibinfo{journal}{Nonlinear Dyn.} \bibinfo{volume}{86}.
\bibitem[{Haro and de~la Llave(2006)}]{haro2006}
\bibinfo{author}{Haro, A.}, \bibinfo{author}{de~la Llave, R.},
  \bibinfo{year}{2006}.
\newblock \bibinfo{title}{A parameterization method for the computation of
  invariant tori and their whiskers in quasi-periodic maps: rigorous results}.
\newblock \bibinfo{journal}{J. Differ. Equ.} \bibinfo{volume}{228},
  \bibinfo{pages}{530--579}.
\bibitem[{{Hassler Whitney}(1944)}]{whitney_embed}
\bibinfo{author}{{Hassler Whitney}}, \bibinfo{year}{1944}.
\newblock \bibinfo{title}{{The Self-Intersections of a Smooth n-Manifold in
  2n-Space}}.
\newblock \bibinfo{journal}{Ann. Math.} \bibinfo{volume}{45},
  \bibinfo{pages}{220--246}.
\bibitem[{Jain and Haller(2022)}]{shobit_2022}
\bibinfo{author}{Jain, S.}, \bibinfo{author}{Haller, G.}, \bibinfo{year}{2022}.
\newblock \bibinfo{title}{{How to compute invariant manifolds and their reduced
  dynamics in high-dimensional finite element models}}.
\newblock \bibinfo{journal}{Nonlinear Dyn.} \bibinfo{volume}{107},
  \bibinfo{pages}{1--34}.
\bibitem[{Jain et~al.(2023)Jain, Li, Thurnher and Haller}]{shobit_2023}
\bibinfo{author}{Jain, S.}, \bibinfo{author}{Li, M.},
  \bibinfo{author}{Thurnher, T.}, \bibinfo{author}{Haller, G.},
  \bibinfo{year}{2023}.
\newblock \bibinfo{title}{{SSMTool: Computation of invariant manifolds in
  high-dimensional mechanics problems}}.
\bibitem[{Jain et~al.(2017)Jain, Tiso and Haller}]{shobit_2017}
\bibinfo{author}{Jain, S.}, \bibinfo{author}{Tiso, P.},
  \bibinfo{author}{Haller, G.}, \bibinfo{year}{2017}.
\newblock \bibinfo{title}{{Exact Nonlinear Model Reduction for a von Karman
  beam: Slow-Fast Decomposition and Spectral Submanifolds}}.
\newblock \bibinfo{journal}{J. Sound Vib.} \bibinfo{volume}{423}.
\bibitem[{Kasz{\'a}s et~al.(2022)Kasz{\'a}s, Cenedese and Haller}]{balint_2022}
\bibinfo{author}{Kasz{\'a}s, B.}, \bibinfo{author}{Cenedese, M.},
  \bibinfo{author}{Haller, G.}, \bibinfo{year}{2022}.
\newblock \bibinfo{title}{{Dynamics-based machine learning of transitions in
  Couette flow}}.
\newblock \bibinfo{journal}{Phys. Rev. Fluids} \bibinfo{volume}{7}.
\bibitem[{K{\"u}pper(2008)}]{kupper_2008}
\bibinfo{author}{K{\"u}pper, T.}, \bibinfo{year}{2008}.
\newblock \bibinfo{title}{{Invariant cones for non-smooth dynamical systems}}.
\newblock \bibinfo{journal}{Math. Comput. Simul.} \bibinfo{volume}{79},
  \bibinfo{pages}{1396--1408}.
\bibitem[{Leine and Nijmeijer(2004)}]{leine}
\bibinfo{author}{Leine, R.}, \bibinfo{author}{Nijmeijer, H.},
  \bibinfo{year}{2004}.
\newblock \bibinfo{title}{{Dynamics and Bifurcations of Non-Smooth Mechanical
  Systems}}. volume~\bibinfo{volume}{18}.
\bibitem[{Li and Haller(2022)}]{mingwu_2022_II}
\bibinfo{author}{Li, M.}, \bibinfo{author}{Haller, G.}, \bibinfo{year}{2022}.
\newblock \bibinfo{title}{{Nonlinear analysis of forced mechanical systems with
  internal resonance using spectral submanifolds, Part II: Bifurcation and
  quasi-periodic response}}.
\newblock \bibinfo{journal}{Nonlinear Dyn.} \bibinfo{volume}{110},
  \bibinfo{pages}{1045--1080}.
\bibitem[{Li et~al.(2022)Li, Jain and Haller}]{mingwu_2022_I}
\bibinfo{author}{Li, M.}, \bibinfo{author}{Jain, S.}, \bibinfo{author}{Haller,
  G.}, \bibinfo{year}{2022}.
\newblock \bibinfo{title}{{Nonlinear analysis of forced mechanical systemswith
  internal resonance using spectral submanifolds, Part I: Periodic response and
  forced response curve}}.
\newblock \bibinfo{journal}{Nonlinear Dyn.} \bibinfo{volume}{110},
  \bibinfo{pages}{1005--1043}.
\bibitem[{Mikhlin and Avramov(2023)}]{mikhlin_2023}
\bibinfo{author}{Mikhlin, Y.}, \bibinfo{author}{Avramov, K.V.},
  \bibinfo{year}{2023}.
\newblock \bibinfo{title}{Nonlinear normal modes of vibrating mechanical
  systems: 10 years of progress.}
\newblock \bibinfo{journal}{Appl. Mech. Rev.} , \bibinfo{pages}{1--57}.
\bibitem[{Ponsioen et~al.(2020)Ponsioen, Jain and Haller}]{Sten_2020}
\bibinfo{author}{Ponsioen, S.}, \bibinfo{author}{Jain, S.},
  \bibinfo{author}{Haller, G.}, \bibinfo{year}{2020}.
\newblock \bibinfo{title}{{Model reduction to spectral submanifolds and
  forced-response calculation in high-dimensional mechanical systems}}.
\newblock \bibinfo{journal}{J. Sound Vib.} \bibinfo{volume}{488},
  \bibinfo{pages}{115640}.
\bibitem[{Rowley and Dawson(2017)}]{rowley2017}
\bibinfo{author}{Rowley, C.W.}, \bibinfo{author}{Dawson, S.T.},
  \bibinfo{year}{2017}.
\newblock \bibinfo{title}{Model reduction for flow analysis and control}.
\newblock \bibinfo{journal}{Annu. Rev. Fluid Mech.} \bibinfo{volume}{49},
  \bibinfo{pages}{387--417}.
\bibitem[{Shaw and Pierre(1993)}]{shaw_pierre_93}
\bibinfo{author}{Shaw, S.}, \bibinfo{author}{Pierre, C.}, \bibinfo{year}{1993}.
\newblock \bibinfo{title}{Normal modes for non-linear vibratory systems}.
\newblock \bibinfo{journal}{J. Sound Vib.} \bibinfo{volume}{164},
  \bibinfo{pages}{85--124}.
\bibitem[{Shaw and Pierre(1994)}]{shaw_pierre_94}
\bibinfo{author}{Shaw, S.}, \bibinfo{author}{Pierre, C.}, \bibinfo{year}{1994}.
\newblock \bibinfo{title}{Normal modes of vibration for non-linear continuous
  systems}.
\newblock \bibinfo{journal}{J. Sound Vib.} \bibinfo{volume}{169},
  \bibinfo{pages}{319--347}.
\bibitem[{Shaw et~al.(1999)Shaw, Pierre and Pesheck}]{shaw_pierre_99}
\bibinfo{author}{Shaw, S.W.}, \bibinfo{author}{Pierre, C.},
  \bibinfo{author}{Pesheck, E.}, \bibinfo{year}{1999}.
\newblock \bibinfo{title}{Modal analysis-based reduced-order models for
  nonlinear structures : An invariant manifold approach}.
\newblock \bibinfo{journal}{Shock Vib. Dig.} \bibinfo{volume}{31},
  \bibinfo{pages}{3--16}.
\bibitem[{Szalai(2019)}]{szalai_2019}
\bibinfo{author}{Szalai, R.}, \bibinfo{year}{2019}.
\newblock \bibinfo{title}{{Model Reduction of Non-densely Defined
  Piecewise-Smooth Systems in Banach Spaces}}.
\newblock \bibinfo{journal}{J. Nonlinear Sci.} \bibinfo{volume}{29}.
\bibitem[{Taira et~al.(2017)Taira, Brunton, Dawson, Rowley, Colonius, McKeon,
  Schmidt, Gordeyev, Theofilis and Ukeiley}]{taira2017}
\bibinfo{author}{Taira, K.}, \bibinfo{author}{Brunton, S.L.},
  \bibinfo{author}{Dawson, S.T.}, \bibinfo{author}{Rowley, C.W.},
  \bibinfo{author}{Colonius, T.}, \bibinfo{author}{McKeon, B.J.},
  \bibinfo{author}{Schmidt, O.T.}, \bibinfo{author}{Gordeyev, S.},
  \bibinfo{author}{Theofilis, V.}, \bibinfo{author}{Ukeiley, L.S.},
  \bibinfo{year}{2017}.
\newblock \bibinfo{title}{Modal analysis of fluid flows: An overview}.
\newblock \bibinfo{journal}{AIAA J.} \bibinfo{volume}{55},
  \bibinfo{pages}{4013--4041}.
\bibitem[{Touz{\'e} et~al.(2021)Touz{\'e}, Vizzaccaro and Thomas}]{touze2021}
\bibinfo{author}{Touz{\'e}, C.}, \bibinfo{author}{Vizzaccaro, A.},
  \bibinfo{author}{Thomas, O.}, \bibinfo{year}{2021}.
\newblock \bibinfo{title}{Model order reduction methods for geometrically
  nonlinear structures: a review of nonlinear techniques}.
\newblock \bibinfo{journal}{Nonlinear Dyn.} \bibinfo{volume}{105},
  \bibinfo{pages}{1141--1190}.
\bibitem[{Weiss et~al.(2012)Weiss, K{\"u}pper and Hosham}]{weiss_2012}
\bibinfo{author}{Weiss, D.}, \bibinfo{author}{K{\"u}pper, T.},
  \bibinfo{author}{Hosham, H.}, \bibinfo{year}{2012}.
\newblock \bibinfo{title}{{Invariant manifolds for nonsmooth systems}}.
\newblock \bibinfo{journal}{Phys. D: Nonlinear Phenom.} \bibinfo{volume}{241},
  \bibinfo{pages}{1895--1902}.
\bibitem[{Weiss et~al.(2015)Weiss, K{\"u}pper and Hosham}]{weiss_2015}
\bibinfo{author}{Weiss, D.}, \bibinfo{author}{K{\"u}pper, T.},
  \bibinfo{author}{Hosham, H.}, \bibinfo{year}{2015}.
\newblock \bibinfo{title}{{Invariant manifolds for nonsmooth systems with
  sliding mode}}.
\newblock \bibinfo{journal}{Math. Comput. Simul.} \bibinfo{volume}{110},
  \bibinfo{pages}{15--32}.

\end{thebibliography}


%
%
%
\end{document}